\documentclass[preprint,authoryear,3p,times]{elsarticle}

\usepackage{amscd,amssymb,amsmath,amsthm,diagbox}
\usepackage{graphicx,verbatim,lscape,rotating}
\usepackage{color,float,caption,subcaption,geometry}
\usepackage{booktabs,tabularx,multirow,array}
\usepackage{dsfont,here}
\newcommand{\bfzero}{\mathbf {0}}

\newcommand{\bfa}{\mathbf {a}}
\newcommand{\bfb}{\mathbf {b}}
\newcommand{\bfc}{\mathbf {c}}

\newcommand{\bfe}{\mathbf {e}}

\newcommand{\bfr}{\mathbf {r}}

\newcommand{\bfz}{\mathbf {z}}
\newcommand{\bfA}{\mathbf {A}}
\newcommand{\bfB}{\mathbf {B}}
\newcommand{\bfC}{\mathbf {C}}
\newcommand{\bfD}{\mathbf {D}}
\newcommand{\bfE}{\mathbf {E}}
\newcommand{\bfF}{\mathbf {F}}
\newcommand{\bfG}{\mathbf {G}}
\newcommand{\bfH}{\mathbf {H}}
\newcommand{\bfI}{\mathbf {I}}
\newcommand{\bfJ}{\mathbf {J}}
\newcommand{\bfK}{\mathbf {K}}
\newcommand{\bfL}{\mathbf {L}}
\newcommand{\bfM}{\mathbf {M}}
\newcommand{\bfN}{\mathbf {N}}

\newcommand{\bfP}{\mathbf {P}}
\newcommand{\bfR}{\mathbf {R}}

\newcommand{\bfV}{\mathbf {V}}
\newcommand{\bfQ}{\mathbf {Q}}
\newcommand{\bfW}{\mathbf {W}}
\newcommand{\bfX}{\mathbf {X}}
\newcommand{\bfY}{\mathbf {Y}}
\newcommand{\bfZ}{\mathbf {Z}}
\newcommand{\cd}{\stackrel{d}{\to}}
\newcommand{\cp}{\stackrel {\mathbb{P}} {\to}}
\newcommand{\cas}{\stackrel{a.s.} {\to}}

\newcommand{\bfbeta}{\mbox{\boldmath$\beta$}}

\newcommand{\bfepsilon}{\mbox{\boldmath$\epsilon$}}
\newcommand{\bfUpsilon}{\mbox{\boldmath$\Upsilon$}}
\newcommand{\bfgamma}{\mbox{\boldmath$\gamma$}}

\newcommand{\bfmu}{\mbox{\boldmath$\mu$}}
\newcommand{\bfeta}{\mbox{\boldmath$\eta$}}

\newcommand{\bfrho}{\mbox{\boldmath$\rho$}}
\newcommand{\bfsigma}{\mbox{\boldmath$\sigma$}}

\newcommand{\bfxi}{\mbox{\boldmath$\xi$}}
\newcommand{\bfSigma}{\mbox{\boldmath$\Sigma$}}
\newcommand{\bfDelta}{\mathbf{\Delta}}
\newcommand{\bfGamma}{\mbox{\boldmath$\Gamma$}}

\newcommand{\bfXi}{\mbox{\boldmath$\Xi$}}
\newcommand{\bfnabla}{\mbox{\boldmath$\nabla$}}
\newcommand{\bfOmega}{\mathbf{\Omega}}
\newcommand{\bfPhi}{\mathbf{\Phi}}
\newcommand{\bfPsi}{\mathbf{\Psi}}
\newcommand{\bfTheta}{\mathbf{\Theta}}

\newcommand{\tr}{\mbox{tr}}

\newcommand{\cov}{\mbox{cov}}
\newcommand{\diag}{\mbox{diag}}

\newtheorem{theorem}{Theorem}
\newtheorem{proposition}{Proposition}
\newtheorem{corollary}{Corollary}

\newtheorem{remark}{Remark}
\numberwithin{theorem}{section}
\numberwithin{corollary}{section}
\numberwithin{lemma}{section}
\numberwithin{remark}{section}

\def\vec{\mbox{\textnormal{vec}}}

\usepackage{lineno}
\modulolinenumbers[5]

\journal{Journal of Multivariate Analysis}

\begin{document}

\begin{frontmatter}



\title{Diagnostic checking of  periodic vector autoregressive time series models with dependent errors}

\author[uphf,franche]{Yacouba Boubacar Ma\"{i}nassara \corref{cor1}} \ead{Yacouba.BoubacarMainassara@uphf.fr}
\cortext[cor1]{Corresponding author}
\address[uphf]{ Univ. Polytechnique Hauts-de-France, INSA Hauts-de-France, CERAMATHS - Laboratoire de
Matériaux Céramiques et de Mathématiques, F-59313 Valenciennes, France}
\address[franche]{ Universit\'e Bourgogne Franche-Comt\'e, Laboratoire de math\'ematiques de Besançon, UMR CNRS 6623, 16 Route de Gray, 25030 Besançon, France}
\author[gretha,ecreb]{Eugen Ursu}\ead{eugen.ursu@u-bordeaux.fr}
\address[gretha]{Universit\'e de Bordeaux, Bordeaux School of Economics, 16 Avenue L\'eon Duguit, B\^{a}t. H2, 33608 Pessac CEDEX, France}
\address[ecreb]{ECREB, West University of Timisoara, Romania}

\begin{abstract}
In this article, we study the asymptotic behaviour of the residual autocorrelations for periodic vector autoregressive time series models (PVAR henceforth) with uncorrelated but dependent innovations ({\em i.e.}, weak PVAR). We then deduce the asymptotic distribution of the Ljung-Box-McLeod modified Portmanteau statistics for weak PVAR models. In Monte Carlo experiments, we illustrate that the proposed test statistics have reasonable finite sample performance. When the innovations exhibit conditional heteroscedasticity or other forms of dependence, it appears that the standard test statistics (under independent and identically distributed innovations) are generally nonreliable, overrejecting, or underrejecting severely, while the proposed test statistics offer satisfactory levels. An illustrative application on real data  is also proposed.

\end{abstract}


\end{frontmatter}


\section{Introduction}
\noindent Periodic vector time series provide an alternative approach that describes seasonal time series. The correlations in a stationary seasonal model do not vary with season, but there are many applications where model parameters need to vary periodically. These periodic models are non-stationary and they are designed to model time series data that display periodic statistical structure (for vector time series, see, e.g.,~\citet{FP04},~\citet{Lu05} and~\citet{UD09}). Seasonal and periodic time series models are quite different, and~\citet{LB99} present an interesting comparison between them. More recently,~\cite{DHW16} present a review of existing work on time series analysis of periodic ARMA (PARMA) models when periodic parameters are expressed as Fourier series.

It should be noted that VAR models may be seen as special cases of PVAR models. The same is true for seasonal VARMA (Vector ARMA) models, which admit VARMA representations with particular autoregressive and moving average structures. Consequently, several basic properties of VAR models could likely be derived from those established in the class of PVAR models. It is usually recognized that periodic time series models rely on a large number of parameters. For example, a periodic VAR model of order one for bivariate monthly data involves $48$ independent parameters, and finding parsimonious representation is not always an easy task (see, e.g.,~\citet{UD09} or \cite{BMU23}).

In PVAR models, the choice of order $p(\nu)$, $\nu \in \{1,\ldots,s\}$ ($s$ denotes the length of the seasonal period) is particularly important because the number of parameters ($d^2 \sum_{\nu=1}^{s} p(\nu)$, where $d$ is the number of series considered) quickly increases with $p(\nu)$, which entails statistical difficulties.
After identification and estimation of the PVAR processes, the next important step in the modeling process consists of checking if the estimated model fits the data satisfactorily. This adequacy checking step allows to validate or invalidate the choice of the orders $p(\nu)$. Thus, it is important to check the validity of a PVAR model, for a given order $p(\nu)$. The adequacy checking step is often used together with model selection criteria such as the Akaike information criterion (AIC) and the Bayesian information criterion (BIC). These two proceedings complement each other.~\cite{McL1994} focuses on the model selection of PAR models either by examining plots of the periodic partial autocorrelation or by using an information criterion. For a comprehensive treatment of time series model selection, see~\cite{MT98} among others.

Residual autocovariances and autocorrelations are appropriate tools for diagnosing many time series models.~\cite{bp70} have proposed a goodness-of-fit test, the so-called portmanteau test, for univariate ARMA models with independent and identically distributed ({\em i.e.} iid) errors term. A modification of their test has been proposed by~\cite{lb}. A more powerful test, based on the $m$th root of the determinant of the $m$th autocorrelation matrix, was proposed by~\cite{PR02}. Nowadays, it is one of the most popular diagnostic checking tools in ARMA modeling of time series. An overview of diagnostic checking for univariate and multivariate time series models is presented in the monograph of~\citet{Li04}. Inspired by the univariate portmanteau statistics defined in~\cite{bp70} and~\cite{lb}, \cite{C1974} and~\cite{Ho1980} have introduced  the  multivariate versions of the  portmanteau statistics.~\cite{H81} gave several equivalent forms of this statistic. To test simultaneously whether all residual autocorrelations at lags $1,\dots,M$ (for a fixed integer $M$) of a PVAR model are equal to zero for a specified period $\nu$,  the portmanteau test introduced by~\cite{C1974} and~\cite{Ho1980} have been adapted by~\cite{UD09},~\cite{DL13}. The proposed statistics have the same asymptotic distribution. Under the assumption that the noise sequence is iid, the standard test procedure consists of rejecting the null hypothesis of a PVAR model if the statistics tests are larger than a certain quantile of a chi-squared distribution.

The aforementioned test procedures for PVAR models were established under the assumption of independent errors (strong PVAR). Of course, this assumption is not satisfied for multivariate nonlinear processes that admit a weak PVAR representation (the errors are uncorrelated but dependent) such as  the multivariate periodic  generalized autoregressive conditional heteroscedastic (MPGARCH) model (see for instance  \cite{Bi18}) or even when the (possibly dependent) error is subject to unknown multivariate conditional heteroscedasticity. 
From a practical point of view, it appears important to relax the assumption of independent errors as financial data~\citep{Tsay2005} or daily and monthly streamflow processes~\citep{WGVM05} exhibit heavier tails than those coming from the normal distribution. Another argument in favor of considering the weak PVAR models comes from the fact that, in general, temporal aggregation or systematic sampling of a strong periodic model yield a weak periodic model, see~\cite{RS08}. All these examples have important practical meanings and emphasize the need for taking into account an eventual dependence of the errors when testing the adequacy of a PVAR model. For weak time series models, the asymptotic distributions of the statistics are no longer chi-squared distributions but a mixture of chi-squared distributions, weighted by eigenvalues of  the asymptotic covariance matrix of the vector of autocorrelations (see for instance~\cite{frz05,FR07,BM11,BMIA23}).

The main goal of this paper is to derive the asymptotic distribution of the residual autocovariance and autocorrelation matrices in the class of weak PVAR models. It is shown that the obtained asymptotic variance matrices in such class of models are the long-run variances, which may be consistently estimated by a nonparametric kernel estimator, also called  heteroscedastic autocorrelation consistent estimator (see, e.g.,~\cite{A91,newey}) or by a parametric approach (see, for instance,~\cite{B74,haan}). The drawback of needing to estimate additional quantities is that the proposed asymptotic properties require a fairly large sample size to perform properly. Applications of the result include diagnostic checking with portmanteau test. Three different portmanteau tests are considered in a small simulation study: the standard portmanteau test proposed by~\cite{UD09}, the test based on the quantiles of the asymptotic distribution proposed by~\cite{DL13} and a modified test for PVAR models with nonindependent innovations. It is shown that the standard portmanteau test can be quite misleading in the framework of dependent errors. Our results generalize the literature in several directions. They extend previous theorems establishing the asymptotic distributions of residual autocovariance and autocorrelation matrices in weak PVAR models. Our modified portmanteau test statistics generalize the standard portmanteau test of~\cite{UD09} for diagnostic checking of PVAR with independent errors (strong PVAR). Furthermore, our asymptotic results provide multivariate generalizations of theorems obtained by~\cite{BMIA23} when the moving average order is null. When the period is one, we retrieve the results obtained by~\cite{CD08},~\cite{FR07} and~\cite{Ho1980}.

This article is organized as follows. Section~\ref{model} presents the parametrization and Assumptions used in what follows. In Section~\ref{estimation}, we recall the results on the least squares estimators of the weak PVAR model obtained by~\cite{BMU23}.
Section~\ref{result} presents our main results. Then we derive the limiting distribution of the residual  autocovariances and autocorrelations matrices in the framework of  PVAR models with dependent errors under the assumption of unconstrained and also of linear constraints  on the parameters of a given season in Section~\ref{result1}. In Section~\ref{result2}, it is shown how the standard  portmanteau test must be adapted in the case of PVAR models with nonindependent innovations. From the results of Section~\ref{result1}, we propose modified test statistics for checking the adequacy of a weak PVAR model. In Section~\ref{simul}, some simulation results are reported, and an illustrative application on real data is presented in Section~\ref{real}. Finally, Section~\ref{conclusion} offers some concluding remarks. The proofs of the main results are collected in the appendix.

\section{Weak periodic vector autoregressive time series models and Assumptions}\label{model}
\noindent In this section, we present principal results on the least squares estimators in the unconstrained and constrained cases of the PVAR models with dependent errors.

Let $\bfY = \{ \bfY_t, t \in \mathbb{Z} \}$ be a stochastic process, where $$\bfY_t = (Y_t(1),\ldots,Y_t(d))^{\top}$$ represents a random vector of dimension $d$.
The process $\bfY$ is a PVAR process of order $p(\nu)$, $\nu \in \{1,\ldots,s\}$ ($s$ is a predetermined value), if there exist
$d\times d$ matrices $\bfPhi_k(\nu) = \left(  \Phi_{k,ij}(\nu) \right)_{i,j=1,\ldots,d}$, $k=1,\ldots,p(\nu)$ such that
\begin{equation}
\label{pvar}
 \bfY_{ns+\nu} = \sum_{k=1}^{p(\nu)} \bfPhi_k(\nu) \bfY_{ns+\nu-k} + \bfepsilon_{ns+\nu}.
\end{equation}
The process $\bfepsilon:=(\bfepsilon_{t})_t=(\bfepsilon_{ns+\nu})_{n\in \mathbb Z}$ can be interpreted as in \cite{frs11} as the linear innovation of $\bfY:=(\bfY_t)_t=(\bfY_{ns+\nu})_{n\in \mathbb Z}$, \textit{i.e.} $\bfepsilon_t=\bfY_t-\mathbb{E}[\bfY_t|\mathcal{H}_\bfY(t-1)]$, where $\mathcal{H}_\bfY(t-1)$ is the Hilbert space generated by $(\bfY_u, u<t)$. The innovation process $\bfepsilon$ is assumed to be a stationary sequence satisfying
\begin{itemize}
\item[\hspace*{1em} {\bf (A0):}] \hspace*{1em}
$\mathbb{E}\left[ \bfepsilon_t\right]=0, \ \mathrm{Var}\left(\bfepsilon_t\right)=\Sigma_{\bfepsilon}(\nu) \text{ and } \mathrm{Cov}\left(\bfepsilon_t,\bfepsilon_{t-h}\right)=0$ for all $t\in\mathbb{Z}$ and all $h\neq 0$. The covariance matrix $\Sigma_{\bfepsilon}(\nu)$ is assumed to be non-singular.
\end{itemize}
Under the above assumptions, the process $(\bfepsilon_{ns+\nu})_{n\in\mathbb{Z}}$ is called a weak multivariate periodic white noise. It is accepted that  $(\bfY_{ns+\nu})_{n\in\mathbb{Z}}$ admits a strong PVAR representation, if in~\eqref{pvar} $(\bfepsilon_{ns+\nu})_{n\in\mathbb{Z}}$
is a  strong multivariate periodic white noise, namely an iid sequence of random variables with mean $0$ and common variance matrix. In contrast with this previous definition, the representation \eqref{pvar} is called a weak PVAR if no additional assumption is made on $(\bfepsilon_{ns+\nu})_{n\in\mathbb{Z}}$, that is, if $(\bfepsilon_{ns+\nu})_{n\in\mathbb{Z}}$ is only a weak periodic white noise (not necessarily iid).

To derive some basic properties, it is convenient to write the model~(\ref{pvar}) in VAR representation:
\begin{equation}
\label{VARrepre}
\bfPhi_0^{\ast} \bfY_n^{\ast} = \sum_{k=1}^{p^*} \bfPhi_k^{\ast} \bfY_{n-k}^{\ast} +
                                \bfepsilon_n^{\ast},
\end{equation}
where
$\bfY_n^{\ast}= ( \bfY_{ns+s}^{\top}, \bfY_{ns+s-1}^{\top}, \ldots, \bfY_{ns+1}^{\top} )^{\top}$
and
$\bfepsilon_n^{\ast}= ( \bfepsilon_{ns+s}^{\top}, \bfepsilon_{ns+s-1}^{\top}, \ldots, \bfepsilon_{ns+1}^{\top} )^{\top}$
are
$(ds) \times 1$
random vectors.
The autoregressive model order in~(\ref{VARrepre})
is given by
$p^{\ast} = \lceil p/s \rceil$,
where
$\lceil x \rceil$
denotes the smallest integer greater than or equal to the real number $x$.
The
matrix
$\bfPhi_0^{\ast}$,
and the
autoregressive coefficients
$\bfPhi_k^{\ast}$, $k=1,\ldots,p^{\ast}$,
all of dimension
$(ds) \times (ds)$,
are given by the non-singular matrix:
\[
\bfPhi_0^{\ast}=
  \left [ \begin{array}{cccccc}
\bfI_d & -\bfPhi_1(s) & -\bfPhi_2(s) & \ldots & -\bfPhi_{s-2}(s) & -\bfPhi_{s-1}(s)\\
\bfzero & \bfI_d & -\bfPhi_1(s-1) & \ldots & -\bfPhi_{s-3}(s-1) & -\bfPhi_{s-2}(s-1)\\
\vdots &  & & \ddots & & \vdots  \\
\bfzero & \bfzero & \bfzero & \ldots & \bfI_d & -\bfPhi_1(2)\\
\bfzero & \bfzero & \bfzero & \ldots & \bfzero & \bfI_d
          \end{array}
  \right ],
\]
where $\bfI_d$ denotes the $d \times d$ identity matrix, and:
\[
\bfPhi_k^{\ast}=
   \left [ \begin{array}{cccc}
\bfPhi_{ks}(s) & \bfPhi_{ks+1}(s) & \ldots & \bfPhi_{ks+s-1}(s) \\
\bfPhi_{ks-1}(s-1) & \bfPhi_{ks}(s-1) & \ldots & \bfPhi_{ks+s-2}(s-1)\\
\vdots & & \ddots & \vdots \\
\bfPhi_{ks-s+1}(1) & \bfPhi_{ks-s+2}(1) &  \ldots & \bfPhi_{ks}(1)
           \end{array}
    \right ],
\]
where $k = 1,2,\ldots,p^{\ast}$ and $\bfPhi_k(\nu) = \bfzero$, $k > p$.

Let $\det (\bfA)$ be the determinant of the squared matrix $\bfA$. Using general properties of VAR models, it follows that the multivariate stochastic process $\{ \bfY_t^{\ast} \}$ is causal if:
\begin{itemize}
\item[\hspace*{1em} {\bf (A1):}]
\hspace*{1em} $\det \left(\bfPhi_0^{\ast} - \bfPhi_1^{\ast} z - \ldots - \bfPhi_{p^{\ast}} z^{p^{\ast}} \right) \neq 0,$
for all complex numbers $z$ satisfying the condition
$|z|\leq 1$.
\end{itemize}

Under Assumption {\bf (A1)}, there exists a sequence of constant matrices $(\bfF_i(\nu))_{i\geq0}$ such that, for $\nu=1,2,\dots,s$,
$\sum_{i=0}^{\infty}\|\bfF_i(\nu)\|<\infty$ with $\bfF_0(\nu)=\bfI_d$ and
\begin{equation}
\label{MAinf}
\bfY_{ns+\nu}=\sum_{i=0}^{\infty}\bfF_i(\nu)\bfepsilon_{ns+\nu-i},
\end{equation}
where the sequence of matrices  $\|\bfF_i(\nu)\|\to 0$ at a geometric rate as $i\to\infty$. The $\|\bfA\|$ denotes the Euclidian norm of the matrix $\bfA$, that is $\|\bfA\| = \{  \tr (\bfA \bfA^{\top}) \}^{1/2}$, with $\tr(\bfB)$ being the trace of the squared matrix $\bfB$.

Using the algebraic equivalence between multivariate stationarity and periodic correlation~\citep{Gl61}, the $ds$-dimensional process $\{ \bfY_n^{\ast} \}$ is stationary if and only if the $d$-dimensional process $\{ \bfY_t \}$ is periodic stationary with period $s$, in the sense that:
\[
   \cov(\bfY_{n+s},\bfY_{m+s}) = \cov(\bfY_n,\bfY_m),
\]
for all integers $n$ and $m$.~\cite{UD09} show that the autocovariances of $\bfY_{ns+\nu}$ generally depend on $s$ so that the process is not stationary as defined by~\cite{Re97}.

To establish the consistency  of the least squares estimators, an additional assumption is needed.
\begin{itemize}
\item[\hspace*{1em} {\bf (A2):}]
\hspace*{1em}  The $ds$-dimensional process $\left(\boldsymbol{\epsilon}^{\ast}_n\right)_{n\in \mathbb{Z}}$ is ergodic and strictly  stationary.
\end{itemize}
Note that Assumption {\bf (A2)} is entailed by an iid assumption on  $\boldsymbol{\epsilon}^{\ast}_n$, but not by Assumption {\bf (A0)}.

For the asymptotic normality of  least squares estimators, additional assumptions are also required. To control the serial dependence of the stationary process $(\boldsymbol{\epsilon}^{\ast}_n)_{n\in\mathbb{Z}}$, we introduce the strong mixing coefficients $\alpha_{\boldsymbol{\epsilon}^{\ast}}(h)$ defined by
$$\alpha_{\boldsymbol{\epsilon}^{\ast}}\left(h\right)=\sup_{A\in\mathcal F^n_{-\infty},B\in\mathcal F_{n+h}^{+\infty}}\left|\mathbb{P}\left(A\cap
B\right)-\mathbb{P}(A)\mathbb{P}(B)\right|,$$
where $\mathcal F_{-\infty}^n=\sigma (\boldsymbol{\epsilon}^{\ast}_u, u\leq n )$ and $\mathcal F_{n+h}^{+\infty}=\sigma (\boldsymbol{\epsilon}^{\ast}_u, u\geq n+h )$.
We will make an integrability assumption on the moment of the noise and a summability condition on the strong mixing coefficients $(\alpha_{\boldsymbol{\epsilon}^{\ast}}(k))_{k\ge 0}$.
\begin{itemize}
\item[\hspace*{1em} {\bf (A3):}]
\hspace*{1em}
$\text{We have }\mathbb{E}\|{\boldsymbol{\epsilon}}^{\ast}_n\|^{4+2\kappa}<\infty\text{ and }
\sum_{k=0}^{\infty}\left\{\alpha_{{\boldsymbol{\epsilon}}^{\ast}}(k)\right\}^{\frac{\kappa}{2+\kappa}}<\infty\text{ for some } \kappa>0.$
\end{itemize}


\section{Estimating weak PVAR models}\label{estimation}
In this section, we recall the technical problems of the least squares estimation method for weak PVAR models as well as the asymptotic properties of the least squares estimators. For periodic models, there are at least two more choices in terms of methods of estimation: the maximum likelihood procedure~\citep{LB00} and the Yule-Walker equations~\citep{Pa78}.

Consider the time series data
$\bfY_{ns+\nu}$, $n=0,1,\ldots,N-1$, $\nu = 1,\ldots,s$,
giving a sample size equal to $n = Ns$. Let
\begin{eqnarray}
\label{Znu}
\bfZ(\nu) &=& \left( \bfY_{\nu}, \bfY_{s+\nu}, \ldots, \bfY_{(N-1)s+\nu} \right), \\
\label{Enu}
\bfE(\nu) &=& \left( \bfepsilon_{\nu}, \bfepsilon_{s+\nu}, \ldots, \bfepsilon_{(N-1)s+\nu} \right),\\
\label{Xnu}
\bfX(\nu) &=& \left( \bfX_0(\nu), \ldots, \bfX_{N-1}(\nu) \right),
\end{eqnarray}
be $d \times N$, $d \times N$ and $\{ d p(\nu) \} \times N$ random matrices, where
\[
\bfX_n(\nu) = (\bfY_{ns+\nu-1}^{\top},\ldots,\bfY_{ns+\nu-p(\nu)}^{\top})^{\top},\]
$n=0,1,\ldots,N-1,$ denote $\{ d p(\nu) \} \times 1$ random vectors. The PVAR model can be reformulated as:
\begin{equation}
\label{modelZnu}
 \bfZ(\nu) = \bfB(\nu) \bfX(\nu) + \bfE(\nu), \; \nu = 1,\ldots,s,
\end{equation}
where the model parameters are collected in the $d \times \{ dp(\nu) \}$ matrix $\bfB(\nu)$ which is defined as:
\begin{equation}
\label{Bnu}
\bfB(\nu) = \left( \bfPhi_1(\nu), \ldots, \bfPhi_{p(\nu)}(\nu) \right).
\end{equation}
In what follows, the symbols '$\cd$' and '$\cas$' stand for convergence in distribution and almost surely, respectively, and $\mathcal{N}_d(\bfmu, \bfSigma)$ denotes a $d$-dimensional normal distribution with mean $\bfmu$ and covariance matrix $\bfSigma$. Denote by $\vec(\bfA)$ the vector obtained by stacking the columns of $\bfA$.
\subsection{Unconstrained least squares estimators}\label{LSnoncont}
\noindent The least squares estimators based on the generalized least squares criteria are obtained equivalently by minimizing the ordinary least squares (a similar result holds for VAR models, see~\citet[p.71]{Lu05}):
\begin{equation}
\label{ols}
  S(\bfbeta) = \sum_{\nu=1}^s \bfe^{\top}(\nu)   \bfe(\nu),
\end{equation}
where $\bfbeta = (\bfbeta^{\top}(1),\ldots,\bfbeta^{\top}(s))^{\top}$ is the
$\{ d^2 \sum_{\nu=1}^s p(\nu) \} \times 1$ vector of model parameters and $\bfbeta(\nu) = \vec\{ \bfB(\nu) \}$ is a $\{ d^2 p(\nu) \} \times 1$ vector of parameters for $\nu=1,\ldots,s$.

To obtain the least squares estimators, we differentiate $S(\bfbeta)$ with respect to each parameter $\bfPhi_k(\nu)$, $k=1,\ldots,p(\nu)$, $\nu=1,\ldots,s$. Thus we obtain easily:
\[
  \frac{\partial S(\bfbeta)}{\partial \vec\{\bfPhi_k(\nu) \}} =
  -2\sum_{n=0}^{N-1} (\bfY_{ns+\nu-k} \otimes \bfepsilon_{ns+\nu}), \; k=1,\ldots,p(\nu), \; \nu=1,\ldots,s.
\]
Setting the derivatives equal to zero, $k=1,\ldots,p(\nu)$, gives the following system for a given season $\nu$:
\[
 \sum_{n=0}^{N-1} \left\{ \bfX_n(\nu) \otimes \bfepsilon_{ns+\nu} \right\} = \bfzero,
\]
where $\bfzero$ is the $\{ d^2 p(\nu) \} \times 1$ null vector. Since
$\bfepsilon_{ns+\nu} = \bfY_{ns+\nu} - \{ \bfX_n^{\top}(\nu) \otimes \bfI_d \} \bfbeta(\nu)$,
the normal equations at season $\nu$ are:
\[
 \sum_{n=0}^{N-1} \{ \bfX_n(\nu) \otimes \bfY_{ns+\nu} \} =
 \left[ \sum_{n=0}^{N-1} \left\{ \bfX_n(\nu)\bfX_n^{\top}(\nu) \otimes \bfI_d \right\} \right] \bfbeta(\nu).
\]
Consequently, the least squares estimators of $\bfbeta(\nu)$ satisfy the relation:
\[
  \hat{\bfbeta}(\nu) = \left[ \{ \bfX(\nu)\bfX^{\top}(\nu) \}^{-1}\bfX(\nu) \otimes \bfI_d \right]\bfz(\nu),
\]
and the residuals are
$\hat{\bfepsilon}_{ns+\nu} = \bfY_{ns+\nu} - \{ \bfX_n^{\top}(\nu) \otimes \bfI_d \} \hat{\bfbeta}(\nu)$. Thus the following relation is satisfied:
\begin{align}
\nonumber
  N^{1/2}\{ \hat{\bfbeta}(\nu) - \bfbeta(\nu) \} &=
  \left[ \{ N^{-1}\bfX(\nu) \bfX^{\top}(\nu) \}^{-1} \otimes \bfI_d \right] N^{-1/2} \{ \bfX(\nu) \otimes \bfI_d \} \bfe(\nu),\\ \label{betahatnu}
 & =
  \left[ \{ N^{-1}\bfX(\nu) \bfX^{\top}(\nu) \}^{-1} \otimes \bfI_d \right] N^{-1/2} \vec\{ \bfE(\nu) \bfX^\top(\nu)  \},
\end{align}
Note also that, by using the properties of the $\vec(\cdot)$ operator, it should be noted that an alternative expression for the least squares estimators is given by:
\begin{equation}
\label{hatBnu}
  \hat{\bfB}(\nu) = \bfZ(\nu) \bfX^{\top}(\nu) \{ \bfX(\nu)\bfX^{\top}(\nu) \}^{-1}.
\end{equation}
Under the assumptions {\bf (A0)}, {\bf (A1)}, {\bf (A2)} and {\bf (A3)}, for $\nu=1,\ldots,s$, \cite{BMU23} studied the asymptotic properties of the least squares estimators in the unrestricted case and showed that:
\begin{eqnarray}
\label{th1b}
  \hat{\bfbeta}(\nu) &\cas& \bfbeta(\nu), \\
\label{th1a}
  N^{-1/2} \sum_{n=0}^{N-1} \vec\{ \bfepsilon_{ns+\nu} \bfX_n^{\top}(\nu) \} &\cd&
  \mathcal{N}_{d^2p(\nu)}\left( \bfzero,\bfPsi_{\nu\nu}
  \right), \\ \nonumber
  \bfPsi_{\nu\nu}&=& \sum_{h=-\infty}^{\infty}\mathbb{E}\left(\bfX_n(\nu)\bfX_{n-h}^{\top}(\nu) \otimes \bfepsilon_{ns+\nu}\bfepsilon_{(n-h)s+\nu}^{\top}\right)\\
\label{th1c}
  N^{1/2}\{ \hat{\bfbeta}(\nu) - \bfbeta(\nu) \}
  &\cd& \mathcal{N}_{d^2p(\nu)}\left(\bfzero, \bfTheta_{\nu\nu} \right),
  \\ \nonumber
  \bfTheta_{\nu\nu}&=&\left(\bfOmega^{-1}(\nu) \otimes \bfI_d\right)\bfPsi_{\nu\nu}\left(\bfOmega^{-1}(\nu) \otimes \bfI_d\right)\\
  \label{th1d}
  N^{1/2}\{ \hat{\bfbeta} - \bfbeta\}
  &\cd& \mathcal{N}_{sd^2p(\nu)}\left(\bfzero, \bfTheta \right),
\end{eqnarray}
where $\bfOmega(\nu)=:\mathbb{E}\left(\bfX_n(\nu)\bfX_n^{\top}(\nu)\right)$ corresponds to the
$\{ d p(\nu) \} \times \{ d p(\nu) \}$ covariance matrix of the $\{ d p(\nu) \} \times 1$
random vector $\bfX_n(\nu)$ and where the asymptotic covariance matrix $\bfTheta$ is a block matrix, with blocs given by: $$\bfTheta_{\nu\nu'}:=\left(\bfOmega^{-1}(\nu) \otimes \bfI_d\right)\sum_{h=-\infty}^{\infty}\mathbb{E}\left(\bfX_n(\nu)\bfX_{n-h}^{\top}(\nu') \otimes \bfepsilon_{ns+\nu}\bfepsilon_{(n-h)s+\nu'}^{\top}\right)\left(\bfOmega^{-1}(\nu') \otimes \bfI_d\right),\quad\text{for }\quad\nu, \nu' = 1,\ldots,s.$$

\subsection{Least squares estimation with linear constraint on the parameters.}\label{LScont}
\noindent The PVAR model in~(\ref{pvar}) has $d^2 \sum_{\nu=1}^{s} p(\nu)$ autoregressive parameters $\bfPhi_k(\nu)$, $k=1,\ldots,p(\nu)$, $\nu = 1,\ldots,s$, and $s$ additional $d \times d$ covariance matrices $\bfSigma_{\bfepsilon}(\nu)$, $\nu = 1,\ldots,s$. For multivariate processes, the number of parameters can be quite large; for vector periodic processes, the inflation of parameters is due to the $s$ seasons, which makes some PVAR inference matters unwieldy. Consequently, many authors have investigated parsimonious versions of \eqref{pvar} (see, for instance, \cite{UD09, DL13}).  In this section, we consider estimation in the unrestricted case but also in the situation where the parameters of the same season $\nu$ satisfy the relation:
\begin{equation}
\label{constraints}
\bfbeta(\nu) = \bfR(\nu) \bfxi(\nu) + \bfb(\nu),
\end{equation}
where $\bfR(\nu)$ is a known $\{ d^2 p(\nu) \} \times K(\nu)$ matrix of rank $K(\nu)$, $\bfb(\nu)$ a known $\{ d^2 p(\nu) \} \times 1$ vector and $\bfxi(\nu)$ represents a $K(\nu) \times 1$ vector of unknown parameters. Letting $\bfR(\nu) = \bfI_{d^2p(\nu)}$, $\bfb(\nu) = \bfzero$, $\nu = 1,\ldots,s$ give what we call the full unconstrained case of Section \ref{LSnoncont}. In general, the matrices $\bfR(\nu)$ and the vectors $\bfb(\nu)$ allow for linear constraints on the parameters of the same season $\nu$, $\nu = 1,\ldots,s$.

This linear constraint includes the important special case of parameters set to zero on certain components of $\bfPhi_k(\nu)$, $\nu = 1,\ldots,s$. In practice, a two-step procedure could consist of fitting a full unconstrained model, and, in a second stage of inference, the estimators that are statistically not significant could be considered known zero parameters, providing frequently more parsimonious models.

Vectorizing \eqref{modelZnu}, we obtain:
\begin{eqnarray}
 \bfz(\nu) &=& \{ \bfX^{\top}(\nu) \otimes \bfI_d \} \vec\{ \bfB(\nu) \} + \vec\{ \bfE(\nu) \}, \nonumber \\
           &=& \{ \bfX^{\top}(\nu) \otimes \bfI_d \} \bfbeta(\nu) + \bfe(\nu), \nonumber \\
\label{enu}
           &=& \{ \bfX^{\top}(\nu) \otimes \bfI_d \} \{ \bfR(\nu)\bfxi(\nu) + \bfb(\nu) \} + \bfe(\nu),
\end{eqnarray}
where $\bfz(\nu) = \vec\{ \bfZ(\nu) \}$ and $\bfe(\nu) = \vec\{ \bfE(\nu) \}$. The covariance matrix of the random vector $\bfe(\nu)$ is $\bfI_N \otimes \bfSigma_{\bfepsilon}(\nu)$.

The least squares estimators of $\bfxi(\nu)$, $\nu = 1,\ldots,s$ are obtained by minimizing
the generalized least squares criterion:
\begin{equation}
\label{GLS}
  S_G(\bfxi) = \sum_{\nu=1}^s \bfe^{\top}(\nu) \{ \bfI_N \otimes \bfSigma_{\bfepsilon}(\nu) \}^{-1}  \bfe(\nu),
\end{equation}
where $\bfxi = (\bfxi^{\top}(1),\ldots,\bfxi^{\top}(s))^{\top}$ represents a
$\{ \sum_{\nu=1}^s K(\nu) \} \times 1$ vector. When the parameters satisfy the linear constraint~(\ref{constraints}), the least squares estimators of $\bfxi(\nu)$, $\nu = 1,\ldots,s$, minimize the generalized criterion~(\ref{GLS}), which is not equivalent to~(\ref{ols}), see~\citet{Lu05}, amongst others. Recall that from~(\ref{enu}) we have the following relation:
\[
  \bfe(\nu) = \bfz(\nu) - \{ \bfX^{\top}(\nu) \otimes \bfI_d \} \{ \bfR(\nu) \bfxi(\nu) + \bfb(\nu) \},
\]
which is convenient to derive the asymptotic properties of the least squares estimator of $\bfxi(\nu)$.

\noindent Proceeding as in the previous section, it is possible to show that the least squares estimator $\hat{\bfxi}(\nu)$ of $\bfxi(\nu)$ is given by:
\begin{align*}
\hat{\bfxi}(\nu) &=
  \left[ \bfR^{\top}(\nu) \{ \bfX(\nu)\bfX^{\top}(\nu) \otimes \bfSigma^{-1}_{\bfepsilon}(\nu) \} \bfR(\nu) \right]^{-1}
  \bfR^{\top}(\nu) \{ \bfX(\nu)  \otimes \bfSigma^{-1}_{\bfepsilon}(\nu) \}
 \left[\bfz(\nu) - \{ \bfX^{\top}(\nu) \otimes \bfI_d \}\bfb(\nu) \right].
\end{align*}
Furthermore, the following relation is satisfied:
\begin{align}\nonumber
  N^{1/2} &\{  \hat{\bfxi}(\nu) - \bfxi(\nu) \} \\ \nonumber&=
  N^{1/2} \left[ \bfR^{\top}(\nu) \{ \bfX(\nu)\bfX^{\top}(\nu) \otimes \bfSigma^{-1}_{\bfepsilon}(\nu) \} \bfR(\nu) \right]^{-1}
  \bfR^{\top}(\nu) \{ \bfX(\nu)  \otimes \bfSigma^{-1}_{\bfepsilon}(\nu) \} \bfe(\nu)\\  \label{ecartxi}
  &= \left[ \bfR^{\top}(\nu)\frac{1}{N} \{ \bfX(\nu)\bfX^{\top}(\nu) \otimes \bfSigma^{-1}_{\bfepsilon}(\nu) \} \bfR(\nu) \right]^{-1}
  \bfR^{\top}(\nu)
  \{ \bfI_{dp(\nu)}  \otimes \bfSigma^{-1}_{\bfepsilon}(\nu) \} N^{-1/2} \vec\{ \bfE(\nu) \bfX^\top(\nu)  \}.
\end{align}

Under the assumptions {\bf (A0)}, {\bf (A1)}, {\bf (A2)} and {\bf (A3)}, for $\nu=1,\ldots,s$, \cite{BMU23} showed that, the estimator $\hat{\bfxi}(\nu)$ is consistent for $\bfxi(\nu)$, and $\hat{\bfxi}(\nu)$ follows asymptotically a normal distribution, that is:
\begin{align}
\label{NAcont}
N^{1/2} \{  \hat{\bfxi}(\nu) - \bfxi(\nu) \} \cd
  \mathcal{N}_{K(\nu)}\left(\bfzero,
              \bfTheta^{\bfxi}_{\nu\nu}
            \right),
\end{align}
where
\begin{align*}
  \bfTheta^{\bfxi}_{\nu\nu}&=\left[ \bfR^{\top}(\nu) \{ \bfOmega(\nu) \otimes \bfSigma^{-1}_{\bfepsilon}(\nu) \} \bfR(\nu) \right]^{-1}
\bfR^{\top}(\nu) \{ \bfI_{dp(\nu)} \otimes \bfSigma^{-1}_{\bfepsilon}(\nu) \}  \bfPsi_{\nu\nu} \\&\hspace{6cm}\times\{ \bfI_{dp(\nu)} \otimes \bfSigma^{-1}_{\bfepsilon}(\nu) \}   \bfR(\nu)
 \left( \left[ \bfR^{\top}(\nu) \{ \bfOmega(\nu) \otimes \bfSigma^{-1}_{\bfepsilon}(\nu) \} \bfR(\nu) \right]^{-1}\right)^{\top}
 \\
 &=\left[ \bfR^{\top}(\nu) \{ \bfOmega(\nu) \otimes \bfSigma^{-1}_{\bfepsilon}(\nu) \} \bfR(\nu) \right]^{-1}
\bfR^{\top}(\nu)  \sum_{h=-\infty}^{\infty}\mathbb{E}\left[\bfX_n(\nu)\bfX_{n-h}^{\top}(\nu) \otimes \bfSigma^{-1}_{\bfepsilon}(\nu)\bfepsilon_{ns+\nu}\bfepsilon_{(n-h)s+\nu}^{\top}\bfSigma^{-1}_{\bfepsilon}(\nu)\right]    \\&\hspace{6cm}\qquad\times   \bfR(\nu)
 \left( \left[ \bfR^{\top}(\nu) \{ \bfOmega(\nu) \otimes \bfSigma^{-1}_{\bfepsilon}(\nu) \} \bfR(\nu) \right]^{-1}\right)^{\top}.
\end{align*}
Furthermore, we also have
\begin{eqnarray}\label{th1dCon}
N^{1/2} \{  \hat{\bfxi} - \bfxi \} \cd \mathcal{N}_{sK(\nu)}\left(\bfzero, \bfTheta^{\mathrm{\bfxi}} \right),
\end{eqnarray}
where the asymptotic covariance matrix $\bfTheta^{\mathrm{\bfxi}} = [\bfTheta^{\bfxi}_{\nu\nu'}]_{1\leq \nu,\nu'\leq s}$ is a block matrix, with blocs given by:
\begin{align*}
  \bfTheta^{\bfxi}_{\nu\nu'}&=\left[ \bfR^{\top}(\nu) \{ \bfOmega(\nu) \otimes \bfSigma^{-1}_{\bfepsilon}(\nu) \} \bfR(\nu) \right]^{-1}
\bfR^{\top}(\nu)  \sum_{h=-\infty}^{\infty}\mathbb{E}\left[\bfX_n(\nu)\bfX_{n-h}^{\top}(\nu') \otimes \bfSigma^{-1}_{\bfepsilon}(\nu)\bfepsilon_{ns+\nu}\bfepsilon_{(n-h)s+\nu'}^{\top}\bfSigma^{-1}_{\bfepsilon}(\nu')\right]    \\&\hspace{7cm}\qquad\times   \bfR(\nu')
 \left( \left[ \bfR^{\top}(\nu') \{ \bfOmega(\nu') \otimes \bfSigma^{-1}_{\bfepsilon}(\nu') \} \bfR(\nu') \right]^{-1}\right)^{\top}.
\end{align*}
It should be noted that the estimator $\hat{\bfxi}(\nu)$ is unfeasible in practice, since it relies on the unknown matrix $\bfSigma_{\bfepsilon}(\nu)$. A feasible estimator is given by:
\begin{align*}
\hat{\hat{\bfxi}}(\nu) &=
  \left[ \bfR^{\top}(\nu) \{ \bfX(\nu)\bfX^{\top}(\nu) \otimes \tilde{\bfSigma}^{-1}_{\bfepsilon}(\nu) \} \bfR(\nu) \right]^{-1}
  \bfR(\nu) \{ \bfX(\nu)  \otimes \tilde{\bfSigma}^{-1}_{\bfepsilon}(\nu) \}
\\& \hspace{7cm}\times   [\bfz(\nu) - \{ \bfX^{\top}(\nu) \otimes \bfI_d \}\bfb(\nu)],
\end{align*}
\noindent where $\tilde{\bfSigma}_{\bfepsilon}(\nu)$ denotes a consistent estimator of the covariance matrix $\bfSigma_{\bfepsilon}(\nu)$ for $\nu=1,\ldots,s$.
A possible candidate is obtained from the unconstrained least squares estimators:
\[
\tilde{\bfSigma}_{\bfepsilon}(\nu) = \{ N-d p(\nu) \}^{-1} \left\{ \bfZ(\nu) - \hat{\bfB}(\nu)\bfX(\nu) \right\}
                                     \left\{ \bfZ(\nu) - \hat{\bfB}(\nu)\bfX(\nu) \right\}^{\top},
\]
where $\hat{\bfB}(\nu)$ represents the unconstrained least squares estimators~(\ref{hatBnu}) obtained in Section~\ref{LSnoncont}. The resulting estimator of
$\bfbeta(\nu)$ is given by
$\hat{\hat{\bfbeta}}(\nu) = \bfR(\nu)\hat{\hat{\bfxi}}(\nu) + \bfb(\nu)$, and its asymptotic distribution  is normal:
\begin{eqnarray}
\label{NAcontbeta}
  N^{1/2} \{  \hat{\hat{\bfbeta}}(\nu) - \bfbeta(\nu) \} \cd
  \mathcal{N}_{d^2p(\nu)}\left(\bfzero,
             \bfR(\nu)
              \bfTheta^{\bfxi}_{\nu\nu}
             \bfR^{\top}(\nu)
            \right).
\end{eqnarray}

\section{Diagnostic checking in weak PVAR models}\label{result}
\noindent After the estimation phase, the next important step consists in checking if the estimated model fits satisfactorily the data. In this section we derive the limiting distribution of the residual  autocovariance and autocorrelation matrices in the framework of weak PVAR models with unconstrained and also with linear constraints defined by~\eqref{constraints} on the parameters of a given season. From these results, we propose modified test statistics for checking the adequacy of a weak PVAR model. Note that the results stated in this section extend directly for strong PVAR models. For the class of strong PVAR models with parameter constraints,~\cite{UD09, DL13} derives the asymptotic distribution of the residual autocovariance and autocorrelation matrices, which follow normal distributions.
\subsection{Asymptotic distribution of the residual autocovariance and autocorrelation matrices}\label{result1}
For any
$\dot{\bfbeta}(\nu) = \bfR(\nu) \dot{\bfxi}(\nu) + \bfb(\nu)$,
where
$\dot{\bfbeta}(\nu) = (\vec^{\top}\{ \dot{\bfPhi}_1(\nu) \},\ldots, \vec^{\top}\{ \dot{\bfPhi}_{p(\nu)}(\nu) \} )^{\top}$,
with
general
$d \times d$
matrices
$\dot{\bfPhi}_k(\nu)$,
$k=1,\ldots,p(\nu)$,
we introduce the model residuals:
\[
\label{autocov}
\dot{\bfepsilon}_{ns+\nu} = \left\{
      \begin{array}{ll}
         \bfY_{ns+\nu} - \sum_{k=1}^{p(\nu)} \dot{\bfPhi}_k(\nu) \bfY_{ns+\nu-k}, &
         ns+\nu > p(\nu), \\
         \bfzero, & ns+\nu \leq p(\nu),
      \end{array}
           \right.
\]
which are well-defined for
$n=0,1,\ldots,N-1$.
As in Li and McLeod (1981),
we use the dot notation to designate
the residuals
$\dot{\bfepsilon}_{ns+\nu}$, $n=0,1,\ldots,N-1$,
expressed in function of the general quantities
$\dot{\bfbeta}(\nu)$
and
$\dot{\bfxi}(\nu)$,
$\nu=1,\ldots,s$.
Using the well-known relation
$\vec(\bfA \bfB \bfC) = (\bfC^{\top} \otimes \bfA) \vec(\bfB)$,
where
'$\otimes$'
represents the Kronecker product,
we obtain for
$ns+\nu > p(\nu)$
the relation
$\dot{\bfepsilon}_{ns+\nu} =
\bfY_{ns+\nu} - \sum_{k=1}^{p(\nu)} (\bfY_{ns+\nu-k}^{\top}  \otimes \bfI_d) \vec\{ \dot{\bfPhi}_k(\nu) \}$.
Let
$\bfGamma_{\bfepsilon}(h; \nu) = \cov(\bfepsilon_{ns+\nu}, \bfepsilon_{ns+\nu-h})$
be the lag
$h$
theoretical autocovariance matrix at season
$\nu$
of the error process
$\bfepsilon$.
For a fixed integer $M\geq1$, we let
\begin{equation}
\label{gammanu}
\bfgamma_{\bfepsilon}(\nu)  =
     \left(
       \bfgamma_{\bfepsilon}^{\top}(1;\nu),
       \bfgamma_{\bfepsilon}^{\top}(2;\nu), \ldots,
       \bfgamma_{\bfepsilon}^{\top}(M;\nu)
     \right)^{\top}
\end{equation}
be a $(d^2M) \times 1$ vector of theoretical autocovariances,
where
$\bfgamma_{\bfepsilon}(h;\nu) = \vec\{ \bfGamma_{\bfepsilon}(h;\nu) \}$.
We define
$\bfrho_{\bfepsilon}(h;\nu) = \bfGamma_0^{-1}(\nu)\bfGamma_{\bfepsilon}(h; \nu)\bfGamma_0^{-1}(\nu-h)$
the lag $h$ theoretical autocorrelation matrix at season $\nu$,
using the diagonal matrix
$\bfGamma_0(\nu) = \diag(\bfsigma_{\bfepsilon,11}^{1/2}(\nu),\ldots,\bfsigma_{\bfepsilon,dd}^{1/2}(\nu))$. Similarly, for a fixed integer $M\geq1$, we let
\begin{equation}
\label{rhonu}
\bfrho_{\bfepsilon}(\nu)  =
     \left(
       \left\{\vec\left(\bfrho_{\bfepsilon}(1;\nu)\right)\right\}^{\top},
       \left\{\vec\left(\bfrho_{\bfepsilon}(2;\nu)\right)\right\}^{\top}, \ldots,
       \left\{\vec\left(\bfrho_{\bfepsilon}(M;\nu)\right)\right\}^{\top}
     \right)^{\top}
\end{equation}
be a $(d^2M) \times 1$ vector of theoretical autocorrelation.
We introduce the sample autocovariance matrices
$\bfC_{\dot{\bfepsilon}}(h;\nu) =
\left( C_{\dot{\bfepsilon},ij}(h; \nu) \right)_{i,j=1,\ldots,d}$:
\[
\bfC_{\dot{\bfepsilon}}(h;\nu) = \left\{ \begin{array}{ll}
                                   N^{-1} \sum_{n=h}^{N-1} \dot{\bfepsilon}_{ns+\nu} \dot{\bfepsilon}_{ns+\nu-h}^{\top}, &
                                     h \geq 0, \\
                                  \bfC_{\dot{\bfepsilon}}^{\top}(-h;\nu-h), & h < 0.
                       \end{array}
         \right.
\]
For a fixed integer $M\geq1$, we let
\begin{equation}
\label{cnu}
\bfc_{\dot{\bfepsilon}}(\nu) =
(\bfc_{\dot{\bfepsilon}}^{\top}(1; \nu),\ldots,\bfc_{\dot{\bfepsilon}}^{\top}(M; \nu))^{\top}\quad\text{and}\quad
\bfc_{\dot{\bfepsilon}} =
(\bfc_{\dot{\bfepsilon}}^{\top}(1),\ldots,\bfc_{\dot{\bfepsilon}}^{\top}(M))^{\top}
\end{equation}
be the
$(d^2M) \times 1$
vector of sample autocovariances and the
$(d^2Ms) \times 1$
vector of global sample autocovariances,
where
$\bfc_{\dot{\bfepsilon}}(h; \nu) = \vec\{ \bfC_{\dot{\bfepsilon}}(h; \nu) \}$.
Here
$M$
represents a fixed integer with respect to the sample size
$n=Ns$,
satisfying the relation
$1 \leq M < N$;
this constant is the maximal lag order.
Similarly,
the vector of sample autocorrelations is given by
$\bfr_{\dot{\bfepsilon}}(\nu) =
(\bfr_{\dot{\bfepsilon}}^{\top}(1; \nu),\ldots,\bfr_{\dot{\bfepsilon}}^{\top}(M; \nu))^{\top}$,
where:
\begin{equation}\label{r_hnu}
\bfr_{\dot{\bfepsilon}}(h; \nu) = \vec\left\{ \bfD_{\dot{\bfepsilon}}^{-1}(\nu)\bfC_{\dot{\bfepsilon}}(h; \nu)
                                    \bfD_{\dot{\bfepsilon}}^{-1}(\nu-h) \right\}
                                = \left( \bfD_{\dot{\bfepsilon}}^{-1}(\nu-h)  \otimes
                                    \bfD_{\dot{\bfepsilon}}^{-1}(\nu) \right)
                                     \bfc_{\dot{\bfepsilon}}(h;\nu),
\end{equation}
with
$\bfD_{\dot{\bfepsilon}}(\nu) =
\diag\left( \bfC_{\dot{\bfepsilon},11}^{1/2}(0; \nu),\ldots,\bfC_{\dot{\bfepsilon},dd}^{1/2}(0; \nu) \right)$.
We also define the vector of global sample autocorrelations by \\%
$\bfr_{\dot{\bfepsilon}} =
(\bfr_{\dot{\bfepsilon}}^{\top}(1),\ldots,\bfr_{\dot{\bfepsilon}}^{\top}(M))^{\top}$.

We now discuss the asymptotic distribution of
$N^{1/2} \bfc_{\dot{\bfepsilon}}(\nu)$
where there is no constraint on the parameters. Define the matrix
\begin{equation}\label{Ups_M}\bfUpsilon_M(\nu)=-\mathbb{E}\left\{\left(\begin{array}{c}
\bfepsilon_{ns+\nu-1}\\\vdots\\\bfepsilon_{ns+\nu-M}\end{array}\right)\otimes
\bfX_n^\top(\nu)\otimes \bfI_d\right\}.
\end{equation}
By expanding
$\bfc_{\dot{\bfepsilon}}(\nu)$
in a Taylor expansion around
${\bfbeta}(\nu)$
and evaluating
at the point
$\dot{\bfbeta}(\nu) = \hat{\bfbeta}(\nu)$,
we obtain the following development:
\begin{align}\label{hat_gamma}
\sqrt{N}\bfc_{\hat{\bfepsilon}}(\nu) = \sqrt{N}\bfc_{\bfepsilon}(\nu) +
                    \bfUpsilon_M(\nu)\sqrt{N}
                    \left\{
                    \hat{\bfbeta}(\nu) - \bfbeta(\nu)
                    \right\} + \mathrm{o}_{\mathbb{P}}(1).
\end{align}
Let $\bfD^\ast_{\dot{\bfepsilon}}(\nu) =
\diag\left( \bfGamma_{{\bfepsilon},11}^{1/2}(0; \nu),\ldots,\bfGamma_{{\bfepsilon},dd}^{1/2}(0; \nu) \right)$, it also follows that:
\begin{align}
\sqrt{N}\bfr_{\hat{\bfepsilon}}(\nu)
&=
\left(\diag\left[\left\{ \bfD^\ast_{{\bfepsilon}}(\nu-1),\dots,\bfD^\ast_{{\bfepsilon}}(\nu-M)\right\}  \otimes
                                    \bfD^\ast_{{\bfepsilon}}(\nu) \right]\right)^{-1}
\sqrt{N}\bfc_{\hat{\bfepsilon}}(\nu)+\mathrm{o}_{\mathbb{P}}(1).\label{hat_rho}
\end{align}
Thus from \eqref{hat_rho} the asymptotic distribution of the residual autocorrelations $\sqrt{N}\bfr_{\hat{\bfepsilon}}(\nu)$ depends on the
 distribution of $\sqrt{N}\bfc_{\hat{\bfepsilon}}(\nu)$.
In view of \eqref{hat_gamma} the asymptotic distribution of the residual autocovariances $\sqrt{N}\bfc_{\hat{\bfepsilon}}(\nu)$
will be obtained from the joint asymptotic behavior of
$N^{1/2} \left( \hat{\bfbeta}(\nu) - \bfbeta(\nu)\right)$
and
$N^{1/2} \bfc_{\bfepsilon}(\nu)$.
From~(\ref{betahatnu}), we note that we have the following relation:
\begin{align}
\label{betahatnuOmega}
  N^{1/2}\{ \hat{\bfbeta}(\nu) - \bfbeta(\nu) \} =
  \{  \bfOmega^{-1}(\nu) \otimes \bfI_d \} N^{-1/2} \left( \vec\{ \bfE(\nu) \bfX^\top(\nu)  \} \right) + \mathrm{o}_\mathbb{P}(1).
\end{align}
Thus the next proposition states the joint asymptotic distribution of $\{  \bfOmega^{-1}(\nu) \otimes \bfI_d \}N^{-1/2} \vec\{ \bfE(\nu)\bfX^{\top}(\nu) \}$
and
$N^{1/2} \bfc_{\bfepsilon}(\nu)$.
\begin{proposition}\label{loijointeUn}
  Let
  $\bfE(\nu)$
  and
  $\bfX(\nu)$
  be defined by~(\ref{Enu}) and~(\ref{Xnu}),
  respectively.
  Consider a vector of $M$ sample autocovariances collected in the vector
  $\bfc_{\bfepsilon}(\nu)$ given by~(\ref{cnu}).
  Suppose that
  $\{ \bfY_t \}$
  denotes a PVAR process satisfying~(\ref{pvar}). Under the assumptions {\bf (A0)}, {\bf (A1)}, {\bf (A2)} and {\bf (A3)}, for a specified season $\nu=1,\dots,s$ we have:
\begin{eqnarray}\nonumber
  N^{1/2} \left( \hat{\bfbeta}^\top(\nu) - \bfbeta^\top(\nu),\bfc_{\bfepsilon}^\top(\nu) \right)^\top
  &\cd& \mathcal{N}_{d^2(p(\nu)+M)}\left(\bfzero, \bfXi_{\nu\nu}^{\mathrm{U}}\right),
\end{eqnarray}
 with
 \begin{align}\label{specWnUnc}
\bfXi_{\nu\nu}^{\mathrm{U}}=\left(\begin{array}{cc}
                           \bfTheta_{\nu\nu} & \bfG^{\mathrm{U}}_{\nu\nu}\\
                           \left(\bfG^{\mathrm{U}}_{\nu\nu}\right)^{\top}& \bfV_{\nu\nu}
                           \end{array}
                     \right) =\sum_{h=-\infty}^{\infty}\mathbb{E}\left[ \bfW_{ns+\nu}^{\mathrm{U}}\left(\bfW_{(n-h)s+\nu}^{\mathrm{U}}\right)^{\top}\right],
\end{align}
where $\bfTheta_{\nu\nu}$, $\bfV_{\nu\nu}$, $\bfG^{\mathrm{U}}_{\nu\nu}$ are given respectively by \eqref{th1c}, \eqref{VnuM}, \eqref{Gnu} and
\begin{equation}\label{WnUnc}
\bfW_{ns+\nu}^{\mathrm{U}}=\begin{pmatrix}\bfW^{\mathrm{U}}_{1,ns+\nu}\vspace{0.2cm}\\
\bfW_{2,ns+\nu}
\end{pmatrix}=\begin{pmatrix} \bfOmega^{-1}(\nu)\bfX_n(\nu) \otimes \bfepsilon_{ns+\nu}\vspace{0.2cm}\\
\left(\bfepsilon_{ns+\nu-1}^\top,\dots,\bfepsilon_{ns+\nu-M}^\top\right)^{\top}\otimes\bfepsilon_{ns+\nu}
\end{pmatrix}.
\end{equation}
Finally, for $\nu,\nu'=1,\dots,s$ the
$(d^2 M) \times (d^2 M)$
  matrix $\bfV_{\nu\nu'}(M)$is given by:
\begin{equation}
\label{VnuM}
  \bfV_{\nu\nu'} = \left[\bfV_{\nu\nu'}(l,l')\right]_{1\leq l,l'\leq M}\text{ where }\bfV_{\nu\nu'}(l,l')= \sum_{h=-\infty}^{\infty}
                            \mathbb{E}(\bfepsilon_{ns+\nu-l} \bfepsilon_{(n-h)s+\nu'-l'}^{\top} \otimes\bfepsilon_{ns+\nu} \bfepsilon_{(n-h)s+\nu'}^{\top})
\end{equation}
and collecting the $\{ d^2 p(\nu) \} \times d^2$ matrix
$\bfG^{\mathrm{U}}_{\nu\nu'}(\cdot,l)$, in a
$\{ d^2 p(\nu) \} \times (Md^2)$
matrix leads to the expression of:
\begin{equation}
\label{Gnu}
 \bfG^{\mathrm{U}}_{\nu\nu'}=: \left[ \bfG^{\mathrm{U}}_{\nu\nu'}(\cdot,1), \ldots, \bfG^{\mathrm{U}}_{\nu\nu'}(\cdot,M) \right]
\end{equation}
where
\begin{align*}
\bfG^{\mathrm{U}}_{\nu\nu'}(\cdot,l) = :\sum_{h=-\infty}^{\infty}\mathbb{E}\left[\bfOmega^{-1}(\nu) \bfX_n(\nu) \bfepsilon_{(n-h)s+\nu'-l}^{\top}
                             \otimes \bfepsilon_{ns+\nu}\bfepsilon_{(n-h)s+\nu'}^\top\right]\quad \text{for}\quad l=1,\ldots,M.
\end{align*}
We also have:
\begin{eqnarray}\nonumber
  N^{1/2} \left( \hat{\bfbeta}^\top - \bfbeta^\top,\bfc_{\bfepsilon}^\top \right)^\top
  &\cd& \mathcal{N}_{sd^2(p(\nu)+M)}\left(\bfzero, \bfXi^{\mathrm{U}}=:\left(\begin{array}{cc}
                           \bfTheta & \bfG^{\mathrm{U}}\\
                           \left(\bfG^{\mathrm{U}}\right)^{\top} & \bfV
                           \end{array}
                     \right) \right),
\end{eqnarray}
where the asymptotic covariance matrix $\bfTheta$ is given by \eqref{th1d}, $\bfV = \left[\bfV_{\nu\nu'}(M)\right]_{1\leq \nu,\nu'\leq s}$ and $\bfG^{\mathrm{U}} = \left[\bfG^{\mathrm{U}}_{\nu\nu'}\right]_{1\leq \nu,\nu'\leq s}$.
The exponent $\mathrm{U}$  stands for Unconstrained parameter.
\end{proposition}
The proof of Proposition~\ref{loijointeUn} is postponed to Section Appendix.
\begin{remark} \label{VuseMAinfty}
Using the infinite moving average representation~(\ref{MAinf}), for $l=1,\dots,s$ and $m=1,\dots,p(\nu)$ observe that:
\begin{align*}
\sum_{h=-\infty}^{\infty}
\mathbb{E}\left[ \bfY_{ns+\nu-m} \bfepsilon_{(n-h)s+\nu'-l}^{\top}
                             \otimes \bfepsilon_{ns+\nu}\bfepsilon_{(n-h)s+\nu'}^\top\right]
       &=\sum_{k=0}^{\infty}\left( \bfF_k(\nu-m)\otimes \bfI_d\right)\sum_{h=-\infty}^{\infty} \mathbb{E}\left[ \bfepsilon_{ns+\nu-k-m} \bfepsilon_{(n-h)s+\nu'-l}^{\top}
                             \otimes \bfepsilon_{ns+\nu}\bfepsilon_{(n-h)s+\nu'}^\top\right], \\
      &= \sum_{k=0}^{\infty}\left( \bfF_k(\nu-m)\otimes \bfI_d\right)\bfV_{\nu\nu'}(m+k,l),
\end{align*}
where $\bfF_0(\nu) = \bfI_d$ and $\bfF_k(\nu) = \bfzero$ for $k < 0$, $\nu=1,\ldots,s$. It follows that:
\begin{align*}
\bfG^{\mathrm{U}}_{\nu\nu'}(\cdot,l) &=\left[ \bfOmega^{-1}(\nu)\otimes \bfI_d\right)\sum_{k=0}^{\infty} \left(
                    \begin{array}{c}
                       \left(\bfF_{k}(\nu-1)\otimes \bfI_d\right)\bfV_{\nu\nu'}(1+k,l) \\
                       \vdots \\
                       \left(\bfF_{k}(\nu-p(\nu))\otimes \bfI_d\right)\bfV_{\nu\nu'}(p(\nu)+k,l) \\
                    \end{array}
                 \right),
\end{align*}
where $\bfF_0(\nu) = \bfI_d$
and $\bfF_k(\nu) = \bfzero$ for $k < 0$, $\nu=1,\ldots,s$. For $\ell,\ell'=1,\dots,p(\nu)$ observe also that:
\begin{align*}
\sum_{h=-\infty}^{\infty}
\mathbb{E}&\left[ \bfOmega^{-1}(\nu)\bfY_{ns+\nu-\ell} \bfY_{(n-h)s+\nu'-\ell'}^{\top} \bfOmega^{-1}(\nu')
                             \otimes \bfepsilon_{ns+\nu}\bfepsilon_{(n-h)s+\nu'}^\top\right]\\
       &=\sum_{i,j=0}^{\infty}\sum_{h=-\infty}^{\infty} \mathbb{E}\left[\bfOmega^{-1}(\nu)\bfF_i(\nu-\ell) \bfepsilon_{ns+\nu-i-\ell} \bfepsilon_{(n-h)s+\nu'-j-\ell'}^{\top} \bfF^\top_j(\nu'-\ell')\bfOmega^{-1}(\nu')
                             \otimes \bfepsilon_{ns+\nu}\bfepsilon_{(n-h)s+\nu'}^\top\right], \\
      &= \sum_{i,j=0}^{\infty}\left( \bfOmega^{-1}(\nu)\bfF_i(\nu-\ell)\otimes \bfI_d\right)\bfV_{\nu\nu'}(i+\ell,j+\ell')\left( \bfF^\top_j(\nu'-\ell')\bfOmega^{-1}(\nu')\otimes \bfI_d\right):=\bfTheta_{\nu\nu'}(\ell,\ell').
\end{align*}
It follows that:
$\bfTheta_{\nu\nu'} =\left[ \bfTheta_{\nu\nu'}(\ell,\ell')\right]_{1\leq\ell,\ell'\leq p(\nu)}$
It is clear that the existence of the matrices $\bfXi_{\nu\nu}^{\mathrm{U}}$ and $\bfXi^{\mathrm{U}}$ is ensured by the existence of
$\bfV_{\nu\nu'}(l',l)$, which is the consequence of Assumption \textbf{(A3)} and the \cite{D68} inequality.
\end{remark}

\begin{remark}

\begin{itemize}
\item[i)] When $d=1$ and if the moving average order is null, we retrieve the result on weak PARMA obtained by \cite{BMIA23}. \item[ii)] When $s=1$, we retrieve the well-known result on weak VAR obtained by \cite{FR07}. Finally, if the moving average order is null, we obtain the results on weak VARMA of~\cite{BM11} and  those of~\cite{frz05} when $d=1$.
\end{itemize}
\end{remark}

\begin{remark}
\label{remUn}
In the standard strong PVAR case, i.e. when {\bf (A2)} is replaced by the assumption that $\left(\boldsymbol{\epsilon}^{\ast}_n\right)_{n\in \mathbb{Z}}$ is  an iid sequence, we
have $$\bfPsi_{\nu\nu} =\bfOmega(\nu) \otimes \bfSigma_{\bfepsilon}(\nu).$$
In view of Remark \ref{VuseMAinfty}, independence of the $\bfepsilon_{ns+\nu}$'s implies that only the terms for $h=0$ and $\nu=\nu'$ are non-zero in $\bfG^{\mathrm{U}}_{\nu\nu'}(\cdot,l)$, for $l=1,\dots,M$.
Thus the asymptotic covariance matrices $\bfTheta_{\nu\nu}$, $\bfV_{\nu\nu}$ and $\bfG^{\mathrm{U}}_{\nu\nu}$ are reduced respectively as
\begin{align*}
\bfTheta_{\nu\nu}&:= \bfOmega^{-1}(\nu) \otimes \bfSigma_{\bfepsilon}(\nu),
\\
\bfV_{\nu\nu} &:= \left(
    \begin{array}{ccccc}
    \bfSigma_{\bfepsilon}(\nu-1) & \bfzero        & \bfzero & \hdots & \bfzero \\
    \bfzero         & \bfSigma_{\bfepsilon}(\nu-2) & \bfzero & \hdots & \bfzero \\
    \vdots          &                &          & \ddots  &  \vdots   \\
    \bfzero         & \bfzero        & \bfzero & \hdots &\bfSigma_{\bfepsilon}(\nu-M)
    \end{array}
    \right) \otimes \bfSigma_{\bfepsilon}(\nu),
\\ \bfG^{\mathrm{U}}_{\nu\nu}&=\left[\bfOmega^{-1}(\nu)\mathbb{E}\left(\bfX_n(\nu)\bfepsilon^\top_{ns+\nu-1}\right)\otimes \bfSigma_{\bfepsilon}(\nu),\dots,\bfOmega^{-1}(\nu)\mathbb{E}\left(\bfX_n(\nu)\bfepsilon^\top_{ns+\nu-M}\right)\otimes \bfSigma_{\bfepsilon}(\nu) \right].
\end{align*}
Let
\begin{equation}
\label{Anu}
 \bfA_M(\nu)=\left[ \bfA_{\nu}(1), \ldots, \bfA_{\nu}(M) \right]
\end{equation}
be the $\{ d p(\nu) \} \times (Md)$
matrix where   the $\{ d p(\nu) \} \times d$ matrix
$\bfA_{\nu}(l)$ is given by:
\begin{align*}
\bfA_{\nu}(l) = \mathbb{E}\left[ \bfX_n(\nu) \bfepsilon_{ns+\nu-l}^{\top}\right]\quad \text{for}\quad l=1,\ldots,M.
\end{align*}
Therefore, we obtain that:
\begin{align*}
\bfG^{\mathrm{U}}_{\nu\nu}(\cdot,l) &=\left[ \bfOmega^{-1}(\nu)\bfA_{\nu}(l)\right]\otimes \bfSigma_{\bfepsilon}(\nu)\quad\text{where}\quad
\bfA_{\nu}(l)=\left(
                    \begin{array}{c}
                       \bfF_{l-1}(\nu-1)\bfSigma_{\bfepsilon}(\nu-l) \\
                       \vdots \\
                       \bfF_{l-p(\nu)}(\nu-p(\nu))\bfSigma_{\bfepsilon}(\nu-l) \\
                    \end{array}
                 \right).
\end{align*}
Then we deduce that: $\bfG^{\mathrm{U}}_{\nu\nu} = \left[\bfOmega^{-1}(\nu)\otimes\bfI_d\right]\left[\bfA_{M}(\nu)\otimes \bfSigma_{\bfepsilon}(\nu)\right]$
and we retrieve the results given in \cite{UD09}.
\end{remark}
We now investigate a similar tractable expression when the parameters satisfy the linear constraints~(\ref{constraints}). In view of \eqref{ecartxi}, we deduce the following relation:
\begin{equation}\label{ecartxiOmega}
  N^{1/2} \{  \hat{\bfxi}(\nu) - \bfxi(\nu) \}  
  = \left[ \bfR^{\top}(\nu) \{ \bfOmega(\nu) \otimes \bfSigma^{-1}_{\bfepsilon}(\nu) \} \bfR(\nu) \right]^{-1}
  \bfR^{\top}(\nu)
  \{ \bfI_{dp(\nu)}  \otimes \bfSigma^{-1}_{\bfepsilon}(\nu) \} N^{-1/2} \vec\{ \bfE(\nu) \bfX^\top(\nu)  \}+\mathrm{o}_{\mathbb{P}}(1).
\end{equation}
The next proposition gives the joint asymptotic distribution of \[\left[ \bfR^{\top}(\nu) \{ \bfOmega(\nu) \otimes \bfSigma^{-1}_{\bfepsilon}(\nu) \} \bfR(\nu) \right]^{-1}
  \bfR^{\top}(\nu)
  \{ \bfI_{dp(\nu)}  \otimes \bfSigma^{-1}_{\bfepsilon}(\nu) \} N^{-1/2} \vec\{ \bfE(\nu) \bfX^\top(\nu)  \}\]
and
$N^{1/2} \bfc_{\bfepsilon}(\nu)$.
\begin{proposition}\label{loijointeCon}
Under the assumptions of Proposition~\ref{loijointeUn} and~(\ref{constraints}), for a specified season $\nu=1,\dots,s$ we have:
\begin{eqnarray}\nonumber
  N^{1/2} \left( \hat{\bfxi}^\top(\nu) - \bfxi^\top(\nu),\bfc_{\bfepsilon}^\top(\nu) \right)^\top
  &\cd& \mathcal{N}_{K(\nu)+Md^2}\left(\bfzero, \bfXi_{\nu\nu}^{\mathrm{R}}\right),
\end{eqnarray}
 with
 \begin{align}\label{specWnCon}
\bfXi_{\nu\nu}^{\mathrm{R}}=\left(\begin{array}{cc}
                           \bfTheta_{\nu\nu}^{\bfxi} & \bfG^{\mathrm{R}}_{\nu\nu}\\
                           \left(\bfG_{\nu\nu}^{\mathrm{R}}\right)^{\top} & \bfV_{\nu\nu}
                           \end{array}
                     \right) =\sum_{h=-\infty}^{\infty}\mathbb{E}\left[ \bfW_{ns+\nu}^{\mathrm{R}}\left(\bfW_{(n-h)s+\nu}^{\mathrm{R}}\right)^{\top}\right],
\end{align}
where $\bfTheta_{\nu\nu}^{\bfxi}$, $\bfV_{\nu\nu}$, $\bfG^{\mathrm{R}}_{\nu\nu}$ are given respectively by \eqref{NAcont}, \eqref{VnuM}, \eqref{GnuCon} and 
\begin{equation}\label{WnCon}
\bfW_{ns+\nu}^{\mathrm{R}}=\begin{pmatrix}\bfW^{\mathrm{R}}_{1,ns+\nu}\vspace{0.2cm}\\
\bfW_{2,ns+\nu}
\end{pmatrix}=\begin{pmatrix} \bfH_{\bfepsilon}(\nu)\{\bfX_n(\nu) \otimes \bfepsilon_{ns+\nu}\}\vspace{0.2cm}\\
\left(\bfepsilon_{ns+\nu-1}^\top,\dots,\bfepsilon_{ns+\nu-M}^\top\right)^{\top}\otimes\bfepsilon_{ns+\nu}
\end{pmatrix}.
\end{equation}
with $\bfH_{\bfepsilon}(\nu)=
                             \left[ \bfR^{\top}(\nu) \{ \bfOmega(\nu) \otimes \bfSigma^{-1}_{\bfepsilon}(\nu) \} \bfR(\nu) \right]^{-1}
  \bfR^{\top}(\nu)
  \{ \bfI_{dp(\nu)}  \otimes \bfSigma^{-1}_{\bfepsilon}(\nu) \}.$
Finally, collecting the $\{ K(\nu) \} \times d^2$ matrix
$\bfG^{\mathrm{R}}_{\nu\nu'}(\cdot,l)$, for $l=1,\ldots,M$ and $\nu,\nu'=1,\dots,s$, in a
$\{ K(\nu) \} \times (Md^2)$
matrix leads to the expression of:
\begin{equation}
\label{GnuCon}
 \bfG^{\mathrm{R}}_{\nu\nu'}=: \left[ \bfG^{\mathrm{R}}_{\nu\nu'}(\cdot,1), \ldots, \bfG^{\mathrm{R}}_{\nu\nu'}(\cdot,M) \right].
\end{equation}
where
\begin{align*}
\bfG^{\mathrm{R}}_{\nu\nu'}(\cdot,l) &= \bfH_{\bfepsilon}(\nu) \sum_{h=-\infty}^{\infty}\mathbb{E}\left[\bfX_n(\nu) \bfepsilon_{(n-h)s+\nu'-l}^{\top}
                             \otimes \bfepsilon_{ns+\nu}\bfepsilon_{(n-h)s+\nu'}^\top\right].
\end{align*}
We also have:
\begin{eqnarray}\nonumber
  N^{1/2} \left( \hat{\bfxi}^\top - \bfxi^\top,\bfc_{\bfepsilon}^\top \right)^\top
  &\cd& \mathcal{N}_{s(K(\nu)+Md^2)}\left(\bfzero, \bfXi^{\mathrm{R}}=:\left(\begin{array}{cc}
                           \bfTheta^{\bfxi} & \bfG^{\mathrm{R}}\\
                           \left(\bfG^{\mathrm{R}}\right)^{\top} & \bfV
                           \end{array}
                     \right) \right),
\end{eqnarray}
where the asymptotic covariance matrix $\bfTheta^{\bfxi}$ is given by~\eqref{th1dCon}, $\bfV = \left[\bfV_{\nu\nu'}(M)\right]_{1\leq \nu,\nu'\leq s}$ and\\%
$\bfG^{\mathrm{R}} = \left[\bfG^{\mathrm{R}}_{\nu\nu'}\right]_{1\leq \nu,\nu'\leq s}$.
The exponent $\mathrm{R}$  stands for restricted parameter.
\end{proposition}
The proof of Proposition~\ref{loijointeCon} follows using the same kind of arguments as Proposition~\ref{loijointeUn} and it is therefore omitted.
\begin{remark}
\label{remCon} In view of Remark \ref{remUn}, in the standard strong  PVAR case, the asymptotic covariance matrices $\bfTheta_{\nu\nu}$, $\bfV_{\nu\nu}$ and $\bfG^{\mathrm{R}}_{\nu\nu}$ are reduced respectively as:
\begin{align*}
\bfTheta_{\nu\nu}^{\bfxi}&:= \left[ \bfR^{\top}(\nu) \{ \bfOmega(\nu) \otimes \bfSigma^{-1}_{\bfepsilon}(\nu) \} \bfR(\nu) \right]^{-1},
\\
\bfV_{\nu\nu} &:= \left(
    \begin{array}{ccccc}
    \bfSigma_{\bfepsilon}(\nu-1) & \bfzero        & \bfzero & \hdots & \bfzero \\
    \bfzero         & \bfSigma_{\bfepsilon}(\nu-2) & \bfzero & \hdots & \bfzero \\
    \vdots          &                &          & \ddots  &  \vdots   \\
    \bfzero         & \bfzero        & \bfzero & \hdots &\bfSigma_{\bfepsilon}(\nu-M)
    \end{array}
    \right) \otimes \bfSigma_{\bfepsilon}(\nu),
\\ \bfG^{\mathrm{R}}_{\nu\nu}&=\bfH_{\bfepsilon}(\nu)\left[\bfA_{M}(\nu) \otimes \bfSigma_{\bfepsilon}(\nu)\right],
=
\left[ \bfR^{\top}(\nu) \{ \bfOmega(\nu) \otimes \bfSigma^{-1}_{\bfepsilon}(\nu) \} \bfR(\nu) \right]^{-1}
  \bfR^{\top}(\nu)
  \left(\bfA_{M}(\nu)  \otimes \bfI_{d} \right).
\end{align*}
\end{remark}

\noindent The following theorem, which is an extension of the results given in~\cite{FR07}, provides the  limit distribution of the residual autocovariances and autocorrelations of weak PVAR models for a specified season $\nu$. We also  provide the  limit distribution of the global residual autocovariances and autocorrelations of weak PVAR models across all the seasons.
\begin{theorem}
\label{loi_res_gamUn} Under the assumptions of Proposition~\ref{loijointeUn} we have
\begin{equation}
\label{chatgammaU}
  N^{1/2}\bfc_{\hat{\bfepsilon}}(\nu) \cd \mathcal{N}_{d^2M}\left(\bfzero, \bfDelta^{\mathrm{U}}_{\nu\nu}\right),
\end{equation}
where
\begin{align*}
\bfDelta^{\mathrm{U}}_{\nu\nu}=\bfV_{\nu\nu}+\bfUpsilon_M(\nu)\bfTheta_{\nu\nu}\bfUpsilon^{\top}_M(\nu)+ \bfUpsilon_M(\nu)\bfG^{\mathrm{U}}_{\nu\nu}+\left(\bfG^{\mathrm{U}}_{\nu\nu}\right)^{\top}\bfUpsilon^{\top}_M(\nu)
\text{ when }p(\nu)>0\quad\text{and }
\bfDelta^{\mathrm{U}}_{\nu\nu}&=\bfV_{\nu\nu}\quad\text{ when }\quad p(\nu)=0.
\end{align*}
Then:
\begin{equation}
\label{chatgammaglobalU}
  N^{1/2}\bfc_{\hat{\bfepsilon}} \cd \mathcal{N}_{sd^2M}\left(\bfzero, \bfDelta^{\mathrm{U}}\right), 
\end{equation}
where the asymptotic covariance matrix $\bfDelta^{\mathrm{U}}$ is a block matrix, with the asymptotic variances given by $\bfDelta^{\mathrm{U}}_{\nu\nu}$, for  $\nu=1,\dots,s$  and the asymptotic covariances given by
$$\bfDelta^{\mathrm{U}}_{\nu\nu'}=\bfV_{\nu\nu'}+\bfUpsilon_M(\nu)\bfTheta_{\nu\nu'}\bfUpsilon^{\top}_M(\nu')+ \bfUpsilon_M(\nu)\bfG^{\mathrm{U}}_{\nu\nu'}+\left(\bfG^{\mathrm{U}}_{\nu\nu'}\right)^{\top}\bfUpsilon^{\top}_M(\nu'),\quad \nu\neq\nu'.$$
We also have
\begin{equation}
\label{chatrhoU}
  N^{1/2}\bfr_{\hat{\bfepsilon}}(\nu) \cd \mathcal{N}_{d^2M}\left(\bfzero, \bfnabla^{\mathrm{U}}_{\nu\nu}\right),
\end{equation}
where
\begin{eqnarray*}
\bfnabla^{\mathrm{U}}_{\nu\nu} = \left(\diag\left[\left\{ \bfD^{\ast}_{{\bfepsilon}}(\nu-1),\dots,\bfD^{\ast}_{{\bfepsilon}}(\nu-M)\right\}  \otimes
                                    \bfD^{\ast}_{{\bfepsilon}}(\nu) \right]\right)^{-1}\bfDelta^{\mathrm{U}}_{\nu\nu}
\left(\diag\left[\left\{ \bfD^{\ast}_{{\bfepsilon}}(\nu-1),\dots,\bfD^{\ast}_{{\bfepsilon}}(\nu-M)\right\}  \otimes
                                    \bfD^{\ast}_{{\bfepsilon}}(\nu) \right]\right)^{-1},
\end{eqnarray*}
with $\bfD^{\ast}_{{\bfepsilon}}(\nu) =
\diag\left( \bfGamma_{{\bfepsilon},11}^{1/2}(0; \nu),\ldots,\bfGamma_{{\bfepsilon},dd}^{1/2}(0; \nu) \right)$. Then:
\begin{equation}
\label{chatrhoglobalU}
  N^{1/2}\bfr_{\hat{\bfepsilon}} \cd \mathcal{N}_{sd^2M}\left(\bfzero, \bfnabla^{\mathrm{U}}\right), 
\end{equation}
where the asymptotic covariance matrix $\bfnabla^{\mathrm{U}}= [\bfnabla^{\mathrm{U}}_{\nu\nu'}]_{1\leq \nu,\nu'\leq s}$ is a block matrix, with blocs given by
\begin{eqnarray*}
\bfnabla^{\mathrm{U}}_{\nu\nu'} = \left(\diag\left[\left\{ \bfD^{\ast}_{{\bfepsilon}}(\nu-1),\dots,\bfD^{\ast}_{{\bfepsilon}}(\nu-M)\right\}  \otimes
                                    \bfD^{\ast}_{{\bfepsilon}}(\nu) \right]\right)^{-1}\bfDelta^{\mathrm{U}}_{\nu\nu'}
\left(\diag\left[\left\{ \bfD^{\ast}_{{\bfepsilon}}(\nu'-1),\dots,\bfD^{\ast}_{{\bfepsilon}}(\nu'-M)\right\}  \otimes
                                    \bfD^{\ast}_{{\bfepsilon}}(\nu') \right]\right)^{-1}.
\end{eqnarray*}
\end{theorem}

The proof of Theorem~\ref{loi_res_gamUn} is postponed to Section Appendix.
\begin{remark} Observe that, in view of \eqref{Anu}, the matrix $\bfUpsilon_M(\nu)$  can be rewritten as:
\begin{equation}\label{Ups_M_MAinfty}\bfUpsilon_M(\nu)=-\bfA_{M}^\top(\nu)\otimes \bfI_d.
\end{equation}
When $s=1$, we retrieve the well-known result on weak VAR obtained by \cite{FR07}. Finally, if the moving average order is null in \cite{BM11} and \cite{frz05}, we obtain the results on weak VARMA of \cite{BM11} and  those of \cite{frz05} when $d=1$.
\end{remark}
\begin{remark}
\label{remrhoUn}
In view of Remark \ref{remUn} and \eqref{Ups_M_MAinfty}, in the standard strong PVAR case, the asymptotic covariance matrices $\bfDelta^{\mathrm{U}}_{\nu\nu}$ and $\bfnabla^{\mathrm{U}}_{\nu\nu}$ are reduced respectively as
\begin{align*}
\bfDelta^{\mathrm{U}}_{\mathrm{S}}(\nu,\nu)=:\bfV_{\nu\nu} -\left[\bfA^\top_{M}(\nu)\bfOmega^{-1}(\nu) \bfA_{M}(\nu)\right]\otimes \bfSigma_{\bfepsilon}(\nu)
\end{align*}
and
\begin{eqnarray*}
\bfnabla^{\mathrm{U}}_{\mathrm{S}}(\nu,\nu) =: \left(\diag\left[\left\{ \bfD^{\ast}_{{\bfepsilon}}(\nu-1),\dots,\bfD^{\ast}_{{\bfepsilon}}(\nu-M)\right\}  \otimes
                                    \bfD^{\ast}_{{\bfepsilon}}(\nu) \right]\right)^{-1} \bfDelta^{\mathrm{U}}_{\mathrm{S}}(\nu,\nu)
\left(\diag\left[\left\{ \bfD^{\ast}_{{\bfepsilon}}(\nu-1),\dots,\bfD^{\ast}_{{\bfepsilon}}(\nu-M)\right\}  \otimes
                                    \bfD^{\ast}_{{\bfepsilon}}(\nu) \right]\right)^{-1},
\end{eqnarray*}
and we retrieve the results given in \cite{UD09}. When $s=1$, we also retrieve the results obtained by \cite{CD08} and \cite{Ho1980}.

Independence of the $\bfepsilon_{ns+\nu}$'s implies that only the terms for  $\nu=\nu'$ are non-zero in $\bfDelta^{\mathrm{U}}$. Thus the matrices $\bfDelta^{\mathrm{U}}$ and $\bfnabla^{\mathrm{U}}$ are reduced to blocs diagonal matrices denoted $\bfDelta_{\mathrm{S}}^{\mathrm{U}}$ and $\bfnabla_{\mathrm{S}}^{\mathrm{U}}$, with blocs given respectively  $\bfDelta^{\mathrm{U}}_{\mathrm{S}}(\nu,\nu)$ and $\bfnabla^{\mathrm{U}}_{\mathrm{S}}(\nu,\nu)$.
\end{remark}
We now investigate a similar tractable expression when the parameters satisfy the linear constraints~(\ref{constraints}). In this case,
similar calculations give:
\begin{align}
\label{Ups_Con}
  \frac{\partial \bfc_{\dot{\bfepsilon}}(\nu)}{\partial \bfxi^{\top}(\nu)} \cp \bfUpsilon_M(\nu) \bfR(\nu),
\end{align}
and using the same kind of arguments we obtain the following theorem.
\begin{theorem}
\label{loi_res_gamCon}
Under the assumptions of Proposition~\ref{loijointeCon} we have
\begin{equation}
\label{chatgammaR}
  N^{1/2}\bfc_{\hat{\bfepsilon}}(\nu) \cd \mathcal{N}_{d^2M}\left(\bfzero, \bfDelta^{\mathrm{R}}_{\nu\nu}\right),
\end{equation}
where
\begin{align*}
\bfDelta^{\mathrm{R}}_{\nu\nu}&=\bfV_{\nu\nu}+\bfUpsilon_M(\nu)\bfR(\nu)\bfTheta^{\bfxi}_{\nu\nu} \left(\bfUpsilon_M(\nu)\bfR(\nu)\right)^{\top}+\bfUpsilon_M(\nu)\bfR(\nu)\bfG^{\mathrm{R}}_{\nu\nu} +\left(\bfG^{\mathrm{R}}_{\nu\nu}\right)^{\top}\left(\bfUpsilon_M(\nu)\bfR(\nu)\right)^{\top}
\text{ when }p(\nu)>0\;\text{and}\\
\bfDelta^{\mathrm{R}}_{\nu\nu}&=\bfV_{\nu\nu}\;\text{ when}\; p(\nu)=0.
\end{align*}
Then:
\begin{equation}
\label{chatgammaglobalR}
  N^{1/2}\bfc_{\hat{\bfepsilon}} \cd \mathcal{N}_{sd^2M}\left(\bfzero, \bfDelta^{\mathrm{R}}\right), 
\end{equation}
where the asymptotic covariance matrix $\bfDelta^{\mathrm{R}}$ is a block matrix, with the asymptotic variances given by $\bfDelta^{\mathrm{R}}_{\nu\nu}$, for  $\nu=1,\dots,s$  and the asymptotic covariances given by
\begin{eqnarray*}
\bfDelta^{\mathrm{R}}_{\nu\nu'}=\bfV_{\nu\nu'}+\bfUpsilon_M(\nu)\bfR(\nu)\bfTheta^{\bfxi}_{\nu\nu'} \left(\bfR(\nu)\bfUpsilon_M(\nu')\bfR(\nu')\right)^{\top}+\bfUpsilon_M(\nu)\bfR(\nu)\bfG^{\mathrm{R}}_{\nu\nu'}
+ \left(\bfG^{\mathrm{R}}_{\nu\nu'}\right)^{\top}\left(\bfUpsilon_M(\nu')\bfR(\nu')\right)^{\top},\quad \nu\neq\nu'.
\end{eqnarray*}
We also have
\begin{equation}
\label{chatrhoR}
  N^{1/2}\bfr_{\hat{\bfepsilon}}(\nu) \cd \mathcal{N}_{d^2M}\left(\bfzero, \bfnabla^{\mathrm{R}}_{\nu\nu}\right),
\end{equation}
where
\begin{eqnarray*}
\bfnabla^{\mathrm{R}}_{\nu\nu} = \left(\diag\left[\left\{ \bfD^{\ast}_{{\bfepsilon}}(\nu-1),\dots,\bfD^{\ast}_{{\bfepsilon}}(\nu-M)\right\}  \otimes
                                    \bfD^{\ast}_{{\bfepsilon}}(\nu) \right]\right)^{-1}\bfDelta^{\mathrm{R}}_{\nu\nu}
\left(\diag\left[\left\{ \bfD^{\ast}_{{\bfepsilon}}(\nu-1),\dots,\bfD^{\ast}_{{\bfepsilon}}(\nu-M)\right\}  \otimes
                                    \bfD^{\ast}_{{\bfepsilon}}(\nu) \right]\right)^{-1}.
\end{eqnarray*}
Then:
\begin{equation}
\label{chatrhoglobalR}
  N^{1/2}\bfr_{\hat{\bfepsilon}} \cd \mathcal{N}_{sd^2M}\left(\bfzero, \bfnabla^{\mathrm{R}}\right), 
\end{equation}
where the asymptotic covariance matrix $\bfnabla^{\mathrm{R}}= [\bfnabla^{\mathrm{R}}_{\nu\nu'}]_{1\leq \nu,\nu'\leq s}$ is a block matrix, with blocs given by
\begin{eqnarray*}
\bfnabla^{\mathrm{R}}_{\nu\nu'} = \left(\diag\left[\left\{ \bfD^{\ast}_{{\bfepsilon}}(\nu-1),\dots,\bfD^{\ast}_{{\bfepsilon}}(\nu-M)\right\}  \otimes
                                    \bfD^{\ast}_{{\bfepsilon}}(\nu) \right]\right)^{-1}\bfDelta^{\mathrm{R}}_{\nu\nu'}
\left(\diag\left[\left\{ \bfD^{\ast}_{{\bfepsilon}}(\nu'-1),\dots,\bfD^{\ast}_{{\bfepsilon}}(\nu'-M)\right\}  \otimes
                                    \bfD^{\ast}_{{\bfepsilon}}(\nu') \right]\right)^{-1}.
\end{eqnarray*}
\end{theorem}
The proof of Theorem~\ref{loi_res_gamCon} follows using the same kind of arguments as Theorem~\ref{loi_res_gamUn} and it is therefore omitted.
\begin{remark}
\label{remrhoCon}
In view of Remark \ref{remCon} and \eqref{Ups_M_MAinfty}, in the standard strong PVAR case, the asymptotic covariance matrices $\bfDelta^{\mathrm{R}}_{\nu\nu}$ and $\bfnabla^{\mathrm{R}}_{\nu\nu}$ are reduced respectively as
\begin{eqnarray*}
\bfDelta^{\mathrm{R}}_{\mathrm{S}}(\nu,\nu)&=:\bfV_{\nu\nu}-
  \left(\bfA^{\top}_{M}(\nu)  \otimes \bfI_{d} \right)\bfR(\nu)\left[ \bfR^{\top}(\nu) \{ \bfOmega(\nu) \otimes \bfSigma^{-1}_{\bfepsilon}(\nu) \} \bfR(\nu) \right]^{-1}
  \bfR^{\top}(\nu)
  \left(\bfA_{M}(\nu)  \otimes \bfI_{d} \right)
\end{eqnarray*}
and
\begin{eqnarray*}
\bfnabla^{\mathrm{R}}_{\mathrm{S}}(\nu,\nu) =: \left(\diag\left[\left\{ \bfD^{\ast}_{{\bfepsilon}}(\nu-1),\dots,\bfD^{\ast}_{{\bfepsilon}}(\nu-M)\right\}  \otimes
                                    \bfD^{\ast}_{{\bfepsilon}}(\nu) \right]\right)^{-1}\bfDelta^{\mathrm{R}}_{\mathrm{S}}(\nu,\nu)
\left(\diag\left[\left\{ \bfD^{\ast}_{{\bfepsilon}}(\nu-1),\dots,\bfD^{\ast}_{{\bfepsilon}}(\nu-M)\right\}  \otimes
                                    \bfD^{\ast}_{{\bfepsilon}}(\nu) \right]\right)^{-1},
\end{eqnarray*}
and we retrieve the results given in \cite{UD09} and \cite{DL13}.

Independence of the $\bfepsilon_{ns+\nu}$'s implies that only the terms for  $\nu=\nu'$ are non-zero in $\bfDelta^{\mathrm{R}}$. Thus the matrices $\bfDelta^{\mathrm{R}}$ and $\bfnabla^{\mathrm{R}}$ are reduced to blocs diagonal matrices denoted $\bfDelta_{\mathrm{S}}^{\mathrm{R}}$ and $\bfnabla_{\mathrm{S}}^{\mathrm{R}}$, with blocs given respectively  $\bfDelta^{\mathrm{R}}_{\mathrm{S}}(\nu,\nu)$ and $\bfnabla^{\mathrm{R}}_{\mathrm{S}}(\nu,\nu)$.
\end{remark}

\subsection{Modified version of the portmanteau test}\label{result2}
\noindent Based on the residual empirical autocorrelations,  \cite{bp70} have proposed a goodness-of-fit test, the so-called portmanteau test, for strong ARMA models. A modification of their test has been proposed by \cite{lb} which is nowadays one of the most popular diagnostic checking tools in ARMA modeling of time series. To test simultaneously whether all residual autocorrelations at lags $1,\dots,M$ of an univariate periodic ARMA model are equal to zero for a specified period $\nu$, the portmanteau of~\cite{bp70,lb} can be adapted as proposed by~\cite{McL1994,McL1994corrig} (see also~\cite{hipelMcleod1994}).

Being inspired by the univariate portmanteau statistics defined in~\cite{bp70} and~\cite{lb},~\cite{C1974} and~\cite{Ho1980}) have introduced  the  multivariate versions of the  portmanteau statistics. \cite{H81} gave several equivalent forms of this statistic.  Basic forms are:
\begin{eqnarray}\label{bpmulti}
{{\bfQ}}_M^{\textsc{c}}=N\sum_{h=1}^M\tr\left(\bfC_{\hat{\bfepsilon}}^\top(h)\bfC_{\hat{\bfepsilon}}^{-1}(0)
\bfC_{\hat{\bfepsilon}}(h)\bfC_{\hat{\bfepsilon}}^{-1}(0)\right)\text{ and
}{{\bfQ}}_M^{\textsc{h}}=\sum_{h=1}^M\frac{N^2}{(N-h)}\tr\left(\bfC_{\hat{\bfepsilon}}^\top(h)\bfC_{\hat{\bfepsilon}}^{-1}(0)
\bfC_{\hat{\bfepsilon}}(h)\bfC_{\hat{\bfepsilon}}^{-1}(0)\right),
 \end{eqnarray}
where $\bfC_{\hat{\bfepsilon}}(h)$ is the residual autocovariances matrices function of  a multivariate process $\bfY$.

To test simultaneously whether all residual autocorrelations at lags $1,\dots,M$ of a PVAR model are equal to zero for a specified period $\nu$,  the portmanteau test \eqref{bpmulti} can be adapted.  Based on the above results, for instance,   \cite{UD09,DL13} suggested the following portmanteau statistics defined by
\begin{align}
\nonumber
\bfQ_M(\nu) &= N\sum_{l=1}^{M}
              \tr\left(\bfC_{\hat{\bfepsilon}}^{\top}(l; \nu)
                   \bfC^{-1}_{\hat{\bfepsilon}}(0; \nu)
                   \bfC_{\hat{\bfepsilon}}(l; \nu)
                   \bfC^{-1}_{\hat{\bfepsilon}}(0; \nu-l) \right),
\\\nonumber &=   N\sum_{l=1}^{M}
             \bfc_{\hat{\bfepsilon}}^{\top}(l;\nu)
             \left\{
                    \bfC^{-1}_{\hat{\bfepsilon}}(0; \nu-l) \otimes \bfC^{-1}_{\hat{\bfepsilon}}(0; \nu)
             \right\}
              \bfc_{\hat{\bfepsilon}}(l;\nu),
\\ \nonumber &=   N\bfc_{\hat{\bfepsilon}}^{\top}(\nu)
             \left\{\left(\diag\left[\bfC_{\hat{\bfepsilon}}(0; \nu-1),\dots, \bfC_{\hat{\bfepsilon}}(0; \nu-M)\right]\right) \otimes \bfC_{\hat{\bfepsilon}}(0; \nu)
             \right\}^{-1}
              \bfc_{\hat{\bfepsilon}}(\nu),
\\&=   N\bfr_{\hat{\bfepsilon}}^{\top}(\nu)
             \left\{\left(\diag\left[\bfR_{\hat{\bfepsilon}}(0; \nu-1),\dots, \bfR_{\hat{\bfepsilon}}(0; \nu-M)\right]\right) \otimes \bfR_{\hat{\bfepsilon}}(0; \nu)
             \right\}^{-1}
              \bfr_{\hat{\bfepsilon}}(\nu),
                   \label{QMnu}
\end{align}
where $\bfR_{\hat{\bfepsilon}}(0; \nu) =  \bfD_{\hat{\bfepsilon}}^{-1}(\nu)\bfC_{\hat{\bfepsilon}}(0; \nu)
                                    \bfD_{\hat{\bfepsilon}}^{-1}(\nu)$.
The last equalities are obtained from the following elementary identities on matrix calculus
\[
   \tr(\bfA^{\top} \bfB \bfC \bfD^{\top})= \{ \vec(\bfA) \}^{\top} \{\bfD \otimes \bfB\} \{ \vec(\bfC) \}, \quad (\bfA\otimes \bfB)(\bfC\otimes \bfD)
 =\bfA\bfC\otimes \bfB\bfD
\]
where
$\bfA$,
$\bfB$,
$\bfC$
and
$\bfD$
are any matrices
for which the above product is defined.

As discussed in \cite{McL1994,McL1994corrig}, the \cite{lb}  correction factor
defined by
$$\varrho(l,\nu,N,s)= \left\lbrace
\begin{array}{cc}
(N+2)/(N-l/s) & \mbox{if} \quad l\equiv 0\mod s\\
N/\left(N-\lfloor (l-\nu+s)/s \rfloor\right) &\qquad \mbox{otherwise},
\end{array}\right.$$
is expected to improve the finite-sample properties of the test statistic~(\ref{QMnu}). This leads to the finite-sample corrected test statistic:
\begin{align}\nonumber
\bfQ_M^{\ast}(\nu) &= N \sum_{l=1}^{M} \frac{N}{N-\lfloor (l-\nu+s)/s \rfloor}
                 \tr\left(\bfC_{\hat{\bfepsilon}}^{\top}(l; \nu)
                   \bfC^{-1}_{\hat{\bfepsilon}}(0; \nu)
                   \bfC_{\hat{\bfepsilon}}(l; \nu)
                   \bfC^{-1}_{\hat{\bfepsilon}}(0; \nu-l) \right)
\\ \nonumber &=   N\bfc_{\hat{\bfepsilon}}^{\top}(\nu)\bfP^{1/2}_M(\nu,N,s)
             \left\{\left(\diag\left[\bfC_{\hat{\bfepsilon}}(0; \nu-1),\dots, \bfC_{\hat{\bfepsilon}}(0; \nu-M)\right]\right) \otimes \bfC_{\hat{\bfepsilon}}(0; \nu)
             \right\}^{-1}\bfP^{1/2}_M(\nu,N,s)
              \bfc_{\hat{\bfepsilon}}(\nu),
\\&=   N\bfr_{\hat{\bfepsilon}}^{\top}(\nu)\bfP^{1/2}_M(\nu,N,s)
             \left\{\left(\diag\left[\bfR_{\hat{\bfepsilon}}(0; \nu-1),\dots, \bfR_{\hat{\bfepsilon}}(0; \nu-M)\right]\right) \otimes \bfR_{\hat{\bfepsilon}}(0; \nu)
             \right\}^{-1}\bfP^{1/2}_M(\nu,N,s)
              \bfr_{\hat{\bfepsilon}}(\nu),
\label{QMnuast}
\end{align}
where  collecting the $d \times d$ matrix  $\bfP(l,\nu,N,s)=\diag\left[\varrho(l,\nu,N,s),\dots,\varrho(l,\nu,N,s)\right]$, in a
$(Md^2) \times (Md^2)$
matrix leads to the expression of:
$\bfP_M(\nu,N,s)=\diag\left[\left(\vec{\left\{\bfP(1,\nu,N,s)\right\}}\right)^\top,\dots,\left(\vec{\left\{\bfP(M,\nu,N,s)\right\}}\right)^\top)\right].$

Note that the lag $M$ used in \eqref{QMnu} and \eqref{QMnuast} could be chosen to be different across the seasons but in most applications it is reasonable to use the same value of $M$ for all seasons (see for instance \cite{hipelMcleod1994}).

The test statistics~(\ref{QMnu}) and~(\ref{QMnuast})
are asymptotically independent across the seasons
$\nu = 1, 2, \ldots, s$.
Consequently, global test statistics, which can be used to test the null hypothesis of model adequacy
for all seasons taken simultaneously, can be constructed by summing~(\ref{QMnu}) or~(\ref{QMnuast}) over all seasons:
\begin{eqnarray}
\label{QM}
  \bfQ_M &=& \sum_{\nu=1}^s \bfQ_M(\nu), \\
\label{QMast}
  \bfQ_M^{\ast} &=& \sum_{\nu=1}^s \bfQ_M^{\ast}(\nu).
\end{eqnarray}
The statistics \eqref{QMnu} and \eqref{QMnuast} are usually used to test the joint statistical significance of $\bfc_{\bfepsilon}(l; \nu)$, for $l=1,\ldots,M$. More formally, the null hypothesis of model adequacy is given by:
\[
H_0(\nu): \bfgamma_{\bfepsilon}(\nu) = \bfzero \quad\text{or }\quad
H_0(\nu): \bfrho_{\bfepsilon}(\nu) = \bfzero,\quad\text{for }\quad\nu=1,\dots,s,
\]
where $\bfgamma_{\bfepsilon}(\nu)$ and $\bfrho_{\bfepsilon}(\nu)$
are defined by~(\ref{gammanu}) and \eqref{rhonu} and $\bfzero$
corresponds to the $(d^2M) \times 1$ null vector. The statistics \eqref{QM} and \eqref{QMast} are usually used to test the following null  hypothesis
\[
H_0: \bfgamma_{\bfepsilon}=:\left(\bfgamma^\top_{\bfepsilon}(1),\dots,\bfgamma^\top_{\bfepsilon}(s)\right)^\top = \bfzero \quad\text{or }\quad
H_0: \bfrho_{\bfepsilon}=:\left(\bfrho^\top_{\bfepsilon}(1),\dots,\bfrho^\top_{\bfepsilon}(s)\right)^\top = \bfzero,
\]
where $\bfzero$ corresponds to the $(sd^2M) \times 1$ null vector.

From Theorems \ref{loi_res_gamUn} and \ref{loi_res_gamCon} we can deduce the following results which give  the exact  limiting distribution of the standard portmanteau statistics \eqref{QMnu}, \eqref{QMnuast}, \eqref{QM} and \eqref{QMast} under general assumptions on the innovation process of the fitted PVAR model.
\begin{theorem}
\label{theolimitdistBPUn}
Under Assumptions of Theorem \ref{loi_res_gamUn} and $H_0(\nu)$, $\nu=1,\dots,s$ or $H_0$,  the statistics $\bfQ_M(\nu)$ and $\bfQ_M^{\ast}(\nu)$ defined in \eqref{QMnu} and \eqref{QMnuast} converge in distribution, as $N\rightarrow\infty,$ to $$Z_M^{\nu}(\xi_{Md^2}^{\mathrm{U}}(\nu))=\sum_{i=1}^{Md^2}\xi_{i,Md^2}^{\mathrm{U}}(\nu) Z_i^2$$ where $\xi_{Md^2}^{\mathrm{U}}(\nu)=(\xi_{1,Md^2}^{\mathrm{U}}(\nu),\dots,\xi_{Md^2,Md^2}^{\mathrm{U}}(\nu))^\top$ is the vector of the eigenvalues of the matrix
\begin{align*}
\bfnabla^{\mathrm{U}}_{\nu\nu}(M) = \bfJ^{-1/2}_M (\nu)\bfDelta^{\mathrm{U}}_{\nu\nu}
\bfJ^{-1/2}_M (\nu),
\end{align*}
where $\bfDelta^{\mathrm{U}}_{\nu\nu}$ is given in \eqref{chatgammaU},
$\bfJ_M (\nu)= \left(
    \begin{array}{ccccc}
    \bfSigma_{\bfepsilon}(\nu-1) & \bfzero        & \bfzero & \hdots & \bfzero \\
    \bfzero         & \bfSigma_{\bfepsilon}(\nu-2) & \bfzero & \hdots & \bfzero \\
    \vdots          &                &          & \ddots  &  \vdots   \\
    \bfzero         & \bfzero        & \bfzero & \hdots &\bfSigma_{\bfepsilon}(\nu-M)
    \end{array}
    \right) \otimes \bfSigma_{\bfepsilon}(\nu)$
is a $Md^2\times Md^2$ block diagonal matrix  and $Z_1,\dots,Z_{Md^2}$ are independent $\mathcal{N}(0,1)$ variables.

The asymptotic distribution of the global portmanteau test statistics (that takes into account all the seasons) $\bfQ_M$ and $\bfQ_M^{\ast}$ defined in \eqref{QM} and \eqref{QMast} are also a weighted sum of chi-square random variables:
$$Z_M(\xi^{\mathrm{U}}_{sMd^2})=\sum_{i=1}^{sMd^2}\xi^{\mathrm{U}}_{i,sMd^2} Z_i^2$$ where
$\xi^{\mathrm{U}}_{sMd^2}=(\xi^{\mathrm{U}}_{1,sMd^2},\dots,\xi^{\mathrm{U}}_{sMd^2,sMd^2})^\top$ denotes the vector of the eigenvalues of the block  matrix
$\bfnabla^{\mathrm{U}}(M)= [\bfnabla^{\mathrm{U}}_{\nu\nu'}(M)]_{1\leq \nu,\nu'\leq s}$, with blocs given by
\begin{align*}
\bfnabla^{\mathrm{U}}_{\nu\nu'}(M) =  \bfJ^{-1/2}_M (\nu)\bfDelta^{\mathrm{U}}
 \bfJ^{-1/2}_M (\nu)
\end{align*}
where
$\bfDelta^{\mathrm{U}}_{\nu\nu'}$ is given in \eqref{chatrhoglobalU}.
\end{theorem}
\begin{remark}
\label{remKhi2Un}
In view of Remark~\ref{remrhoUn} when $M$ is large, the matrix $$\bfnabla^{\mathrm{U}}_{\mathrm{S}}(\nu,\nu)\simeq\bfI_{d^2M} -
\bfJ^{-1/2}_M (\nu)\left\{ \left[\bfA^\top_{M}(\nu)\bfOmega^{-1}(\nu) \bfA_{M}(\nu)\right]\otimes \bfSigma_{\bfepsilon}(\nu)\right\} \bfJ^{-1/2}_M (\nu),
$$ (resp. $\bfnabla_{\mathrm{S}}^{\mathrm{U}}\simeq[\bfnabla^{\mathrm{U}}_{\mathrm{S}}(\nu,\nu)]_{\nu=1,\dots,s}$)  is close  to a projection matrix. Its eigenvalues are therefore equal to 0
and 1. The number of eigenvalues equal to 1 is $\tr(\bfnabla^{\mathrm{U}}_{\mathrm{S}}(\nu,\nu))=\tr(\bfI_{d^2 \{M-p(\nu)\}})=d^2\{M-p(\nu)\}$ (resp. $\tr(\bfnabla_{\mathrm{S}}^{\mathrm{U}}
)=\tr(\bfI_{d^2\sum_{\nu=1}^s \{M-p(\nu)\}})=d^2\sum_{\nu=1}^s \{M-p(\nu)\}$).
Therefore under $H_0(\nu)$ (resp. $H_0$), the asymptotic distributions of the statistics defined in \eqref{QMnu} and \eqref{QMnuast}  can be approximated by a $\chi_{d^2\{M-p(\nu)\}}^2$ (resp. a $\chi_{d^2\sum_{\nu=1}^s \{M-p(\nu)\}}^2$. When $d=1$, we also retrieve the results obtained by \cite{McL1994}.

As a consequence to these results, the chi-squared distributions are not applicable when $M \leq p(\nu)$ in contrast to our results (see Theorem \ref{theolimitdistBPUn}). In fact, the exact distributions of the portmanteau test statistics are better approximated by those of weighted sums of chi-squared random variables. Our results show that in the general PVAR case, the weights in the asymptotic distributions of the portmanteau  procedures may be relatively far from zero, and thus, adjusting the degrees of freedom does not represent a solution in the present framework.
\end{remark}

\begin{theorem}
\label{theolimitdistBPCon}
Under Assumptions of Theorem \ref{loi_res_gamCon} and $H_0(\nu)$, $\nu=1,\dots,s$ or $H_0$,  the
statistics $\bfQ_M(\nu)$ and $\bfQ_M^{\ast}(\nu)$ defined in \eqref{QMnu} and \eqref{QMnuast} converge in distribution, as $N\rightarrow\infty,$ to
$$Z_M^{\nu}(\xi_M^{\mathrm{R}}(\nu))=\sum_{i=1}^{Md^2}\xi_{i,M}^{\mathrm{R}}(\nu) Z_i^2$$ where
$\xi_{Md^2}^{\mathrm{R}}(\nu)=(\xi_{1,Md^2}^{\mathrm{R}}(\nu),\dots,\xi_{Md^2,Md^2}^{\mathrm{R}}(\nu))^\top$ is the vector of the eigenvalues of the matrix
\begin{align*}
\bfnabla^{\mathrm{R}}_{\nu\nu}(M) = \bfJ^{-1/2}_M (\nu)\bfDelta^{\mathrm{R}}_{\nu\nu}
\bfJ^{-1/2}_M (\nu),
\end{align*}
where $\bfDelta^{\mathrm{R}}_{\nu\nu}$ is given in \eqref{chatgammaR}  and $Z_1,\dots,Z_{Md^2}$ are independent $\mathcal{N}(0,1)$ variables.

The asymptotic distribution of the global portmanteau test statistics (that takes into account all the seasons) $\bfQ_M$ and $\bfQ_M^{\ast}$ defined in \eqref{QM} and \eqref{QMast} is also a weighted sum of chi-squared random variables:
$$Z_M(\xi^{\mathrm{R}}_{sMd^2})=\sum_{i=1}^{sMd^2}\xi^{\mathrm{R}}_{i,sMd^2} Z_i^2$$ where
$\xi^{\mathrm{R}}_{sMd^2}=(\xi^{\mathrm{R}}_{1,sMd^2},\dots,\xi^{\mathrm{R}}_{sMd^2,sMd^2})^\top$ denotes the vector of the eigenvalues of the  block  matrix
$\bfnabla^{\mathrm{R}}(M)= [\bfnabla^{\mathrm{R}}_{\nu\nu'}(M)]_{1\leq \nu,\nu'\leq s}$, with blocs given by
\begin{align*}
\bfnabla^{\mathrm{R}}_{\nu\nu'}(M) =  \bfJ^{-1/2}_M (\nu)\bfDelta^{\mathrm{R}}_{\nu\nu'}
 \bfJ^{-1/2}_M (\nu),
\end{align*}
where $\bfDelta^{\mathrm{R}}_{\nu\nu'}$ is given in \eqref{chatrhoglobalR}.
\end{theorem}
\begin{remark}
\label{remKhi2Con}

Under the linear constraints \eqref{constraints}, in view of Remark~\ref{remrhoCon} when $M$ is large, the matrix $$\bfnabla^{\mathrm{R}}_{\mathrm{S}}(\nu,\nu)\simeq\bfI_{d^2M} -
\bfJ^{-1/2}_M (\nu)\left\{ \left(\bfA^{\top}_{\nu}  \otimes \bfI_{d} \right)\bfR(\nu)\left[ \bfR^{\top}(\nu) \{ \bfOmega(\nu) \otimes \bfSigma^{-1}_{\bfepsilon}(\nu) \} \bfR(\nu) \right]^{-1}
  \bfR^{\top}(\nu)
  \left(\bfA_{\nu}  \otimes \bfI_{d} \right)\right\} \bfJ^{-1/2}_M (\nu),
$$ (resp. $\bfnabla_{\mathrm{S}}^{\mathrm{R}}\simeq[\bfnabla^{\mathrm{R}}_{\mathrm{S}}(\nu,\nu)]_{\nu=1,\dots,s}$)  is close  to a projection matrix. Its eigenvalues are therefore equal to $0$ and $1$. The number of eigenvalues equal to $1$ is $\tr(\bfnabla^{\mathrm{R}}_{\mathrm{S}}(\nu,\nu))=\tr(\bfI_{d^2M-K(\nu)})=d^2M-K(\nu)$ (resp. $\tr(\bfnabla_{\mathrm{S}}^{\mathrm{R}})=\tr(\bfI_{\sum_{\nu=1}^s \{d^2M-K(\nu)\}})=\sum_{\nu=1}^s \{d^2M-K(\nu)\}$). Therefore under $H_0(\nu)$ (resp. $H_0$), the asymptotic distributions of the statistics defined in \eqref{QMnu} and \eqref{QMnuast}  can be approximated by a $\chi_{d^2M-K(\nu)}^2$ (resp. a $\chi_{\sum_{\nu=1}^s \{d^2M-K(\nu)\}}^2$. We retrieve the results given in \cite{UD09}.

As a consequence to these results, the chi-squared distribution  are not applicable when $M \leq K(\nu)$ in contrast to our results (see Theorem \ref{theolimitdistBPCon}), which corresponds to the results obtained by \cite{DL13} for parsimonious strong PVAR models.
\end{remark}
\begin{corollary}\label{bobo}
Under the following alternative hypothesis  
\[ 
H_1(\nu):\text{ there exists  } h\in\{1,\ldots,M\} \text{ such that } \bfgamma_{\bfepsilon}(h;\nu) \neq \bfzero \quad\text{or }\quad
H_1(\nu): \vec\left\{\bfrho_{\bfepsilon}(h;\nu)\right\} \neq \bfzero,\quad\text{for }\quad\nu=1,\dots,s,
\]
one may prove that under $H_1(\nu)$
\begin{align*}
\frac{\bfQ_M(\nu)}{N}\quad\text{or}\quad \frac{\bfQ^\ast_M(\nu)}{N}
 & =   \left(\bfc_{\hat{\bfepsilon}}(\nu)  - \bfgamma_{\bfepsilon}(\nu) \right)^\top \bfK^{-1}(\nu){\bfN}^{-1}(\nu)\bfK^{-1}(\nu) \left(\bfc_{\hat{\bfepsilon}}(\nu)  - \bfgamma_{\bfepsilon}(\nu) \right)  \\&\qquad + 2 \left(\bfc_{\hat{\bfepsilon}}(\nu)  - \bfgamma_{\bfepsilon}(\nu) \right)^\top \bfK^{-1}(\nu){\bfN}^{-1}(\nu)\bfK^{-1}(\nu)  \bfgamma_{\bfepsilon}(\nu)  +  \bfgamma^\top_{\bfepsilon}(\nu) \bfK^{-1}(\nu){\bfN}^{-1}(\nu)\bfK^{-1}(\nu) \bfgamma_{\bfepsilon}(\nu)  + \mathrm{o}_{\mathbb P} \left(\frac{1}{N}\right)  \\ 
 & \xrightarrow[N\to\infty]{\mathbb{P}} \bfgamma^\top_{\bfepsilon}(\nu) \bfK^{-1}(\nu){\bfN}^{-1}(\nu)\bfK^{-1}(\nu) \bfgamma_{\bfepsilon}(\nu)=\bfrho^\top_{\bfepsilon}(\nu) {\bfN}^{-1}(\nu) \bfrho_{\bfepsilon}(\nu),
\end{align*}
where $\bfK(\nu)= \diag\left[\left\{ \bfD^\ast_{{\bfepsilon}}(\nu-1),\dots,\bfD^\ast_{{\bfepsilon}}(\nu-M)\right\}  \otimes
                                    \bfD^\ast_{{\bfepsilon}}(\nu) \right]$ and $\bfN(\nu)=\left(\diag\left[\bfR^\ast_{{\bfepsilon}}(0; \nu-1),\dots, \bfR^\ast_{{\bfepsilon}}(0; \nu-M)\right]\right) \otimes \bfR^\ast_{{\bfepsilon}}(0; \nu)$ with 
$\bfR^\ast_{\bfepsilon}(0; \nu)=\left(\bfD^\ast_{\bfepsilon}(\nu)\right)^{-1}\bfGamma_{\bfepsilon}(0; \nu)
                                    \left(\bfD^\ast_{\bfepsilon}(\nu)\right)^{-1}$.     
Therefore the test statistics $\bfQ_M(\nu)$ and $\bfQ_M^{\ast}(\nu)$  are consistent in detecting $H_1(\nu)$.
\end{corollary}
The proof of this corollary is also postponed to Section Appendix.

Theorems \ref{theolimitdistBPUn} and \ref{theolimitdistBPCon} show that, for the asymptotic distributions of \eqref{QMnu} and \eqref{QMnuast} (resp. of \eqref{QM} and \eqref{QMast}), the chi-squared approximations are no longer valid in the framework of weak PVAR models.

The true asymptotic distributions $\bfnabla^{\mathrm{U}}_{\nu\nu'}(M)$ and $\bfnabla^{\mathrm{R}}_{\nu\nu'}(M)$ (resp. $\bfnabla^{\mathrm{U}}(M)$ and $\bfnabla^{\mathrm{R}}(M)$) depend on the periodic nuisance parameters involving $\bfSigma_{{\bfepsilon}}(\nu-h)$ for $h=0,1,\dots,M$ and the elements of $\bfXi^{\mathrm{U}}_{\nu\nu'}$ or $\bfXi^{\mathrm{R}}_{\nu\nu'}$. Consequently, in order to obtain the asymptotic distribution of the portmanteau statistics \eqref{QMnu} and \eqref{QMnuast} (resp. \eqref{QM} and \eqref{QMast}) under weak assumptions on the multivariate periodic  noise, one needs a consistent estimator of the asymptotic covariance matrix $\bfnabla^{\mathrm{U}}_{\nu\nu'}(M)$ and $\bfnabla^{\mathrm{R}}_{\nu\nu'}(M)$ (resp. $\bfnabla^{\mathrm{U}}(M)$ and $\bfnabla^{\mathrm{R}}(M)$).
Matrices  $\bfUpsilon_M(\nu)$ and $\bfSigma_{{\bfepsilon}}(\nu-h)$ for $h=0,1,\dots,M$ can be estimated by its empirical counterpart, respectively by
$$\hat{\bfUpsilon}_M(\nu)=-\frac{1}{N}\sum_{n=0}^{N-1}
\left\{\left(\hat{\bfepsilon}^\top_{ns+\nu-1},\dots,\hat{\bfepsilon}^\top_{ns+\nu-M}\right)^\top\otimes\bfX^\top_{n}(\nu)\otimes\bfI_d\right\}\text{ and }\hat{\bfSigma}_{{\bfepsilon}}(\nu-h)=\frac{1}{N}\sum_{n=0}^{N-1}\hat{\bfepsilon}_{ns+\nu-h}\hat{\bfepsilon}_{ns+\nu-h}^\top.$$
The matrix $\bfOmega(\nu)$ can also be  estimated empirically by the squared matrix $\hat{\bfOmega}_N(\nu)$ of order $dp(\nu)$ defined by:
$\hat{\bfOmega}_N(\nu)={N}^{-1} \bfX_n(\nu) \bfX_n^{\top}(\nu)$.

The estimation of the long-run variance (LRV) matrix $\bfXi^{\mathrm{U}}_{\nu\nu'}$ or $\bfXi^{\mathrm{R}}_{\nu\nu'}$ is more complicated. In the literature, two types of estimators are generally employed: heteroskedasticity and autocorrelation consistent (HAC) estimators based on kernel methods (see \cite{newey} and \cite{A91} for general references, and \cite{FZ07} for an application to testing strong linearity in weak ARMA models) and the spectral density (SP) estimators (see e.g. \citep{B74} and \cite{haan} for a general reference; see also \cite{BMF11} for an application to a weak VARMA model). Thus following the arguments developed in \cite{BMU23}, the matrix $\bfXi^{\mathrm{U}}_{\nu\nu'}$ (resp. $\bfXi^{\mathrm{R}}_{\nu\nu'}$)  can be weakly estimated by HAC or SP. Let $\hat{\bfXi}^{\mathrm{U}}_{\nu\nu'}$ (resp. $\hat{\bfXi}^{\mathrm{R}}_{\nu\nu'}$) be an estimator obtained by only replacing the process $\left(\bfW_{ns+\nu}\right)_{n\in\mathbb{Z}}$  in \citet[Theorems 4.1 or Theorems 4.2]{BMU23}, by $\left(\bfW_{ns+\nu}^{\mathrm{U}}\right)_{n\in\mathbb{Z}}$ (resp. $\left(\bfW_{ns+\nu}^{\mathrm{R}}\right)_{n\in\mathbb{Z}}$).

We let $\hat{\bfnabla}^{\mathrm{U}}_{\nu\nu'}(M)$ and $\hat{\bfnabla}^{\mathrm{R}}_{\nu\nu'}(M)$ (resp. $\hat{\bfnabla}^{\mathrm{U}}(M)$ and $\hat{\bfnabla}^{\mathrm{R}}(M)$) the matrices obtained by replacing $\bfXi^{\mathrm{U}}_{\nu\nu'}$ and $\bfXi^{\mathrm{R}}_{\nu\nu'}$ by $\hat{\bfXi}^{\mathrm{U}}_{\nu\nu'}$ and $\hat{\bfXi}^{\mathrm{R}}_{\nu\nu'}$, $\bfUpsilon_M(\nu)$ by $\hat{\bfUpsilon}_M(\nu)$ and $\bfSigma_{{\bfepsilon}}(\nu-h)$ by $\hat{\bfSigma}_{{\bfepsilon}}(\nu-h)$ in $\bfnabla^{\mathrm{U}}_{\nu\nu'}(M)$ and $\bfnabla^{\mathrm{R}}_{\nu\nu'}(M)$ (resp. $\bfnabla^{\mathrm{U}}(M)$ and $\bfnabla^{\mathrm{R}}(M)$). Denote by
\begin{eqnarray*}
\hat{\xi}_{Md^2}^{\mathrm{U}}(\nu)=(\hat{\xi}_{1,Md^2}^{\mathrm{U}}(\nu),\dots,\hat{\xi}_{Md^2,Md^2}^{\mathrm{U}}(\nu))^\top\;
&\text{ and }&\; \hat{\xi}_{Md^2}^{\mathrm{R}}(\nu)=(\hat{\xi}_{1,Md^2}^{\mathrm{R}}(\nu),\dots,\hat{\xi}_{Md^2,Md^2}^{\mathrm{R}}(\nu))^\top\\
(\text{resp. } \hat{\xi}^{\mathrm{U}}_{sMd^2}=(\hat{\xi}^{\mathrm{U}}_{1,sMd^2},\dots,\hat{\xi}^{\mathrm{U}}_{sMd^2,sMd^2})^\top \; &\text{ and }&\; \hat{\xi}^{\mathrm{R}}_{sMd^2}=(\hat{\xi}^{\mathrm{R}}_{1,sMd^2},\dots,\hat{\xi}^{\mathrm{R}}_{sMd^2,sMd^2})^\top)
\end{eqnarray*}
the vector of the eigenvalues of
$\hat{\bfnabla}^{\mathrm{U}}_{\nu\nu'}(M)$ and $\hat{\bfnabla}^{\mathrm{R}}_{\nu\nu'}(M)$ (resp. $\hat{\bfnabla}^{\mathrm{U}}(M)$ and $\hat{\bfnabla}^{\mathrm{R}}(M)$). At the asymptotic level $\alpha$ and for a specified season $\nu$, the modified test consists in rejecting the adequacy of the weak PVAR model when
$$\lim_{N\to\infty}\mathbb{P}\left(\bfQ_M(\nu)>S_M^\nu(1-\alpha)\right)=\lim_{N\to\infty}\mathbb{P}\left(\bfQ_M^{\ast}(\nu)>S_M^\nu(1-\alpha)\right)=\alpha,$$
where $S_M^\nu(1-\alpha)$ is such that $\mathbb{P}\left\{Z_{Md^2}^\nu(\hat\xi_{Md^2}^{\mathrm{U}\text{ or }\mathrm{R}})(\nu)>S_M^\nu(1-\alpha)\right\}=\alpha$.
We emphasize the fact that the proposed modified versions of the  statistics are more difficult to implement because their critical values have to be computed from the data while the  critical values of the standard method are simply deduced from a $\chi^2$-table.  We shall evaluate the $p$-values $$\mathbb{P}\left\{Z_{Md^2}^\nu(\hat\xi_{Md^2}^{\mathrm{U}\text{ or }\mathrm{R}})(\nu)>{\bfQ}_M(\nu)\right\}\;\mbox{ and }\;
\mathbb{P}\left\{Z_{Md^2}^\nu(\hat\xi_{Md^2}^{\mathrm{U}\text{ or }\mathrm{R}})(\nu)>{\bfQ}_M^{\ast}(\nu)\right\},$$
 with $Z_{Md^2}^\nu(\hat\xi_{Md^2}^{\mathrm{U}\text{ or }\mathrm{R}})(\nu)=
\sum_{i=1}^{Md^2}\hat{\xi}_{i,Md^2}^{\mathrm{U}\text{ or }\mathrm{R}}(\nu) Z_i^2,$ by  means of the Imhof algorithm (see \cite{I61}) or other exact methods.

The test procedures for ${\bfQ}_M$ and ${\bfQ}_M^{\ast}$ defined in \eqref{QM} and \eqref{QMast} are similar but they are based on the $sd^2M$ empirical eigenvalues of  $\hat{\bfnabla}^{\mathrm{U}}(M)$ and $\hat{\bfnabla}^{\mathrm{R}}(M)$.

\begin{remark}\label{hatrho_band} Let $\bfr^i_{\hat{\bfepsilon}}(\nu)$ be the $i$-th component of the $(Md^2)\times 1$ vector $\bfr_{\hat{\bfepsilon}}(\nu)$ of the residual autocorrelation.
For any $1\leq h\leq Md^2$ and for a specified season $\nu$, a $100(1-\alpha)\%$ confidence region for $\bfr^h_{{\bfepsilon}}(\nu)$ is given by
$$-u_\alpha\frac{\sqrt{\hat{\bfnabla}^{\mathrm{U}}_{\nu\nu}(M)(h,h)}}{\sqrt{N}}\leq\bfr^h_{\hat{\bfepsilon}}(\nu)\leq u_\alpha\frac{\sqrt{\hat{\bfnabla}^{\mathrm{U}}_{\nu\nu}(M)(h,h)}}{\sqrt{N}}\quad
\left(\text{resp. } -u_\alpha\frac{\sqrt{\hat{\bfnabla}^{\mathrm{R}}_{\nu\nu}(M)(h,h)}}{\sqrt{N}}\leq\bfr^h_{\hat{\bfepsilon}}(\nu)\leq u_\alpha\frac{\sqrt{\hat{\bfnabla}^{\mathrm{R}}_{\nu\nu}(M)(h,h)}}{\sqrt{N}}\right)$$
where $u_\alpha$ denotes the quantile of order $1-\alpha$ of the $\mathcal{N}(0,1)$ distribution.
\end{remark}

In the next section, the test statistics proposed in this section are
illustrated in a small empirical study.\\

\section{Simulations}\label{simul}
\noindent By means of a small Monte Carlo experiment, we investigate the finite sample properties of the modified version of the portmanteau test. The following data generating process (DGP) are used:
\begin{eqnarray*}
\mbox{DGP$_{1}$} &:& \bfY_{ns+\nu} = \bfPhi(\nu)\bfY_{ns+\nu-1} + \bfepsilon_{ns+\nu},\\
\mbox{DGP$_{2}$} &:& \bfY_{ns+\nu} = \bfPhi_\mathcal{C}(\nu) \bfY_{ns+\nu-1} + \bfepsilon_{ns+\nu},\\
\mbox{DGP$_{3}$} &:& \bfY_{ns+\nu} = \bfPhi_1(\nu)\bfY_{ns+\nu-1} + \bfPhi_2(\nu)\bfY_{ns+\nu-2} + \bfepsilon_{ns+\nu}.
\end{eqnarray*}
We considered the case of four seasons, that is $s=4$.
For the purposes of our illustration, no parameter constraints were hypothesized for DGP$_1$ (Table~\ref{DGP1}), but for DGP$_2$ (Table~\ref{DGP2}) it was assumed that the zero-valued parameters in $\bfPhi_{\mathcal{C}}(1)$, $\bfPhi_{\mathcal{C}}(2)$, $\bfPhi_{\mathcal{C}}(3)$ and $\bfPhi_{\mathcal{C}}(4)$ were known.

\begin{table}[H]
\caption{Parameters of DGP$_1$ models used in the simulation}\label{DGP1}
\centering
\begin{tabular}{lcccccccccccc}
  \toprule
  MODEL & & \multicolumn{2}{c}{$\bfPhi(1)$} & & \multicolumn{2}{c}{$\bfPhi(2)$} & & \multicolumn{2}{c}{$\bfPhi(3)$} & &\multicolumn{2}{c}{$\bfPhi(4)$} \\
  \midrule
  DGP$_1$ & & 0.50 & 0.30 & & 0.42 & 0.24 & & -0.80 & 0.20 & & -0.30 & 0.50\\
          & & 0.10 & 0.20 & &-0.20 & 0.50 & & 0.60 & 0.70 & & 0.90 & -0.20\\
  \bottomrule
\end{tabular}
\end{table}

\begin{table}[H]
\caption{Parameters of DGP$_2$ models used in the simulation}\label{DGP2}
\centering
\begin{tabular}{lcccccccccccc}
  \toprule
  MODEL & & \multicolumn{2}{c}{$\bfPhi_{\mathcal{C}}(1)$} & & \multicolumn{2}{c}{$\bfPhi_{\mathcal{C}}(2)$} & & \multicolumn{2}{c}{$\bfPhi_{\mathcal{C}}(3)$} & &\multicolumn{2}{c}{$\bfPhi_{\mathcal{C}}(4)$} \\
  \midrule
  DGP$_2$ & & 0.95 & 0.00 & & -0.90 & 0.00 & & -0.85 & 0.00 & & -0.95 & 0.00\\
          & & 0.00 & 0.90 & & 0.00 & 0.95 & & 0.00 & 0.90 & & 0.00 & -0.95\\
  \bottomrule
\end{tabular}
\end{table}

\begin{table}[H]
\caption{Parameters of DGP$_3$ models used in the simulation}\label{DGP3}
\centering
\begin{tabular}{lcccccccccccc}
  \toprule
  MODEL & & \multicolumn{2}{c}{$\bfPhi_1(1)$} & & \multicolumn{2}{c}{$\bfPhi_1(2)$} & & \multicolumn{2}{c}{$\bfPhi_1(3)$} & &\multicolumn{2}{c}{$\bfPhi_1(4)$} \\
  \midrule
  DGP$_3$ & & 0.60 & 0.30 & & -0.30 & 0.40 & & 0.30 & 0.30 & & -0.40 & -0.40\\
          & & 0.60 & 0.20 & & 0.20 & 0.40 & & 0.30 & 0.20 & & 0.30 & -0.40\\
  \midrule
        & & \multicolumn{2}{c}{$\bfPhi_2(1)$} & & \multicolumn{2}{c}{$\bfPhi_2(2)$} & & \multicolumn{2}{c}{$\bfPhi_2(3)$} & &\multicolumn{2}{c}{$\bfPhi_2(4)$} \\
  \midrule
         & & 0.40 & 0.60 & & -0.30 & 0.30 & & 0.20 & 0.50 & &  0.30 & 0.50\\
         & & -0.20 & 0.50 & & 0.30 & -0.40 & & -0.30 & -0.30 & & 0.50 & 0.30\\
  \bottomrule
\end{tabular}
\end{table}

\begin{table}[H]
\caption{Error covariance matrices used in the simulation}\label{error}
\centering
\begin{tabular}{cccccccccccc}
  \toprule
    & \multicolumn{2}{c}{$\bfSigma_{\bfepsilon}(1)$} & & \multicolumn{2}{c}{$\bfSigma_{\bfepsilon}(2)$} & & \multicolumn{2}{c}{$\bfSigma_{\bfepsilon}(3)$} & &\multicolumn{2}{c}{$\bfSigma_{\bfepsilon}(4)$} \\
  \midrule
           & 1.00 & 0.50 & & 1.00 & 0.30 & & 1.00 & 0.20 & & 1.00 & 0.10\\
           & 0.50 & 1.00 & & 0.30 & 1.00 & & 0.20 & 1.00 & & 0.10 & 1.00\\
  \bottomrule
\end{tabular}
\end{table}

We studied the empirical frequencies of rejection of the null hypothesis of adequacy at the $5\%$ and $10\%$ nominal levels, for each of three series length ($N = 200$, $1000$ and $5000$). For each experiment, $1000$ independent realizations were generated. For these nominal levels, the empirical  relative frequency of rejection size over the $1000$ independent replications should vary  respectively within the confidence intervals $[3.6\%, 6.4\%]$  and $[8.1\%,11.9\%]$ with probability 95\% and  $[3.3\%, 6.9\%]$ and $[7.6\%, 12.5\%]$ with probability 99\% under the assumption that the true probabilities of rejection are respectively $\alpha=5\%$ and $\alpha=10\%$. When the relative rejection frequencies  are outside the significant limits with probability 95\%, they are displayed in bold type. For each realization of the DGP defined by DGP$_1$, a PVAR model of order one was estimated by least squares estimators, as described in~Section~3.1 of~\cite{BMU23}.
When the DGP was given by DGP$_2$, the zero-valued parameters in
$\bfPhi_{\mathcal{C}}(\nu)$ were taken into account by properly defining the constraint matrix $\bfR(\nu)$, $\nu=1,2,3,4$, and the parameters were estimated using the procedure described in~Section~3.2 of~\cite{BMU23}.

For each residual time series, the portmanteau test statistics $Q_M^{\ast}(\nu)$~(\ref{QMnuast}), $\nu = 1,\ldots, s$ and $Q_M^{\ast}$~(\ref{QMast}) were calculated using the Ljung-Box-McLeod correction factor, for $M\in\{1,2,3,6,8,10\}$. We indicate the conventions that we adopt in the discussion and in the tables:
\begin{itemize}
\item $Q^1_M(\nu)$ and $GQ^1_M(\nu)$ refer to the modified $Q^*_M(\nu)$ and global modified tests as defined in~\citet{UD09}.
\item $Q^2_M(\nu)$ and $GQ^2_M(\nu)$ refer to the modified $Q^*_M(\nu)$ and global modified tests as defined in~\citet{DL13}.
\item $Q^3_M(\nu)$ and $GQ^3_M(\nu)$ refer to the modified test using the statistics~(\ref{QMnuast}) and~(\ref{QMast}).
\end{itemize}
Compared to the modified portmanteau test proposed by~\citet{DL13} and~\citet{UD09}, we use a VAR spectral estimator approach to estimate the asymptotic covariance matrix of a vector autocorrelations residuals. The implementation of this method requires a choice of the VAR order $r$. In the strong PVAR cases, we fixed $r = 1$ as it can be shown that, for strong models, the VAR spectral estimator is consistent with any fixed value of $r$. For weak PVAR models, the VAR order $r$ is set as $r=1,2,3$ and is automatically selected by Akaike Information Criterion using the function \textit{VARselect()} of the \textit{vars} R package. Imhof's algorithm~\citep{I61} has been used to obtain the critical values; we used the R package \textit{CompQuadForm} available from CRAN. The results concerning the standard tests are not presented here, because they are similar to those of the modified tests. 

\subsection{Empirical size}
\subsubsection{Unconstrained strong PVAR model}
\noindent First we consider the unconstrained strong PVAR defined by DGP$_1$ and where the stochastic process $\bfepsilon = \{ \bfepsilon_t, t \in \mathbb{Z} \}$ corresponds to a zero mean periodic white noise with the error covariance matrix $\bfSigma_{\nu}$ given in Table~\ref{error}.

The number of rejections of the null hypothesis of adequacy is reported in Tables~\ref{test_portmanteau_strong_unc} and~\ref{global_test_portmanteau_strong_unc}.
The results presented in Table~\ref{test_portmanteau_strong_unc} indicate that, for the test proposed by~\citet{UD09}, some rejection has been observed for small values of $M$ at both significance levels. This is not surprising because the $\chi^2_{d^2(M-1)}$ approximation is better for larger $M$. Note that for $M=1$, the empirical size is not available (n.a.) for the test proposed by~\citet{UD09} because the test is not applicable to $M\leq p$. On the other hand, the error of first kind is well controlled by the test proposed by~\citet{DL13} and our test even when $M$ is small. The results presented in Table~\ref{test_portmanteau_strong_unc} indicate that some under-rejection has been observed for large values of $M$ at both significance levels when $N = 200$ for our test. No wonder that $N$ has to be large enough as our test requires the estimation of the covariance structure of a high dimension multivariate process.

The empirical levels of the global portmanteau tests statistics are given in Table~\ref{global_test_portmanteau_strong_unc}. For our test some under-rejections occurred for large values of $M$ when $N=200$. When $N = 1000$, the results were generally reasonable at both significance levels. It is not surprising that large sample sizes are needed for the global test statistics given the high dimension of the process.

\begin{table}[H]
\caption{\scriptsize{Empirical levels (number of rejections of the null hypothesis over $1000$ replications) at the $5\%$ and $10\%$ significance levels for the modified portmanteau test statistics $Q_M(\nu)$ for $\nu=1,2,3,4$ and $M = 1,2,3,6,8,10$. The simulated model is a strong PVAR given by DGP$_1$. The numbers of years are equal to $N = 200$, $1000$ and $5000$}. $\mathrm{{Q}}_M^{1}(\nu)$ is the portmanteau test statistics defined by~(35) in~\citet{UD09}. $\mathrm{{Q}}_M^{2}(\nu)$ is the portmanteau test statistics defined by~(26) in~\citet{DL13}. $\mathrm{{Q}}_M^{3}(\nu)$ is the portmanteau test statistics defined by~(\ref{QMnuast}).
}\label{test_portmanteau_strong_unc}
\begin{center}
\centering
\tiny
\setlength{\tabcolsep}{0.8mm}
\begin{tabular}{|c|cccc|cccc|cccc|cccc|cccc|cccc|cccc|cccc|cccc|}
\hline
& \multicolumn{36}{c|}{$\alpha=0.05$} \\
& \multicolumn{12}{c}{$N=200$} & \multicolumn{12}{c}{$N=1000$} & \multicolumn{12}{c|}{$N=5000$}\\
\cline{2-37}
&\multicolumn{4}{c}{$Q^1_M(\nu)$} &  \multicolumn{4}{c|}{$Q^2_M(\nu)$} & \multicolumn{4}{c|}{$Q^3_M(\nu)$} & \multicolumn{4}{c}{$Q^1_M(\nu)$} &  \multicolumn{4}{c|}{$Q^2_M(\nu)$} & \multicolumn{4}{c|}{$Q^3_M(\nu)$} & \multicolumn{4}{c}{$Q^1_M(\nu)$} &  \multicolumn{4}{c|}{$Q^2_M(\nu)$} & \multicolumn{4}{c|}{$Q^3_M(\nu)$}\\
\hline
\diagbox{$M$}{$\nu$}& 1 & 2 & 3 & 4 & 1 & 2 & 3 & 4 & 1 & 2 & 3 & 4 & 1 & 2 & 3 & 4 & 1 & 2 & 3 & 4 & 1 & 2 & 3 & 4 & 1 & 2 & 3 & 4 & 1 & 2 & 3 & 4 & 1 & 2 & 3 & 4\\
\hline
1 & n.a. & n.a. & n.a. & n.a. & \textbf{6.9} & 5.5 & 5.1 & 5.8 & 4.4 & \textbf{3.5} & 4.7 & 3.7
& n.a. & n.a. & n.a. & n.a. & 4.9 & 5.6 & 6.4 & 5.5 & 5.8 & 5.7 & 5.0 & 5.6 &
n.a. & n.a. & n.a. & n.a. & 4.6 & 4.5 & 5.7 & 4.3 & 6.1 & 5.2 & 4.1 & 4.5\\
2 & \textbf{11.9} & 6.1 & 5.6 & 5.5 & 6.2 & 4.9 & 5.4 & 5.0 & \textbf{3.3} & 4.9 & 4.4 & 4.1 &
\textbf{13.5} & 4.5 & 4.7 & 6.1 & 6.1 & 4.6 & \textbf{6.5} & 5.0 & 5.3 & 5.4 & 5.0 & 4.0 &
\textbf{13.5} & 4.8 & \textbf{6.5} & 6.0 & 4.4 & 3.8 & 5.1 & 4.4 & 6.0 & 4.6 & 6.4 & 4.8\\
3 & 6.3 & 5.3 & 6.0 & 6.4 & 5.1 & 3.9 & 5.4 & 5.1 & \textbf{3.3} & 3.6 & 4.2 & 3.7 &
6.1 & 5.2 & 5.3 & 6.0 & \textbf{6.5} & 5.3 & 4.9 & 5.4 & 5.1 & 5.4 & 5.8 & 4.7 &
5.0 & 4.8 & 4.4 & 5.3 & 4.9 & 3.6 & 5.7 & 5.4 & 4.4 & 4.7 & 4.4 & 5.3\\
6 & 5.2 & 4.8 & 4.7 & 5.2 & 6.1 & 4.3 & 5.6 & 4.8 & \textbf{1.6} & \textbf{1.9} & \textbf{1.4} & \textbf{2.0} &
5.1 & 5.4 & 5.8 & 5.4 & 5.6 & 4.7 & 6.3 & 5.6 & 4.0 & 4.6 & 4.3 & 4.0 &
5.9 & 4.7 & 5.1 & 3.9 & 5.1 & 5.5 & 6.0 & 4.2 & 5.3 & 4.4 & 5.0 & 3.6\\
8 & 3.6 & 5.1 & 4.1 & 4.2 & 5.1 & 5.4 & 6.4 & 5.9 & \textbf{1.6} & \textbf{2.0} & \textbf{0.9} & \textbf{1.4} &
5.5 & 5.3 & 4.5 & 5.0 & 5.6 & 5.5 & 5.8 & 5.6 & 4.6 & 4.5 & 4.5 & 4.3 &
4.4 & 5.4 & 4.2 & 6.0 & 4.3 & 5.4 & 5.4 & 4.0 & 4.0 & 5.2 & 3.9 & 6.0\\
10 & 5.3 & 5.5 & 5.1 & 5.8 & 5.4 & 4.1 & 5.9 & 5.2 & \textbf{1.3} & \textbf{2.0} & \textbf{2.3} & \textbf{1.6} &
\textbf{6.5} & 4.7 & 5.7 & 4.2 & 5.4 & 6.4 & 6.4 & 6.4 & \textbf{3.4} & \textbf{3.5} & 4.7 & 3.8 &
5.4 & 5.0 & 5.0 & 5.6 & 5.5 & 4.6 & 5.5 & 4.9 & 5.3 & 4.6 & 4.8 & 5.1\\ \hline
& \multicolumn{36}{c|}{$\alpha=0.10$} \\
& \multicolumn{12}{c}{$N=200$} & \multicolumn{12}{c}{$N=1000$} & \multicolumn{12}{c|}{$N=5000$}\\
\cline{2-37}
&\multicolumn{4}{c}{$Q^1_M(\nu)$} &  \multicolumn{4}{c|}{$Q^2_M(\nu)$} &  \multicolumn{4}{c|}{$Q^3_M(\nu)$} & \multicolumn{4}{c}{$Q^1_M(\nu)$} &  \multicolumn{4}{c|}{$Q^2_M(\nu)$} & \multicolumn{4}{c|}{$Q^3_M(\nu)$} & \multicolumn{4}{c}{$Q^1_M(\nu)$} &  \multicolumn{4}{c|}{$Q^2_M(\nu)$} & \multicolumn{4}{c|}{$Q^3_M(\nu)$}\\
\hline
\diagbox{$M$}{$\nu$} & 1 & 2 & 3 & 4 & 1 & 2 & 3 & 4 & 1 & 2 & 3 & 4 & 1 & 2 & 3 & 4 & 1 & 2 & 3 & 4 & 1 & 2 & 3 & 4 & 1 & 2 & 3 & 4 & 1 & 2 & 3 & 4 & 1 & 2 & 3 & 4\\
\hline
1 & n.a. & n.a. & n.a. & n.a. & \textbf{12.3} & 11.0 & 10.9 & 10.7 & 9.8 & 10.5 & 8.5 & 10.1 & n.a. & n.a. & n.a. & n.a. & 10.1 & 10.5 & 11.6 & 10.8 & 10.7 & 10.1 & 10.1 & 11.6 &
n.a. & n.a. & n.a. & n.a. & 9.8 & 9.6 & 11.1 & 9.2 & 11.3 & 8.9 & 9.4 & 8.5\\
2 & \textbf{20.8} & 11.8 & \textbf{12.5} & 11.9 & 11.5 & 9.5 & 10.3 & 8.8 & \textbf{7.3} & 9.6 & 10.7 & \textbf{7.7} & \textbf{22.9} & 9.3 & 10.5 & 11.9 & 11.0 & 9.8 & 10.7 & 9.8 & 11.3 & 11.1 & 10.6 & 8.9 &
\textbf{23.9} & 10.4 & \textbf{12.5} & \textbf{12.4} & 8.2 & 8.4 & 10.3 & 9.8 & 9.8 & 9.9 & \textbf{12.0} & 9.1\\
3 & \textbf{12.5} & 10.0 & 10.9 & 11.1 & 10.2 & 10.4 & 9.9 & 9.4 & 8.9 & 8.5 & 8.9 & 9.3 & \textbf{12.4} & 10.7 & 9.3 & 11.2 & 10.8 & 11.1 & 9.4 & 9.9 & 11.1 & 9.7 & 10.4 & 9.7 &
11.3 & 8.8 & 9.8 & 9.8 & 10.2 & \textbf{7.8} & 10.3 & 10.2 & 9.6 & 8.2 & 9.6 & 9.8\\
6 & 11.3 & 9.6 & 9.6 & 11.8 & 9.6 & 10.0 & 10.1 & 11.3 & \textbf{5.6} & \textbf{5.5} & \textbf{5.8} & \textbf{5.8} &
9.9 & 9.4 & 11.4 & 11.0 & 11.1 & 10.3 & 11.5 & 9.9 & 8.7 & 9.5 & 9.3 & \textbf{7.8} & 11.3 & 9.9 & 10.9 & 9.1 & 10.3 & 9.8 & 10.7 & 9.5 & 10.7 & 9.5 & 10.1 & 9.0\\
8 & 9.7 & 10.9 & 8.8 & 10.2 & 10.0 & 10.1 & 11.3 & 11.0 & \textbf{4.7} & \textbf{5.1} & \textbf{3.9} & \textbf{3.9} & 10.1 & 10.7 & 10.2 & 8.8 & 11.3 & 11.3 & 10.4 & 10.1 & 9.2 & 10.2 & 8.7 & 9.8 &
10.1 & 10.0 & 9.1 & 10.3 & 9.3 & 10.2 & 10.9 & 9.9 & 10.1 & 10.0 & 8.7 & 10.1\\
10 & 10.1 & 10.8 & 10.5 & 10.4 & 10.5 & 10.0 & \textbf{12.3} & 11.5 & \textbf{4.6} & \textbf{7.2} & \textbf{6.4} & \textbf{5.1} & \textbf{12.4} & 10.9 & 11.3 & 8.9 & 11.0 & 11.4 & 11.2 & 10.7 & \textbf{7.4} & \textbf{7.9} & 10.1 & 8.2 &
11.9 & 8.8 & 10.3 & 11.5 & 10.6 & 9.1 & 11.2 & 9.0 & 11.6 & 8.7 & 10.0 & 11.4\\ \hline
\end{tabular}
\end{center}
\end{table}

\begin{table}[H]
\caption{\scriptsize{Empirical levels (number of rejections of the null hypothesis over $1000$ replications) at the $5\%$ and $10\%$ significance levels for the modified portmanteau test statistics $Q_M(\nu)$ for $\nu=1,2,3,4$ and $M = 1,2,3,6,8,10$. The simulated model is a strong PVAR given by DGP$_3$. The numbers of years are equal to $N = 200$, $1000$ and $5000$}. $\mathrm{{Q}}_M^{1}(\nu)$ is the portmanteau test statistics defined by~(35) in~\citet{UD09}. $\mathrm{{Q}}_M^{2}(\nu)$ is the portmanteau test statistics defined by~(26) in~\citet{DL13}. $\mathrm{{Q}}_M^{3}(\nu)$ is the portmanteau test statistics defined by~(\ref{QMnuast}).
}\label{test_portmanteau_strong_unc_2}
\begin{center}
\centering
\tiny
\setlength{\tabcolsep}{0.8mm}
\begin{tabular}{|c|cccc|cccc|cccc|cccc|cccc|cccc|cccc|cccc|cccc|}
\hline
& \multicolumn{36}{c|}{$\alpha=0.05$} \\
& \multicolumn{12}{c}{$N=200$} & \multicolumn{12}{c}{$N=1000$} & \multicolumn{12}{c|}{$N=5000$}\\
\cline{2-37}
&\multicolumn{4}{c}{$Q^1_M(\nu)$} &  \multicolumn{4}{c|}{$Q^2_M(\nu)$} & \multicolumn{4}{c|}{$Q^3_M(\nu)$} & \multicolumn{4}{c}{$Q^1_M(\nu)$} &  \multicolumn{4}{c|}{$Q^2_M(\nu)$} & \multicolumn{4}{c|}{$Q^3_M(\nu)$} & \multicolumn{4}{c}{$Q^1_M(\nu)$} &  \multicolumn{4}{c|}{$Q^2_M(\nu)$} & \multicolumn{4}{c|}{$Q^3_M(\nu)$}\\
\hline
\diagbox{$M$}{$\nu$}& 1 & 2 & 3 & 4 & 1 & 2 & 3 & 4 & 1 & 2 & 3 & 4 & 1 & 2 & 3 & 4 & 1 & 2 & 3 & 4 & 1 & 2 & 3 & 4 & 1 & 2 & 3 & 4 & 1 & 2 & 3 & 4 & 1 & 2 & 3 & 4\\
\hline
1 & n.a. & n.a. & n.a. & n.a. & 3.9 & 6.0 & 5.6 & 6.3 & 4.9 & 5.0 & 5.5 & 4.2
& n.a. & n.a. & n.a. & n.a. & 4.4 & 5.2 & 3.6 & 5.2 & 5.3 & 4.7 & 4.4 & 4.6 &
n.a. & n.a. & n.a. & n.a. & 5.8 & 5.2 & 5.8 & 5.1 & 5.6 & 5.1 & 4.8 & 6.1\\
2 & n.a. & n.a. & n.a. & n.a & 4.8 & \textbf{6.7} & 5.0 & 5.6 & 4.9 & 5.0 & 4.5 & 4.1 &
n.a. & n.a. & n.a. & n.a. & 4.9 & 5.7 & 4.5 & 5.5 & 4.5 & 5.0 & 3.6 & 4.9 &
n.a. & n.a. & n.a. & n.a. & 4.6 & 4.6 & 4.6 & 5.0 & 4.8 & 5.8 & 5.9 & 5.8\\
3 & \textbf{15.4} & \textbf{17.1} & \textbf{11.2} & \textbf{10.4} & 4.3 & \textbf{6.6} & 4.9 & 6.3 & 4.1 & \textbf{3.2} & 5.0 & 3.7 &
\textbf{14.8} & \textbf{16.7} & \textbf{8.8} & \textbf{9.2} & 5.5 & 4.9 & 3.8 & 4.4  & 4.3 & 5.4 & 5.2 & 5.3 &
\textbf{14.4} & \textbf{13.5} & \textbf{11.1} & \textbf{9.6} & 5.7 & 4.4 & 5.4 & 4.8 & 6.3 & 4.5 & 5.4 & 5.7\\
6 & 6.1 & 5.7 & 4.6 & 5.4 & 5.8 & 5.7 & 4.4 & 5.2 & \textbf{1.8} & \textbf{2.4} & \textbf{3.5} & \textbf{3.3} &
5.8 & 5.2 & 5.6 & 5.5 & 5.0 & 4.8 & 5.4 & 5.1 & 4.6 & 4.7 & 5.6 & \textbf{3.4} &
5.6 & 4.8 & 4.9 & 5.1 & 5.0 & 4.6 & 4.8 & 5.1 & 4.0 & 4.4 & 5.2 & 5.0\\
8 & 4.7 & 5.8 & 5.2 & 5.5 & 4.8 & 6.0 & 5.2 & 6.1 & 4.5 & 5.8 & 5.2 & 6.0 &
5.0 & 3.7 & 5.1 & 5.2 & 5.1 & 3.7 & 4.8 & 5.2 & 4.7 & 3.6 & 4.8 & 4.1 &
4.9 & 4.2 & 4.5 & 4.5 & 4.9 & 4.2 & 4.1 & 4.5 & 4.9 & 6.0 & 6.1 & 5.0\\
10 & 5.0 & 5.3 & 4.7 & 5.4 & 5.5 & 5.7 & 5.2 & 5.7 & \textbf{10.9} & \textbf{12.1} & \textbf{11.8} & \textbf{10.6} &
5.3 & \textbf{3.5} & 4.6 & 5.7 & 5.3 & \textbf{3.5} & 4.6 & 5.7 & 4.2 & 4.3 & 5.1 & 4.5 &
4.1 & 3.7 & 4.6 & 5.2 & 4.1 & 3.7 & 4.6 & 5.1 & 5.5 & 4.7 & 4.9 & 3.7\\ \hline
& \multicolumn{36}{c|}{$\alpha=0.10$} \\
& \multicolumn{12}{c}{$N=200$} & \multicolumn{12}{c}{$N=1000$} & \multicolumn{12}{c|}{$N=5000$}\\
\cline{2-37}
&\multicolumn{4}{c}{$Q^1_M(\nu)$} &  \multicolumn{4}{c|}{$Q^2_M(\nu)$} &  \multicolumn{4}{c|}{$Q^3_M(\nu)$} & \multicolumn{4}{c}{$Q^1_M(\nu)$} &  \multicolumn{4}{c|}{$Q^2_M(\nu)$} & \multicolumn{4}{c|}{$Q^3_M(\nu)$} & \multicolumn{4}{c}{$Q^1_M(\nu)$} &  \multicolumn{4}{c|}{$Q^2_M(\nu)$} & \multicolumn{4}{c|}{$Q^3_M(\nu)$}\\
\hline
\diagbox{$M$}{$\nu$} & 1 & 2 & 3 & 4 & 1 & 2 & 3 & 4 & 1 & 2 & 3 & 4 & 1 & 2 & 3 & 4 & 1 & 2 & 3 & 4 & 1 & 2 & 3 & 4 & 1 & 2 & 3 & 4 & 1 & 2 & 3 & 4 & 1 & 2 & 3 & 4\\
\hline
1 & n.a. & n.a. & n.a. & n.a. & 8.4 & \textbf{12.4} & 10.9 & 11.2 & 10.7 & 11.7 & 10.7 & 10.0 & n.a. & n.a. & n.a. & n.a. & 11.1 & 10.0 & \textbf{7.6} & 10.5 & 11.6 & 9.4 & 8.7 & 9.2 &
n.a. & n.a. & n.a. & n.a. & 10.4 & 10.8 & 10.4 & 10.2 & 11.0 & 9.5 & 9.1 & 11.5\\
2 & n.a. & n.a. & n.a. & n.a. & 10.0 & \textbf{12.2} & 11.2 & 11.3 & 8.4 & 9.5 & 10.4 & 9.9 & n.a. & n.a. & n.a. & n.a. & 10.0 & 11.7 & 8.8 & 10.3 & 9.7 & 9.6 & 8.6 & 10.0 &
n.a. & n.a. & n.a. & n.a. & 9.3 & 8.7 & 9.7 & 10.8 & 9.3 & 10.7 & 9.9 & 11.2\\
3 & \textbf{26.9} & \textbf{27.4} & \textbf{20.7} & \textbf{20.2} & 9.8 & \textbf{12.4} & 10.0 & 11.7 & 8.7 & 8.7 & 9.8 & 10.5 &
\textbf{25.5} & \textbf{26.8} & \textbf{19.0} & \textbf{17.5} & 11.2 & 11.0 & 8.3 & 9.5 & 9.2 & 10.4 & 11.5 & 11.2 &
\textbf{24.5} & \textbf{24.4} & \textbf{20.7} & \textbf{18.0} & 9.7 & 8.8 & 10.8 & 10.1 & 11.6 & 8.8 & 10.7 & 11.6\\
6 & \textbf{13.0} & 11.6 & 9.8 & 11.2 & \textbf{12.0} & 11.3 & 9.5 & 11.1  & \textbf{6.5} & \textbf{7.3} & \textbf{7.6} & 8.4 &
\textbf{12.2} & 10.6 & 10.0 & 11.5 & 11.0 & 10.0 & 9.6 & 10.7 & 10.2 & 11.7 & 10.6 & 8.4 & 
11.0 & 9.5 & 10.7 & 11.5 & 10.1 & 8.9 & 10.0 & 10.8  & 8.1 & 9.8 & 9.8 & 10.9\\
8 & 10.4 & 10.7 & 10.8 & 11.5 & 11.0 & 11.3 & 11.6 & \textbf{12.3} & \textbf{12.8} & \textbf{15.3} & \textbf{13.3} & \textbf{14.2} & 
10.3 & 10.2 & 10.3 & 10.1 & 10.4 & 10.2 & 10.3 & 10.1 & 9.5 & 9.6 & 8.9 & 8.4 &
10.0 & 9.6 & 9.0 & 10.9 & 10.0 & 9.5 & 8.7 & 10.7 & 10.7 & 10.2 & 10.4 & 10.7\\
10 & 10.2 & 10.0 & 10.4 & 11.8 & 10.5 & 10.7 & 11.1 & \textbf{13.1} & \textbf{26.8} & \textbf{26.1} & \textbf{26.7} & \textbf{25.5} & 
9.2 & 8.8 & 10.0 & 10.4 & 9.4 & 8.8 & 10.1 & 10.5 & 8.9 & 10.2 & 9.8 & 10.8 &
9.8 & 8.9 & 9.7 & 10.0 & 9.7 & 8.9 & 9.7 & 10.0 & 9.4 & 9.8 & 10.7 & 8.4\\ \hline
\end{tabular}
\end{center}
\end{table}

\begin{table}[H]
\caption{\scriptsize{Empirical levels at the $5\%$ and $10\%$ significance levels for the global portmanteau test statistics $Q_M$ for $M = 1,2,3,6,8,10$. The simulated model is a strong PVAR given by DGP$_1$. The numbers of years are equal to $N = 200$, $1000$ and $5000$. $\mathrm{{Q}}_M^{1}(\nu)$ is the global portmanteau test statistics defined by~(37) in~\cite{UD09}. $\mathrm{{Q}}_M^{2}(\nu)$ is the global portmanteau test statistics defined by~(28) in~\citet{DL13}. $\mathrm{{Q}}_M^{3}(\nu)$ is the global portmanteau test statistics defined by~(\ref{QMast}).}}\label{global_test_portmanteau_strong_unc}
\begin{center}
\tiny
\setlength{\tabcolsep}{0.8mm}
\begin{tabular}{|c|ccc|ccc|ccc|ccc|ccc|ccc|} \hline
    & \multicolumn{9}{c|}{$\alpha=0.05$} & \multicolumn{9}{c|}{$\alpha=0.10$} \\
    & \multicolumn{3}{c}{$N=200$} & \multicolumn{3}{c}{$N=1000$} & \multicolumn{3}{c|}{$N=5000$}
    & \multicolumn{3}{c}{$N=200$} & \multicolumn{3}{c}{$N=1000$} & \multicolumn{3}{c|}{$N=5000$}\\ \hline
$M$ & $Q^1_M$ & $Q^2_M$ & $Q^3_M$ & $Q^1_M$ & $Q^2_M$ & $Q^3_M$ & $Q^1_M$ & $Q^2_M$ & $Q^3_M$ & $Q^1_M$ & $Q^2_M$ & $Q^3_M$ & $Q^1_M$ & $Q^2_M$ & $Q^3_M$ & $Q^1_M$ & $Q^2_M$ & $Q^3_M$\\ \hline
1 & n.a. & 3.9 & 3.9 & n.a. & 5.7 & 6.0 & n.a. & 4.6 & 5.3 & n.a. & 10.3 & 8.9 & n.a. & 11.2 & 11.5 & n.a. & 10.0 & 11.2\\
2  & \textbf{11.6} & 5.3 & 3.7 & 9.4 & 6.4 & 5.1 & \textbf{11.4} & 4.7 & 5.6 &  \textbf{19.8} & 11.4 & \textbf{7.7} & \textbf{17.7} & 11.6  & 10.1  & \textbf{20.1} & 8.5 & 10.1\\
3   & \textbf{6.9} & 5.6 & 3.6 & \textbf{7.3} & 5.7 & 5.2 & 5.7 & 5.1 & 5.1 &  \textbf{12.6} & 10.7 & 8.7 & \textbf{13.5} & \textbf{12.0} & 11.2  & 10.6 & 9.2 & 9.2 \\
6   & 6.1 & 5.9 & 3.5 & 5.0 & 5.9 & \textbf{3.1} & 4.4 & 4.4 & 4.1 & 10.2 & \textbf{12.2} & \textbf{7.8} & 9.6  & 11.3 & \textbf{6.7}  & 8.6 & 9.9 & 8.4 \\
8   & 5.0 & 5.2 & \textbf{1.3} & 5.7 & 6.4 & \textbf{4.4} & 5.4 & 5.2 & 5.0 & 10.1 & 11.6 & \textbf{3.8} & 10.7 & 11.3 & 10.4  & 10.1 & 9.5 & 9.5 \\
10  & 5.8 & 6.2 & \textbf{1.8} & 6.0 & 6.0 & 3.6 & 5.6 & 4.9 & 5.0 & \textbf{12.0} & \textbf{12.9} & \textbf{6.0} & 11.8 & 11.9 & \textbf{7.8}  & 9.8 & 9.8 & 9.5 \\ \hline
\end{tabular}
\end{center}
\end{table}

\begin{table}[H]
\caption{\scriptsize{Empirical levels at the $5\%$ and $10\%$ significance levels for the global portmanteau test statistics $Q_M$ for $M = 1,2,3,6,8,10$. The simulated model is a strong PVAR given by DGP$_3$. The numbers of years are equal to $N = 200$, $1000$ and $5000$. $\mathrm{{Q}}_M^{1}(\nu)$ is the global portmanteau test statistics defined by~(37) in~\cite{UD09}. $\mathrm{{Q}}_M^{2}(\nu)$ is the global portmanteau test statistics defined by~(28) in~\citet{DL13}. $\mathrm{{Q}}_M^{3}(\nu)$ is the global portmanteau test statistics defined by~(\ref{QMast}).}}\label{global_test_portmanteau_strong_unc_2}
\begin{center}
\tiny
\setlength{\tabcolsep}{0.8mm}
\begin{tabular}{|c|ccc|ccc|ccc|ccc|ccc|ccc|} \hline
    & \multicolumn{9}{c|}{$\alpha=0.05$} & \multicolumn{9}{c|}{$\alpha=0.10$} \\
    & \multicolumn{3}{c}{$N=200$} & \multicolumn{3}{c}{$N=1000$} & \multicolumn{3}{c|}{$N=5000$}
    & \multicolumn{3}{c}{$N=200$} & \multicolumn{3}{c}{$N=1000$} & \multicolumn{3}{c|}{$N=5000$}\\ \hline
$M$ & $Q^1_M$ & $Q^2_M$ & $Q^3_M$ & $Q^1_M$ & $Q^2_M$ & $Q^3_M$ & $Q^1_M$ & $Q^2_M$ & $Q^3_M$ & $Q^1_M$ & $Q^2_M$ & $Q^3_M$ & $Q^1_M$ & $Q^2_M$ & $Q^3_M$ & $Q^1_M$ & $Q^2_M$ & $Q^3_M$\\ \hline
1 & n.a. & 5.4 & 4.9 & n.a. & 5.3 & 4.0 & n.a. & 5.3 & 5.0 & n.a. & 11.5 & 10.8 & n.a. & 11.3 & 9.6 & n.a. & 11.3 & 9.9\\
2 & n.a. & 5.3 & 5.1 & n.a. & 4.7 & 4.4 & n.a. & 4.5 & 6.0 &  n.a. & 11.4 & 10.1 & n.a. & 10.7  & 9.6  & n.a. & 9.8 & 10.0\\
3 & \textbf{30.8} & 5.7 & 4.1 & \textbf{27.8} & 5.3 & 5.3 & \textbf{28.5} & 4.1 & \textbf{6.8} &  \textbf{48.0} & 11.4 & 10.5 & \textbf{42.1} & 9.2 & 11.5  & \textbf{43.5} & 9.4 & 10.0 \\
6 & \textbf{7.6} & \textbf{6.8} & \textbf{3.1} & 5.6 & 5.4 & 3.6 & 5.6 & 4.8 & 4.6 & \textbf{13.7} & \textbf{12.8} & \textbf{6.9} & 11.5  & 9.9 & 9.4  & 11.4 & 10.2 & \textbf{7.6} \\
8 & 6.3 & \textbf{7.2} & \textbf{13.7} & 5.1 & 5.0 & \textbf{3.3} & 5.5 & 4.9 & 5.8 & 
11.6 & \textbf{12.9} & \textbf{25.2} & 11.5 & 9.9 & 8.5  & 10.3 & 10.0 & 10.8 \\
10 & 5.1 & 6.1 & \textbf{38.6} & 4.3 & 4.4 & 4.3 & 4.5 & 4.5 & 4.3 &
\textbf{9.8} & \textbf{12.3} & \textbf{60.6} & 9.5 & 9.9 & 10.1  & 10.2 & 10.2 & 9.7 \\ \hline
\end{tabular}
\end{center}
\end{table}

\subsubsection{Constrained strong PVAR model}
We consider the constrained strong PVAR defined by DGP$_2$ and the innovation process $\bfepsilon = \{ \bfepsilon_t, t \in \mathbb{Z} \}$ is the same as in unconstrained case. The number of rejections of the null hypothesis of adequacy are reported in Tables~(\ref{test_portmanteau_strong_const}) and~(\ref{global_test_portmanteau_strong_const}). Note that the results are very similar to the unconstrained case.
\begin{table}[H]
\caption{\scriptsize{Empirical levels (number of rejections of the null hypothesis over $1000$ replications) at the $5\%$ and $10\%$ significance levels for the modified portmanteau test statistics $Q_M(\nu)$ for $\nu=1,2,3,4$ and $M = 1,2,3,6,8,10$. The simulated model is a strong PVAR given by DGP$_2$. The numbers of years are equal to $N = 200$, $1000$ and $5000$}. $\mathrm{{Q}}_M^{1}(\nu)$ is the portmanteau test statistics defined by~(35) in~\citet{UD09}. $\mathrm{{Q}}_M^{2}(\nu)$ is the portmanteau test statistics defined by~(26) in~\citet{DL13}. $\mathrm{{Q}}_M^{3}(\nu)$ is the portmanteau test statistics defined by~(\ref{QMnuast}).
}\label{test_portmanteau_strong_const}
\begin{center}
\tiny
\setlength{\tabcolsep}{0.8mm}
\begin{tabular}{|c|cccc|cccc|cccc|cccc|cccc|cccc|cccc|cccc|cccc|}
\hline
& \multicolumn{36}{c|}{$\alpha=0.05$} \\
& \multicolumn{12}{c}{$N=200$} & \multicolumn{12}{c}{$N=1000$} & \multicolumn{12}{c|}{$N=5000$}\\
\cline{2-37}
&\multicolumn{4}{c}{$Q^1_M(\nu)$} &  \multicolumn{4}{c|}{$Q^2_M(\nu)$} & \multicolumn{4}{c|}{$Q^3_M(\nu)$} & \multicolumn{4}{c}{$Q^1_M(\nu)$} &  \multicolumn{4}{c|}{$Q^2_M(\nu)$} & \multicolumn{4}{c|}{$Q^3_M(\nu)$} & \multicolumn{4}{c}{$Q^1_M(\nu)$} &  \multicolumn{4}{c|}{$Q^2_M(\nu)$} & \multicolumn{4}{c|}{$Q^3_M(\nu)$}\\
\hline
\diagbox{$M$}{$\nu$}& 1 & 2 & 3 & 4 & 1 & 2 & 3 & 4 & 1 & 2 & 3 & 4 & 1 & 2 & 3 & 4 & 1 & 2 & 3 & 4 & 1 & 2 & 3 & 4 & 1 & 2 & 3 & 4 & 1 & 2 & 3 & 4 & 1 & 2 & 3 & 4\\
\hline
1 & n.a. & n.a. & n.a. & n.a. & 4.2 & 3.6 & 5.4 & 4.2 & \textbf{3.4} & \textbf{3.2} & 4.7 & 4.2
& n.a. & n.a. & n.a. & n.a. & 5.2 & 5.0 & 5.3 & 5.4 & 3.8 & 4.2 & 5.0 & 4.2
& n.a. & n.a. & n.a. & n.a. & 5.2 & 5.2 & 5.3 & 6.0 & 5.2 & 5.3 & 4.5 & 5.3\\
2 & \textbf{10.7} & \textbf{8.9} & \textbf{8.1} & \textbf{10.8} & 5.8 & 4.1 & 6.3 & 4.8 & \textbf{3.5} & \textbf{2.5} & \textbf{3.3} & 3.7 &
\textbf{10.7} & \textbf{7.4} & \textbf{9.3} & \textbf{8.8} & 4.7 & 5.9 & 4.9 & 4.8 & 4.3 & 4.6 & 4.7 & 4.3
& \textbf{12.1} & \textbf{10.5} & \textbf{8.4} & \textbf{9.4} & 6.3 & 5.2 & 4.3 & 5.4 & 6.2 & 4.6 & 5.0 & 4.8\\
3 & \textbf{7.7} & \textbf{7.3} & \textbf{7.1} & \textbf{7.3} & 4.8 & 4.5 & \textbf{6.6} & 5.2 & \textbf{3.0} & \textbf{2.8} & \textbf{3.2} & \textbf{3.1} &
\textbf{7.5} & \textbf{8.8} & 5.5 & 6.4 & 4.5 & 4.9 & 4.3 & 5.0 & 3.9 & 5.6 & 3.6 & 4.3
& \textbf{9.2} & \textbf{6.7} & \textbf{6.9} & 6.3 & 5.3 & 5.8 & 5.1 & 5.2 & 6.1 & 4.0 & 5.1 & 3.8\\
6 & 5.9 & 5.8 & 4.8 & 6.4 & 6.2 & 5.2 & 6.2 & 6.0 & \textbf{2.3} & \textbf{1.8} & \textbf{1.6} & \textbf{1.7} &
6.3 & 5.9 & 6.8 & 5.2 & 4.3 & 5.1 & 5.1 & 5.8 & 4.5 & 3.7 & 4.4 & 3.9
& \textbf{6.5} & 5.3 & 4.4 & 5.6 & 6.0 & 5.0 & 4.5 & 4.6 & 5.4 & 4.5 & \textbf{3.5} & 4.5\\
8 & 5.9 & 4.2 & 5.9 & 5.3 & 6.2 & 5.3 & 6.2 & 5.8 & \textbf{1.3} & \textbf{0.7} & \textbf{1.6} & \textbf{1.0} &
6.2 & 5.3 & 5.8 & 5.1 & 4.4 & 5.8 & 6.0 & 4.9 & 4.8 & 3.7 & 4.6 & \textbf{3.2}
& 5.7 & 5.8 & 5.4 & 5.0 & 5.3 & 5.2 & 4.4 & 5.2 & 4.8 & 5.1 & 4.6 & 4.0\\
10 & 5.1 & 5.0 & 5.3 & 5.0 & 5.9 & 5.6 & 6.0 & \textbf{6.9} & \textbf{1.4} & \textbf{0.8} & \textbf{1.1} & \textbf{1.1} &
4.3 & 4.8 & 5.8 & 5.2 & 4.5 & 6.0 & 5.2 & 4.9 & \textbf{2.9} & \textbf{3.2} & 4.2 & \textbf{3.5}
& 4.8 & 5.6 & 5.5 & 5.3 & 5.8 & 5.4 & 5.9 & 4.2 & 4.3 & 5.2 & 5.1 & 4.2\\
\hline
& \multicolumn{36}{c|}{$\alpha=0.10$} \\
& \multicolumn{12}{c}{$N=200$} & \multicolumn{12}{c}{$N=1000$} & \multicolumn{12}{c|}{$N=5000$}\\
\cline{2-37}
&\multicolumn{4}{c}{$Q^1_M(\nu)$} &  \multicolumn{4}{c|}{$Q^2_M(\nu)$} &  \multicolumn{4}{c|}{$Q^3_M(\nu)$} & \multicolumn{4}{c}{$Q^1_M(\nu)$} &  \multicolumn{4}{c|}{$Q^2_M(\nu)$} & \multicolumn{4}{c|}{$Q^3_M(\nu)$} & \multicolumn{4}{c}{$Q^1_M(\nu)$} &  \multicolumn{4}{c|}{$Q^2_M(\nu)$} & \multicolumn{4}{c|}{$Q^3_M(\nu)$}\\
\hline
\diagbox{$M$}{$\nu$} & 1 & 2 & 3 & 4 & 1 & 2 & 3 & 4 & 1 & 2 & 3 & 4 & 1 & 2 & 3 & 4 & 1 & 2 & 3 & 4 & 1 & 2 & 3 & 4 & 1 & 2 & 3 & 4 & 1 & 2 & 3 & 4 & 1 & 2 & 3 & 4 \\
\hline
1 & n.a. & n.a. & n.a. & n.a. & 9.0 & 9.0 & 11.3 & 8.8 & 8.8 & \textbf{7.6} & 9.6 & 9.5 & n.a. & n.a. & n.a. & n.a. & 9.6 & 9.6 & 10.2 & 10.5 & \textbf{7.5} & 10.3 & 10.3 & 8.8
& n.a. & n.a. & n.a. & n.a. & 9.9 & 10.4 & 10.0 & 9.5 & 10.6 & 9.9 & 9.1 & 9.7\\
2 & \textbf{19.9} & \textbf{17.8} & \textbf{17.8} & \textbf{18.5} & 10.2 & 8.9 & 11.8 & 10.4 & \textbf{8.0} & \textbf{7.5} & \textbf{8.0} & 9.6 & \textbf{19.9} & \textbf{15.2} & \textbf{19.3} & \textbf{17.7} & 10.0 & 9.6 & 10.0 & 10.6 & 9.1 & \textbf{7.9} & 10.6 & 8.8
& \textbf{20.2} & \textbf{18.8} & \textbf{17.9} & \textbf{16.8} & 11.3 & 11.3 & 9.6 & 9.9 & 11.7 & 11.3 & 10.2 & 10.0\\
3 & \textbf{14.4} & \textbf{14.2} & \textbf{14.4} & \textbf{13.6} & 10.0 & 9.7 & \textbf{12.8} & 10.2 & \textbf{6.8} & \textbf{7.0} & \textbf{7.2} & 8.1 & \textbf{14.5} & \textbf{14.5} & 11.5 & \textbf{12.3} & 9.2 & 9.4 & 8.4 & 11.0 & 8.7 & 9.7 & \textbf{7.1} & 8.4
& 16.0 & 12.8 & 13.7 & 13.0 & 10.7 & 10.3 & 8.8 & 9.9 & 10.4 & 8.8 & 8.8 & 8.6\\
6 & \textbf{12.1} & 10.9 & 10.1 & \textbf{12.1} & 11.6 & 10.8 & 11.8 & \textbf{12.4} & \textbf{5.5} & \textbf{5.1} & \textbf{4.4} & \textbf{5.9} &
\textbf{13.0} & \textbf{13.8} & \textbf{12.7} & 10.5 & 8.8 & 10.5 & 9.0 & 10.3 & 9.0 & 9.9 & 10.1 & 8.2
& \textbf{12.8} & 10.6 & 10.6 & 9.7 & 10.9 & 9.9 & 8.8 & 9.1 & 9.9 & 8.6 & 8.8 & \textbf{8.0}\\
8 & 11.1 & 10.2 & 11.4 & 10.5 & 11.6 & 10.8 & 11.8 & 11.5 & \textbf{4.6} & \textbf{3.0} & \textbf{4.6} & \textbf{4.9} & 10.5 & \textbf{12.0} & \textbf{12.2} & 10.2 & 8.3 & 10.7 & 10.5 & 10.6 & \textbf{8.0} & 8.4 & 10.0 & 8.2
& 11.1 & 10.4 & 11.1 & 10.5 & 10.8 & 10.5 & 9.9 & 9.4 & 10.3 & 9.1 & 9.8 & 9.0\\
10 & 9.7 & 10.6 & 10.0 & 10.1 & 11.8 & 10.8 & 11.0 & 11.8 & \textbf{3.5} & \textbf{3.4} & \textbf{3.6} & \textbf{3.0} & 10.3 & 10.3 & 10.1 & 9.4 & 9.2 & 11.0 & 9.8 & 10.2 & \textbf{7.4} & \textbf{8.0} & 8.2 & \textbf{7.7}
& 10.1 & 10.4 & 10.7 & 9.9 & 10.7 & 10.3 & 10.8 & 9.4 & 8.6 & 9.3 & 9.7 & 9.4\\ \hline
\end{tabular}
\end{center}
\end{table}

\begin{table}[H]
\caption{\scriptsize{Empirical levels at the $5\%$ and $10\%$ significance levels for the global portmanteau test statistics $Q_M$ for $M = 1,2,3,6,8,10$. The simulated model is a strong PVAR given by DGP$_2$. The numbers of years are equal to $N = 200$, $1000$ and $5000$. $\mathrm{{Q}}_M^{1}(\nu)$ is the global portmanteau test statistics defined by~(37) in~\cite{UD09}. $\mathrm{{Q}}_M^{2}(\nu)$ is the global portmanteau test statistics defined by~(28) in~\citet{DL13}. $\mathrm{{Q}}_M^{3}(\nu)$ is the global portmanteau test statistics defined by~(\ref{QMast}).}}\label{global_test_portmanteau_strong_const}
\begin{center}
\tiny
\setlength{\tabcolsep}{0.8mm}
\begin{tabular}{|c|ccc|ccc|ccc|ccc|ccc|ccc|} \hline
    & \multicolumn{9}{c|}{$\alpha=0.05$} & \multicolumn{9}{c|}{$\alpha=0.10$} \\
    & \multicolumn{3}{c}{$N=200$} & \multicolumn{3}{c}{$N=1000$} & \multicolumn{3}{c|}{$N=5000$}
    & \multicolumn{3}{c}{$N=200$} & \multicolumn{3}{c}{$N=1000$} & \multicolumn{3}{c|}{$N=5000$}\\ \hline
$M$ & $Q^1_M$ & $Q^2_M$ & $Q^3_M$ & $Q^1_M$ & $Q^2_M$ & $Q^3_M$ & $Q^1_M$ & $Q^2_M$ & $Q^3_M$ & $Q^1_M$ & $Q^2_M$ & $Q^3_M$ & $Q^1_M$ & $Q^2_M$ & $Q^3_M$ & $Q^1_M$ & $Q^2_M$ & $Q^3_M$\\ \hline
1 & n.a. & 5.1 & 3.7 & n.a. & 6.0 & 4.4 & n.a. & 5.7 & 5.9 & n.a. & 10.4 & 8.3 & n.a. & 10.5 & 8.1 & n.a. & 9.9 & 11.5 \\
2 & \textbf{20.0} & 6.3 & 3.8 & \textbf{18.3} & 5.6 & 5.6 & \textbf{20.1} & 4.9 & 5.1 &\textbf{32.8} & 11.8 & 8.9 & \textbf{31.0} & 10.1  & 9.9 & \textbf{33.7} & 10.0 & 10.5\\
3 & \textbf{11.4} & 5.4 & \textbf{3.2} & \textbf{10.3} & 5.2 & 3.6 & \textbf{10.8} & 5.8 & 4.8 &\textbf{20.1} & 10.6 & \textbf{7.4} & \textbf{18.0} & 10.2 & \textbf{7.4} & \textbf{21.4} & 10.3 & 9.4\\
6 & 6.2 & \textbf{6.8} & \textbf{1.6} & \textbf{7.1} & 5.6 & \textbf{3.9} & 6.4 & 5.1 & \textbf{2.6} & \textbf{14.5} & 11.9 & \textbf{4.0} & \textbf{14.6}  & 9.2 & 8.4 & \textbf{12.7} & 10.7 & 8.6\\
8 & \textbf{8.5} & \textbf{7.1} & \textbf{2.0} & 5.9 & 5.0 & 3.6 & 5.5 & 5.0 & 4.2 &\textbf{14.2} & \textbf{13.4} & \textbf{5.0} & \textbf{13.0} & 9.6 & \textbf{7.9} & \textbf{12.9} & 10.4 & 9.4 \\
10 & \textbf{6.6} & \textbf{8.0} & \textbf{0.5} & 5.9 & 4.3 & \textbf{3.3} & 6.1 & 5.2 & 4.6 & \textbf{13.0} & \textbf{14.0} & \textbf{3.1} & 11.4 & 9.9 & \textbf{7.8} & 11.5 & 10.5 & 10.1\\ \hline
\end{tabular}
\end{center}
\end{table}

From these examples we draw the conclusion that, for strong PVAR models diagnostic checks season by season are recommended to use and the test proposed by~\citet{DL13} may be preferable to the other mentioned tests when $N$ is small ($N = 200$). Note that our proposed test can be used safely for $M$ small and $N$ large enough ($N\geq 100$).

\subsubsection{Unconstrained weak PVAR model}
\noindent We repeat the same experiments on a weak PVAR model, meaning that the stochastic process $\bfepsilon$  is defined by
\begin{align}
\label{bruitPT}
\bfepsilon_{ns+\nu} &=\bfM_{\nu}^\top
\begin{pmatrix}
\eta_{1,ns+\nu}\eta_{1,ns+\nu-1}\eta_{1,ns+\nu-2} \\
\eta_{2,ns+\nu}\eta_{2,ns+\nu-1}\eta_{2,ns+\nu-2}
\end{pmatrix}
\end{align}
where $\bfeta_{t} =(\eta_{1,t},\eta_{2,t})^\top$ iid $\mathcal{N}(\bfzero,\bfI_2)$ and $\bfM_{\nu}$ is the upper triangular matrix satisfying the equation $\bfM^{T}_{\nu} \bfM_{\nu} = \bfSigma_{\bfepsilon}(\nu):=\mathbb{E}[\bfepsilon_{ns+\nu}\bfepsilon_{ns+\nu}^\top]$.
This process can be viewed as a multivariate extension of
univariate weak noises considered in~\cite{rt1996}.

\begin{table}[H]
\caption{\scriptsize{Empirical levels (number of rejections of the null hypothesis over $1000$ replications) at the $5\%$ and $10\%$ significance levels for the modified portmanteau test statistics $Q_M(\nu)$ for $\nu=1,2,3,4$ and $M = 1,2,3,6,8,10$. The simulated model is a weak PVAR given by DGP$_1$. The numbers of years are equal to $N = 200$, $1000$ and $5000$}. $\mathrm{{Q}}_M^{1}(\nu)$ is the portmanteau test statistics defined by~(35) in~\citet{UD09}. $\mathrm{{Q}}_M^{2}(\nu)$ is the portmanteau test statistics defined by~(26) in~\citet{DL13}. $\mathrm{{Q}}_M^{3}(\nu)$ is the portmanteau test statistics defined by~(\ref{QMnuast}).
}\label{test_portmanteau_weak_unc}
\begin{center}
\tiny
\setlength{\tabcolsep}{0.8mm}
\begin{tabular}{|c|cccc|cccc|cccc|cccc|cccc|cccc|cccc|cccc|cccc|}
\hline
& \multicolumn{36}{c|}{$\alpha=0.05$} \\
& \multicolumn{12}{c}{$N=200$} & \multicolumn{12}{c}{$N=1000$} & \multicolumn{12}{c|}{$N=5000$}\\
\cline{2-37}
&\multicolumn{4}{c}{$Q^1_M(\nu)$} &  \multicolumn{4}{c|}{$Q^2_M(\nu)$} & \multicolumn{4}{c|}{$Q^3_M(\nu)$} & \multicolumn{4}{c}{$Q^1_M(\nu)$} &  \multicolumn{4}{c|}{$Q^2_M(\nu)$} & \multicolumn{4}{c|}{$Q^3_M(\nu)$} & \multicolumn{4}{c}{$Q^1_M(\nu)$} &  \multicolumn{4}{c|}{$Q^2_M(\nu)$} & \multicolumn{4}{c|}{$Q^3_M(\nu)$}\\
\hline
\diagbox{$M$}{$\nu$}& 1 & 2 & 3 & 4 & 1 & 2 & 3 & 4 & 1 & 2 & 3 & 4 & 1 & 2 & 3 & 4 & 1 & 2 & 3 & 4 & 1 & 2 & 3 & 4 & 1 & 2 & 3 & 4 & 1 & 2 & 3 & 4 & 1 & 2 & 3 & 4\\
\hline
1 & n.a. & n.a. & n.a. & n.a. & \textbf{44.9} & \textbf{20.8} & \textbf{32.3} & \textbf{46.0} & 4.4 & 4.5 & \textbf{2.7} & 4.3 & n.a. & n.a. & n.a. & n.a. & \textbf{50.0} & \textbf{24.8} & \textbf{38.1} & \textbf{55.4 }& 4.5 & 4.7 & 4.8 & 5.1 & n.a. & n.a. & n.a. & n.a. & \textbf{51.1} & \textbf{26.9} & \textbf{43.5} & \textbf{55.8} & \textbf{3.3} & 6.1 & 5.1 & 4.0 \\
2 & \textbf{60.1} & \textbf{29.8} & \textbf{38.2} & \textbf{53.4} & \textbf{46.6} & \textbf{26.5} & \textbf{36.6} & \textbf{47.7} & \textbf{1.8} & 3.6 & \textbf{3.2} & 3.7 & \textbf{68.6} & \textbf{33.3 }& \textbf{44.2} & \textbf{59.5} & \textbf{57.3} & \textbf{32.9} & \textbf{40.9} & \textbf{55.8} & 4.3 & 4.7 & \textbf{3.5} & \textbf{3.2} & \textbf{68.9} & \textbf{34.7} & \textbf{47.1} & \textbf{58.9} & \textbf{58.3} & \textbf{36.0} & \textbf{45.1} & \textbf{55.9} & 4.3 & 4.5 & 6.4 & 4.8\\
3 & \textbf{47.3} & \textbf{23.7} & \textbf{30.3} & \textbf{42.1} & \textbf{43.9} & \textbf{21.3 }& \textbf{32.0} & \textbf{42.8} & \textbf{2.2} & \textbf{1.8} & \textbf{1.9} & \textbf{2.6} & \textbf{55.5} & \textbf{29.3} & \textbf{35.8} & \textbf{50.4} & \textbf{52.9} & \textbf{27.6} & \textbf{35.0} & \textbf{51.7} & 4.2 & \textbf{3.4} & \textbf{2.3} & \textbf{2.9}& \textbf{57.3} & \textbf{29.8} & \textbf{37.1} & \textbf{52.5} & \textbf{56.7} & \textbf{29.6} & \textbf{40.8} & \textbf{52.6} & 5.2 & 3.8 & 4.4 & 5.3 \\
6 & \textbf{33.5} & \textbf{17.5} & \textbf{22.3} & \textbf{32.8} & \textbf{34.7} & \textbf{17.1} & \textbf{25.5} & \textbf{34.5} & \textbf{0.7} & \textbf{0.2} & \textbf{0.7} & \textbf{0.7} & \textbf{41.5} & \textbf{20.4} & \textbf{26.7} & \textbf{41.9} & \textbf{42.8} & \textbf{22.0} & \textbf{25.7} & \textbf{40.1} & \textbf{2.8} & \textbf{2.0 }& \textbf{2.7} & \textbf{3.0}& \textbf{44.1} & \textbf{21.3} & \textbf{29.8} & \textbf{40.3} & \textbf{45.4} & \textbf{21.2} & \textbf{28.9} & \textbf{40.4} & 4.5 & 4.3 & 4.7 & 4.0 \\
8 & \textbf{33.6} & \textbf{17.7} & \textbf{21.0} & \textbf{30.9} & \textbf{30.2} & \textbf{16.5} & \textbf{22.4} & \textbf{31.0} & \textbf{1.1} & \textbf{1.2} & \textbf{1.3} & \textbf{1.1 }& \textbf{39.2 }& \textbf{16.6} & \textbf{22.8} & \textbf{31.8} & \textbf{38.4} & \textbf{19.9} & \textbf{23.2} & \textbf{35.4} & \textbf{3.5} & \textbf{1.5} & \textbf{2.5 }& \textbf{2.3}& \textbf{40.4} & \textbf{19.2} & \textbf{26.5} & \textbf{36.0} & \textbf{39.6} & \textbf{19.1} & \textbf{26.2} & \textbf{36.2} & 4.1 & \textbf{3.3} & 4.2 & \textbf{3.5}\\
10 & \textbf{28.3} & \textbf{15.4} & \textbf{20.2} & \textbf{26.9} & \textbf{29.3} & \textbf{15.2} & \textbf{20.8} & \textbf{28.8} & \textbf{3.3} & 5.5 & 5.2 & 4.0 & \textbf{33.3} & \textbf{18.1} & \textbf{22.5 }& \textbf{34.9} & \textbf{33.7} & \textbf{17.3} & \textbf{22.8} & \textbf{32.9} & \textbf{1.1} & \textbf{1.5} & \textbf{1.8} & \textbf{1.5}& \textbf{34.7} & \textbf{19.9} & \textbf{23.6} & \textbf{34.0} & \textbf{37.3} & \textbf{16.5} & \textbf{23.2} & \textbf{32.6} & 3.8 & 3.8 & \textbf{3.1} & \textbf{2.7} \\
\hline
& \multicolumn{36}{c|}{$\alpha=0.10$} \\
& \multicolumn{12}{c}{$N=200$} & \multicolumn{12}{c}{$N=1000$} & \multicolumn{12}{c|}{$N=5000$}\\
\cline{2-37}
&\multicolumn{4}{c}{$Q^1_M(\nu)$} &  \multicolumn{4}{c|}{$Q^2_M(\nu)$} &  \multicolumn{4}{c|}{$Q^3_M(\nu)$} & \multicolumn{4}{c}{$Q^1_M(\nu)$} &  \multicolumn{4}{c|}{$Q^2_M(\nu)$} & \multicolumn{4}{c|}{$Q^3_M(\nu)$} & \multicolumn{4}{c}{$Q^1_M(\nu)$} &  \multicolumn{4}{c|}{$Q^2_M(\nu)$} & \multicolumn{4}{c|}{$Q^3_M(\nu)$}\\
\hline
\diagbox{$M$}{$\nu$} & 1 & 2 & 3 & 4 & 1 & 2 & 3 & 4 & 1 & 2 & 3 & 4 & 1 & 2 & 3 & 4 & 1 & 2 & 3 & 4 & 1 & 2 & 3 & 4 & 1 & 2 & 3 & 4 & 1 & 2 & 3 & 4 & 1 & 2 & 3 & 4 \\
\hline
1 & n.a. & n.a. & n.a. & n.a. & \textbf{54.8} & \textbf{27.3} & \textbf{40.9} & \textbf{56.0} & 10.9 & 11.4 & 8.7 & 9.9 & n.a. & n.a. & n.a. & n.a. & \textbf{58.5} & \textbf{35.7} & \textbf{47.9} & \textbf{63.5} & 9.4 & 8.8 & 11.8 & 11.6 & n.a. & n.a. & n.a. & n.a. & \textbf{60.8} & \textbf{35.7} & \textbf{52.5} & \textbf{65.1} & 8.0 & 11.6 & 10.1 & 8.7 \\
2 & \textbf{69.3} & \textbf{38.2} & \textbf{46.8} & \textbf{61.6} & \textbf{56.7} & \textbf{34.2} & \textbf{45.7} & \textbf{56.3} & \textbf{6.7} & 8.8 & 8.6 & 9.9 & \textbf{77.1} & \textbf{43.0} & \textbf{53.2} & \textbf{68.7} & \textbf{65.5} & \textbf{41.7} & \textbf{49.3} & \textbf{65.20} & 9.4 & 10.1 & 8.9 & \textbf{7.9} & \textbf{78.4} & \textbf{43.9} & \textbf{56.9} & \textbf{67.7} & \textbf{66.2} & \textbf{44.5} & \textbf{54.5} & \textbf{64.6} & 8.3 & 9.6 & 11.5 & 10.6 \\
3 & \textbf{55.7} & \textbf{31.4} & \textbf{39.4} & \textbf{51.2} & \textbf{52.1} & \textbf{28.5} & \textbf{39.8} & \textbf{51.4} & 8.1 & \textbf{5.6} & \textbf{6.6} & \textbf{7.7} & \textbf{64.4} & \textbf{38.9} & \textbf{46.6} & \textbf{59.0} & \textbf{62.5} & \textbf{36.5} & \textbf{44.2 }& \textbf{59.9} & 8.7 & 8.1 & \textbf{7.2} & \textbf{7.4}& \textbf{65.6} & \textbf{38.8} & \textbf{46.0} & \textbf{62.7} & \textbf{63.7} & \textbf{39.2} & \textbf{50.5} & \textbf{61.5} & 9.5 & 8.1 & 9.4 & 10.9\\
6 & \textbf{43.5} & \textbf{26.4} & \textbf{29.9} & \textbf{42.9} & \textbf{43.4} & \textbf{25.0} & \textbf{33.4} & \textbf{43.8} & \textbf{3.1} & \textbf{2.5} & \textbf{3.2} & \textbf{2.5} & \textbf{50.7} & \textbf{29.3} & \textbf{34.9} & \textbf{51.4} & \textbf{51.7} & \textbf{32.1} & \textbf{34.6} & \textbf{48.7} & \textbf{7.6} & \textbf{6.2 }& \textbf{6.4} & \textbf{7.9}& \textbf{55.1} & \textbf{30.2} & \textbf{38.4} & \textbf{50.1} & \textbf{54.9} & \textbf{30.3} & \textbf{38.7} & \textbf{50.4} & 9.6 & 8.9 & 9.7 & 9.1\\
8 & \textbf{42.0} & \textbf{24.0} & \textbf{28.3} & \textbf{41.0} & \textbf{39.8} & \textbf{24.1} & \textbf{30.1} & \textbf{40.8} & \textbf{3.5} & \textbf{4.6} & \textbf{4.4} & \textbf{3.7} & \textbf{47.2} & \textbf{25.5} & \textbf{31.2} & \textbf{42.1} & \textbf{47.5} & \textbf{28.1} & \textbf{31.6} & \textbf{45.7} & 9.0 & \textbf{4.7} & \textbf{5.9} & \textbf{5.4}& \textbf{51.0} & \textbf{31.2} & \textbf{37.7} & \textbf{46.5} & \textbf{50.0} & \textbf{29.5} & \textbf{36.6} & \textbf{46.2} & 8.7 & \textbf{7.3} & 9.9 & 8.2\\
10 & \textbf{36.0} & \textbf{23.7} & \textbf{27.3} & \textbf{36.0} & \textbf{36.7} & \textbf{23.3} & \textbf{28.0} & \textbf{36.5} & 8.5 & \textbf{12.5} & 11.8 & 9.7 & \textbf{42.5} & \textbf{25.4} & \textbf{31.0} & \textbf{44.5} & \textbf{43.5} & \textbf{26.5} & \textbf{30.1} & \textbf{41.5} & \textbf{4.8} & \textbf{5.8} & \textbf{4.8} & \textbf{4.6}& \textbf{45.8} & \textbf{29.6} & \textbf{33.4} & \textbf{43.8} & \textbf{48.2} & \textbf{26.5} & \textbf{34.9} & \textbf{43.1} & \textbf{7.5} & 8.8 & \textbf{7.4} & \textbf{7.9}\\ \hline
\end{tabular}
\end{center}
\end{table}

\begin{table}[H]
\caption{\scriptsize{Empirical levels (number of rejections of the null hypothesis over $1000$ replications) at the $5\%$ and $10\%$ significance levels for the modified portmanteau test statistics $Q_M(\nu)$ for $\nu=1,2,3,4$ and $M = 1,2,3,6,8,10$. The simulated model is a weak PVAR given by DGP$_3$. The numbers of years are equal to $N = 200$, $1000$ and $5000$}. $\mathrm{{Q}}_M^{1}(\nu)$ is the portmanteau test statistics defined by~(35) in~\citet{UD09}. $\mathrm{{Q}}_M^{2}(\nu)$ is the portmanteau test statistics defined by~(26) in~\citet{DL13}. $\mathrm{{Q}}_M^{3}(\nu)$ is the portmanteau test statistics defined by~(\ref{QMnuast}).
}\label{test_portmanteau_weak_unc}
\begin{center}
\tiny
\setlength{\tabcolsep}{0.8mm}
\begin{tabular}{|c|cccc|cccc|cccc|cccc|cccc|cccc|cccc|cccc|cccc|}
\hline
& \multicolumn{36}{c|}{$\alpha=0.05$} \\
& \multicolumn{12}{c}{$N=200$} & \multicolumn{12}{c}{$N=1000$} & \multicolumn{12}{c|}{$N=5000$}\\
\cline{2-37}
&\multicolumn{4}{c}{$Q^1_M(\nu)$} &  \multicolumn{4}{c|}{$Q^2_M(\nu)$} & \multicolumn{4}{c|}{$Q^3_M(\nu)$} & \multicolumn{4}{c}{$Q^1_M(\nu)$} &  \multicolumn{4}{c|}{$Q^2_M(\nu)$} & \multicolumn{4}{c|}{$Q^3_M(\nu)$} & \multicolumn{4}{c}{$Q^1_M(\nu)$} &  \multicolumn{4}{c|}{$Q^2_M(\nu)$} & \multicolumn{4}{c|}{$Q^3_M(\nu)$}\\
\hline
\diagbox{$M$}{$\nu$}& 1 & 2 & 3 & 4 & 1 & 2 & 3 & 4 & 1 & 2 & 3 & 4 & 1 & 2 & 3 & 4 & 1 & 2 & 3 & 4 & 1 & 2 & 3 & 4 & 1 & 2 & 3 & 4 & 1 & 2 & 3 & 4 & 1 & 2 & 3 & 4\\
\hline
1 & n.a. & n.a. & n.a. & n.a. & \textbf{19.4} & \textbf{13.3} & \textbf{15.3} & \textbf{22.8} & 4.9 & 4.5 & 5.9 & 4.4 & n.a. & n.a. & n.a. & n.a. & \textbf{24.4} & \textbf{16.7} & \textbf{15.8} & \textbf{29.7}
& 6.0 & 4.9 & 5.9 & 4.3 & n.a. & n.a. & n.a. & n.a. & \textbf{23.1} & \textbf{18.1} & \textbf{16.7} & \textbf{29.0} & 5.3 & 5.7 & 5.8 & 5.8 \\
2 & n.a. & n.a. & n.a. & n.a. & \textbf{13.4} & \textbf{9.7} & \textbf{11.9} & \textbf{20.7} & 3.6 & \textbf{2.8} & \textbf{3.4} & 3.8 & n.a. & n.a. & n.a. & n.a. & \textbf{18.8} & \textbf{10.3} & \textbf{12.0} & \textbf{27.0} & 4.8 & 4.3 & \textbf{3.5} & 4.9 & n.a. & n.a. & n.a. & n.a. & \textbf{17.5} & \textbf{11.9} & \textbf{11.7} & \textbf{27.5} & 5.5 & 5.9 & 4.9 & 4.4\\
3 & \textbf{27.5} & \textbf{21.0} & \textbf{18.3} & \textbf{23.6} & \textbf{14.1} & \textbf{11.3}& \textbf{11.7} & \textbf{16.6} & \textbf{2.5} & \textbf{2.0} & \textbf{2.0} & \textbf{3.0} & \textbf{34.0} & \textbf{24.1} & \textbf{19.1} & \textbf{29.4} & \textbf{17.8} & \textbf{12.6} & \textbf{10.8} & \textbf{23.0} & 4.8 & 4.1 & 3.6 & 5 & \textbf{33.2} & \textbf{25.8} & \textbf{17.8} & \textbf{28.1} & \textbf{18.2} & \textbf{13.0} & \textbf{10.5} & \textbf{21.4} & 4.9 & 4.6 & 4.4 & 4.6 \\
6 & \textbf{13.5} & \textbf{11.3} & \textbf{10.7} & \textbf{15.5} & \textbf{12.7} & \textbf{11.1} & \textbf{10.7} & \textbf{15.1} & \textbf{1.0} & \textbf{0.6} & \textbf{1.9} & \textbf{2.2} & \textbf{14.2} & \textbf{12.3} & \textbf{8.8} & \textbf{15.1} & \textbf{13.7} & \textbf{11.9} & \textbf{8.5} & \textbf{14.6} & \textbf{2.7} & \textbf{1.6 }& \textbf{1.8} & \textbf{2.0}& \textbf{14.9} & \textbf{11.0} & \textbf{10.2} & \textbf{15.6} & \textbf{14.0} & \textbf{10.4} & \textbf{10.1} & \textbf{15.1} & 4.4 & 4.4 & \textbf{2.8} & 4.4 \\
8 & \textbf{12.7} & \textbf{11.4} & \textbf{10.9} & \textbf{13.0} & \textbf{12.8} & \textbf{11.6} & \textbf{11.1} & \textbf{13.3} & \textbf{3.5} & \textbf{3.0} & 3.6 & 4.2 & \textbf{12.0 }& \textbf{10.9} & \textbf{8.8} & \textbf{14.4} & \textbf{12.0} & \textbf{10.9} & \textbf{8.7} & \textbf{14.4} & \textbf{2.2} & \textbf{2.0} & \textbf{2.3 }& \textbf{2.1}& \textbf{11.8} & \textbf{11.6} & \textbf{9.0} & \textbf{14.0} & \textbf{11.6} & \textbf{11.3} & \textbf{8.8} & \textbf{13.9} & 3.8 & 4.3 & \textbf{3.1} & \textbf{3.1}\\
10 & \textbf{12.4} & \textbf{12.8} & \textbf{11.1} & \textbf{11.6} & \textbf{13.4} & \textbf{13.4} & \textbf{11.7} & \textbf{11.8} & 3.9 & 4.4 & 4.4 & 5.0 & \textbf{12.7} & \textbf{11.7} & \textbf{9.1 }& \textbf{12.7} & \textbf{12.7} & \textbf{11.7} & \textbf{9.2} & \textbf{12.8} & \textbf{1.6} & \textbf{1.7} & \textbf{2.4} & \textbf{1.8}& \textbf{10.9} & \textbf{10.4} & \textbf{7.7} & \textbf{12.8} & \textbf{10.9} & \textbf{10.3} & \textbf{7.7} & \textbf{12.7} & \textbf{2.8} & \textbf{2.4} & \textbf{2.8} & 3.6 \\
\hline
& \multicolumn{36}{c|}{$\alpha=0.10$} \\
& \multicolumn{12}{c}{$N=200$} & \multicolumn{12}{c}{$N=1000$} & \multicolumn{12}{c|}{$N=5000$}\\
\cline{2-37}
&\multicolumn{4}{c}{$Q^1_M(\nu)$} &  \multicolumn{4}{c|}{$Q^2_M(\nu)$} &  \multicolumn{4}{c|}{$Q^3_M(\nu)$} & \multicolumn{4}{c}{$Q^1_M(\nu)$} &  \multicolumn{4}{c|}{$Q^2_M(\nu)$} & \multicolumn{4}{c|}{$Q^3_M(\nu)$} & \multicolumn{4}{c}{$Q^1_M(\nu)$} &  \multicolumn{4}{c|}{$Q^2_M(\nu)$} & \multicolumn{4}{c|}{$Q^3_M(\nu)$}\\
\hline
\diagbox{$M$}{$\nu$} & 1 & 2 & 3 & 4 & 1 & 2 & 3 & 4 & 1 & 2 & 3 & 4 & 1 & 2 & 3 & 4 & 1 & 2 & 3 & 4 & 1 & 2 & 3 & 4 & 1 & 2 & 3 & 4 & 1 & 2 & 3 & 4 & 1 & 2 & 3 & 4 \\
\hline
1 & n.a. & n.a. & n.a. & n.a. & \textbf{29.0} & \textbf{20.4} & \textbf{23.3} & \textbf{31.9} & 11.8 & 10.3 & \textbf{12.5} & 10.1 & n.a. & n.a. & n.a. & n.a. & \textbf{32.8} & \textbf{24.2} & \textbf{23.5} & \textbf{39.2} & 11.5 & 10.9 & 10.9 & 10.9 & n.a. & n.a. & n.a. & n.a. & \textbf{31.7} & \textbf{28.4} & \textbf{24.4} & \textbf{39.0} & 9.3 & 11.1 & 11.5 & 10.5 \\
2 & n.a. & n.a. & n.a. & n.a. & \textbf{20.8} & \textbf{16.6} & \textbf{20.8} & \textbf{29.3} & 9.2 & 8.1 & 8.7 & 10.1 & n.a. & n.a. & n.a. & n.a. & \textbf{27.6} & \textbf{18.7} & \textbf{18.6} & \textbf{36.4} & 10.6 & 8.9 & 8.6 & 9.6 & n.a. & n.a. & n.a. & n.a. & \textbf{25.6} & \textbf{19.4} & \textbf{20.3} & \textbf{36.3} & 10.5 & 10.6 & 9.3 & 8.9 \\
3 & \textbf{39.6} & \textbf{34.6} & \textbf{28.2} & \textbf{33.0} & \textbf{22.4} & \textbf{16.8} & \textbf{17.4} & \textbf{23.6} & 9.6 & \textbf{7.5} & \textbf{7.5} & 8.4 & \textbf{46.1} & \textbf{34.3} & \textbf{29.2} & \textbf{40.6} & \textbf{27.7} & \textbf{19.5} & \textbf{18.3 }& \textbf{31.0} & 9.8 & 9.0 & 9.6 & 11.1 & \textbf{44.6} & \textbf{37.7} & \textbf{30.3} & \textbf{38.4} & \textbf{27.6} & \textbf{19.5} & \textbf{17.1} & \textbf{29.1} & 9.2 & 9.9 & 8.4 & 9.3\\
6 & \textbf{21.3} & \textbf{18.1} & \textbf{18.0} & \textbf{22.8} & \textbf{20.5} & \textbf{18.1} & \textbf{17.8} & \textbf{22.3} & \textbf{4.5} & \textbf{3.4} & \textbf{5.4} & \textbf{6.2} & \textbf{21.8} & \textbf{19.0} & \textbf{15.6} & \textbf{25.0} & \textbf{20.3} & \textbf{18.7} & \textbf{14.7} & \textbf{24.3} & \textbf{7.7} & \textbf{6.5} & \textbf{5.4} & \textbf{6.6} & \textbf{24.5} & \textbf{19.3} & \textbf{17.5} & \textbf{24.8} & \textbf{22.4} & \textbf{18.7} & \textbf{16.9} & \textbf{23.8} & 8.8 & 9.7 & \textbf{7.9} & 8.4\\
8 & \textbf{19.5} & \textbf{17.0} & \textbf{16.0} & \textbf{19.6} & \textbf{19.9} & \textbf{17.8} & \textbf{16.6} & \textbf{20.6} & 8.6 & 9.2 & \textbf{8.2} & 10.7 & \textbf{18.4} & \textbf{17.5} & \textbf{14.9} & \textbf{21.6} & \textbf{18.5} & \textbf{17.6} & \textbf{14.8} & \textbf{21.6} & \textbf{7.5} & \textbf{5.9} & \textbf{6.0} & \textbf{6.9} & \textbf{21.4} & \textbf{19.2} & \textbf{15.2} & \textbf{21.5} & \textbf{21.4} & \textbf{19.0} & \textbf{15.2} & \textbf{21.4} & 8.7 & 8.8 & 8.8 & \textbf{8.0} \\
10 & \textbf{19.8} & \textbf{18.6} & \textbf{16.8} & \textbf{17.5} & \textbf{21.0} & \textbf{19.1} & \textbf{17.7} & \textbf{18.2} & 11.8 & \textbf{12.7} & 11.7 & \textbf{12.6} & \textbf{20.2} & \textbf{18.2} & \textbf{16.1} & \textbf{20.8} & \textbf{20.2} & \textbf{18.4} & \textbf{16.1} & \textbf{20.9} & \textbf{6.2} & \textbf{6.3} & \textbf{5.4} & \textbf{5.7} & \textbf{18.2} & \textbf{19.1} & \textbf{14.2} & \textbf{18.9} & \textbf{18.2} & \textbf{19.1} & \textbf{14.3} & \textbf{19.0} & 8.7 & \textbf{7.3} & \textbf{7.7} & \textbf{6.8}\\ \hline
\end{tabular}
\end{center}
\end{table}

\begin{table}[H]
\caption{\scriptsize{Empirical levels at the $5\%$ and $10\%$ significance levels for the global portmanteau test statistics $Q_M$ for $M = 1,2,3,6,8,10$. The simulated model is a weak PVAR given by DGP$_1$. The numbers of years are equal to $N = 200$ and $1000$. $\mathrm{{Q}}_M^{1}(\nu)$ is the global portmanteau test statistics defined by~(37) in~\cite{UD09}. $\mathrm{{Q}}_M^{2}(\nu)$ is the global portmanteau test statistics defined by~(28) in~\citet{DL13}. $\mathrm{{Q}}_M^{3}(\nu)$ is the global portmanteau test statistics defined by~(\ref{QMast}).}}\label{global_test_portmanteau_weak_unc}
\begin{center}
\tiny
\setlength{\tabcolsep}{0.8mm}
\begin{tabular}{|c|ccc|ccc|ccc|ccc|ccc|ccc|} \hline
    & \multicolumn{9}{c|}{$\alpha=0.05$} & \multicolumn{9}{c|}{$\alpha=0.10$} \\
    & \multicolumn{3}{c}{$N=200$} & \multicolumn{3}{c}{$N=1000$} & \multicolumn{3}{c|}{$N=5000$}
    & \multicolumn{3}{c}{$N=200$} & \multicolumn{3}{c}{$N=1000$} & \multicolumn{3}{c|}{$N=5000$}\\ \hline
$M$ & $Q^1_M$ & $Q^2_M$ & $Q^3_M$ & $Q^1_M$ & $Q^2_M$ & $Q^3_M$ & $Q^1_M$ & $Q^2_M$ & $Q^3_M$ & $Q^1_M$ & $Q^2_M$ & $Q^3_M$ & $Q^1_M$ & $Q^2_M$ & $Q^3_M$ & $Q^1_M$ & $Q^2_M$ & $Q^3_M$\\ \hline
1  & n.a. & \textbf{70.7} & 3.9 & n.a. & \textbf{81.1} & 4.1 & n.a. & \textbf{83.6} & \textbf{3.5} & n.a. & \textbf{79.7} & 10.7 & n.a. & \textbf{87.3} & 9.9 & n.a. &  \textbf{88.8} & \textbf{7.9}\\
2  & \textbf{87.5} & \textbf{77.2} & \textbf{2.2} & \textbf{92.6} & \textbf{86.7} & \textbf{2.7} & \textbf{92.7} & \textbf{89.1} & 4.7 & \textbf{92.3} & \textbf{84.3} & \textbf{7.4} & \textbf{95.4} & \textbf{90.5} & \textbf{7.2} & \textbf{95.1} & \textbf{92.7} & 10.2\\
3  & \textbf{71.0} & \textbf{70.9} & \textbf{0.8} & \textbf{81.9} & \textbf{81.4} & \textbf{3.2} & \textbf{84.9} & \textbf{83.7} & 4.0 & \textbf{77.4} & \textbf{77.2} & \textbf{4.7} & \textbf{86.8 }&\textbf{ 86.7} & \textbf{6.6} & \textbf{90.6} & \textbf{88.5} & 9.2 \\
6  & \textbf{50.1} & \textbf{52.8} & \textbf{0.2 }& \textbf{65.4} & \textbf{64.7} & \textbf{1.9} & \textbf{70.7} & \textbf{69.2} & 4.1 & \textbf{59.3} & \textbf{61.2} & \textbf{1.2} & \textbf{73.4} & \textbf{73.2} & \textbf{5.6} & \textbf{78.6} & \textbf{79.6} & 9.3\\
8  & \textbf{48.7} & \textbf{44.4} & \textbf{0.9} & \textbf{54.9 }& \textbf{57.2} & \textbf{1.7} & \textbf{65.1} & \textbf{63.7} & 3.9 & \textbf{58.8} & \textbf{54.2} & \textbf{3.9} & \textbf{63.8} & \textbf{66.3} & \textbf{6.9} & \textbf{74.6} & \textbf{74.5} & 8.5\\
10  & \textbf{42.4} & \textbf{41.0} & 4.0 & \textbf{54.4 }& \textbf{53.5} & \textbf{0.8} & \textbf{57.4} & \textbf{58.6} & \textbf{2.8} & \textbf{52.3} & \textbf{51.2} & 10.7 & \textbf{63.6} & \textbf{62.9} & \textbf{4.9} & \textbf{69.8} & \textbf{69.7} & \textbf{6.9}\\ \hline
\end{tabular}
\end{center}
\end{table}

\begin{table}[H]
\caption{\scriptsize{Empirical levels at the $5\%$ and $10\%$ significance levels for the global portmanteau test statistics $Q_M$ for $M = 1,2,3,6,8,10$. The simulated model is a weak PVAR given by DGP$_3$. The numbers of years are equal to $N = 200$ and $1000$. $\mathrm{{Q}}_M^{1}(\nu)$ is the global portmanteau test statistics defined by~(37) in~\cite{UD09}. $\mathrm{{Q}}_M^{2}(\nu)$ is the global portmanteau test statistics defined by~(28) in~\citet{DL13}. $\mathrm{{Q}}_M^{3}(\nu)$ is the global portmanteau test statistics defined by~(\ref{QMast}).}}\label{global_test_portmanteau_weak_unc}
\begin{center}
\tiny
\setlength{\tabcolsep}{0.8mm}
\begin{tabular}{|c|ccc|ccc|ccc|ccc|ccc|ccc|} \hline
    & \multicolumn{9}{c|}{$\alpha=0.05$} & \multicolumn{9}{c|}{$\alpha=0.10$} \\
    & \multicolumn{3}{c}{$N=200$} & \multicolumn{3}{c}{$N=1000$} & \multicolumn{3}{c|}{$N=5000$}
    & \multicolumn{3}{c}{$N=200$} & \multicolumn{3}{c}{$N=1000$} & \multicolumn{3}{c|}{$N=5000$}\\ \hline
$M$ & $Q^1_M$ & $Q^2_M$ & $Q^3_M$ & $Q^1_M$ & $Q^2_M$ & $Q^3_M$ & $Q^1_M$ & $Q^2_M$ & $Q^3_M$ & $Q^1_M$ & $Q^2_M$ & $Q^3_M$ & $Q^1_M$ & $Q^2_M$ & $Q^3_M$ & $Q^1_M$ & $Q^2_M$ & $Q^3_M$\\ \hline
1  & n.a. & \textbf{34.4} & 3.7 & n.a. & \textbf{42.8} & 4.6 & n.a. & \textbf{44.6} & 4.8 & n.a. & \textbf{44.5} & 10.8 & n.a. & \textbf{54.7} & 11.7 & n.a. &  \textbf{55.9} & 8.7\\
2  & n.a. & \textbf{27.0} & 3.9 & n.a. & \textbf{33.0} & 3.8 & n.a. & \textbf{35.0} & 4.3 & n.a. & \textbf{37.9} & 8.7 & n.a. & \textbf{42.9} & 8.6 & n.a. & \textbf{46.2} & 9.4\\
3  & \textbf{54.5} & \textbf{23.8} & \textbf{1.8} & \textbf{59.6} & \textbf{28.6} & 4.5 & \textbf{63.4} & \textbf{30.0} & 3.8 & \textbf{67.5} & \textbf{33.4} & \textbf{7.6} & \textbf{72.8 } &\textbf{ 38.9} & 10.0 & \textbf{75.1} & \textbf{41.2} & \textbf{6.7} \\
6  & \textbf{19.6} & \textbf{18.4} & \textbf{0.6} & \textbf{20.7} & \textbf{19.1} & \textbf{2.4} & \textbf{23.4} & \textbf{20.8} & 4.6 & \textbf{29.2} & \textbf{28.3} & \textbf{3.5} & \textbf{30.4} & \textbf{28.9} & \textbf{5.3} & \textbf{34.5} & \textbf{41.2} & 9.4\\
8  & \textbf{15.7} & \textbf{17.3} & \textbf{3.2} & \textbf{17.6 }& \textbf{17.6} & \textbf{1.9} & \textbf{19.7} & \textbf{19.3} & \textbf{3.4} & \textbf{24.3} & \textbf{25.5} & 10.1 & \textbf{25.9} & \textbf{25.9} & \textbf{6.3} & \textbf{29.1} & \textbf{28.8} & 8.2\\
10  & \textbf{15.3} & \textbf{17.1} & 4.3 & \textbf{16.2 }& \textbf{16.7} & \textbf{2.1} & \textbf{16.5} & \textbf{16.5} & \textbf{3.4} & \textbf{22.5} & \textbf{25.2} & \textbf{12.5} & \textbf{24.1} & \textbf{24.3} & \textbf{5.5} & \textbf{24.8} & \textbf{24.8} & 8.1\\ \hline
\end{tabular}
\end{center}
\end{table}

\subsubsection{Constrained weak PVAR model}

\begin{table}[H]
\caption{\scriptsize{Empirical levels (number of rejections of the null hypothesis over $1000$ replications) at the $5\%$ and $10\%$ significance levels for the modified portmanteau test statistics $Q_M(\nu)$ for $\nu=1,2,3,4$ and $M = 1,2,3,6,8,10$. The simulated model is a weak PVAR given by DGP$_2$. The numbers of years are equal to $N = 200$, $1000$ and $5000$}. $\mathrm{{Q}}_M^{1}(\nu)$ is the portmanteau test statistics defined by~(35) in~\citet{UD09}. $\mathrm{{Q}}_M^{2}(\nu)$ is the portmanteau test statistics defined by~(26) in~\citet{DL13}. $\mathrm{{Q}}_M^{3}(\nu)$ is the portmanteau test statistics defined by~(\ref{QMnuast}).
}\label{test_portmanteau_weak_const}
\centering
\tiny
\setlength{\tabcolsep}{0.8mm}
\begin{tabular}{|c|cccc|cccc|cccc|cccc|cccc|cccc|cccc|cccc|cccc|}
\hline
& \multicolumn{36}{c|}{$\alpha=0.05$} \\
& \multicolumn{12}{c}{$N=200$} & \multicolumn{12}{c}{$N=1000$} & \multicolumn{12}{c|}{$N=5000$}\\
\cline{2-37}
&\multicolumn{4}{c}{$Q^1_M(\nu)$} &  \multicolumn{4}{c|}{$Q^2_M(\nu)$} & \multicolumn{4}{c|}{$Q^3_M(\nu)$} & \multicolumn{4}{c}{$Q^1_M(\nu)$} &  \multicolumn{4}{c|}{$Q^2_M(\nu)$} & \multicolumn{4}{c|}{$Q^3_M(\nu)$} & \multicolumn{4}{c}{$Q^1_M(\nu)$} &  \multicolumn{4}{c|}{$Q^2_M(\nu)$} & \multicolumn{4}{c|}{$Q^3_M(\nu)$}\\
\hline
\diagbox{$M$}{$\nu$}& 1 & 2 & 3 & 4 & 1 & 2 & 3 & 4 & 1 & 2 & 3 & 4 & 1 & 2 & 3 & 4 & 1 & 2 & 3 & 4 & 1 & 2 & 3 & 4 & 1 & 2 & 3 & 4 & 1 & 2 & 3 & 4 & 1 & 2 & 3 & 4\\
\hline
1 & n.a. & n.a. & n.a. & n.a. & \textbf{59.4} & \textbf{48.0} & \textbf{51.0} & \textbf{54.9} & \textbf{2.1} & \textbf{2.7} & \textbf{2.3} & \textbf{3.2}
& n.a. & n.a. & n.a. & n.a. & \textbf{59.5} & \textbf{56.8} & \textbf{60.6} & \textbf{60.1} & 4.1 & 3.9 & \textbf{3.5} & 4.1
& n.a. & n.a. & n.a. & n.a. & \textbf{64.5} & \textbf{59.2} & \textbf{59.8} & \textbf{58.5} & 5.1 & 3.6 & 4.3 & 3.6\\
2 & \textbf{66.0} & \textbf{60.1} & \textbf{59.7} & \textbf{66.5} & \textbf{59.9} & \textbf{51.9} & \textbf{51.5} & \textbf{54.7} & \textbf{1.6} & \textbf{1.1} & \textbf{1.3} & \textbf{2.1} &
\textbf{74.0} & \textbf{68.3} & \textbf{65.3} & \textbf{71.1} & \textbf{62.3} & \textbf{60.3} & \textbf{59.9} & \textbf{61.8} & 3.9 & \textbf{2.4} & \textbf{2.0} & \textbf{3.2}
& \textbf{76.2} & \textbf{71.7} & \textbf{66.7} & \textbf{73.4} & \textbf{68.8} & \textbf{62.0} & \textbf{62.3} & \textbf{62.5} & 5.1 & 3.6 & 4.3 & 3.6\\
3 & \textbf{59.0} & \textbf{51.6} & \textbf{50.6} & \textbf{54.2} & \textbf{56.5} & \textbf{44.9} & \textbf{46.7} & \textbf{49.5} & \textbf{1.4} & \textbf{0.7} & \textbf{1.8} & \textbf{1.4} &
\textbf{64.5} & \textbf{58.8} & \textbf{60.6} & \textbf{61.8} & \textbf{58.5} & \textbf{54.0} & \textbf{55.3} & \textbf{57.5} & 3.6 & \textbf{2.3} & \textbf{3.1} & 3.7
& \textbf{67.0} & \textbf{59.1} & \textbf{65.0} & \textbf{65.5} & \textbf{63.8} & \textbf{57.0} & \textbf{56.7} & \textbf{58.7} & 5.2 & 4.3 & \textbf{3.4} & 4.2 \\
6 & \textbf{47.6} & \textbf{43.9} & \textbf{39.9} & \textbf{45.2} & \textbf{49.8} & \textbf{37.9} & \textbf{41.0} & \textbf{39.6} & \textbf{0.5} &\textbf{ 0.2} & \textbf{0.4} & \textbf{0.5} &
\textbf{53.9} & \textbf{50.7} & \textbf{48.9} & \textbf{52.4} & \textbf{49.9} & \textbf{45.6} & \textbf{48.4} & \textbf{49.2} & \textbf{2.4} & \textbf{1.3} & \textbf{2.3} & \textbf{3.5}
& \textbf{56.5} & \textbf{53.9} & \textbf{49.4} & \textbf{54.3} & \textbf{55.6} & \textbf{49.5} & \textbf{49.2} & \textbf{51.3} & 4.1 & 3.8 & 3.9 & 4.9 \\
8 & \textbf{42.6} & \textbf{38.3} & \textbf{38.0} & \textbf{41.1} & \textbf{44.9} & \textbf{35.7} & \textbf{37.0} & \textbf{37.3} & \textbf{0.2} & \textbf{0.2} & \textbf{0.3} & \textbf{0.2} &
\textbf{51.8} & \textbf{41.6} & \textbf{42.7} & \textbf{45.1} & \textbf{46.9} & \textbf{42.1} & \textbf{43.9} & \textbf{45.0} & \textbf{2.1} & \textbf{1.7} & \textbf{1.6} & \textbf{1.9}
& \textbf{51.2} & \textbf{46.7} & \textbf{48.2} & \textbf{50.7} & \textbf{50.6} & \textbf{45.7} & \textbf{44.2} & \textbf{47.9} & 4.9 & 4.6 & 4.0 & 4.4 \\
10 & \textbf{36.9} & \textbf{33.4} & \textbf{36.2} & \textbf{37.8} & \textbf{42.1} & \textbf{35.0} & \textbf{35.5} & \textbf{36.7} & \textbf{0.2} & \textbf{0.0} & \textbf{0.0} & \textbf{0.1} &
\textbf{43.3} & \textbf{41.4} & \textbf{42.0} & \textbf{44.9} & \textbf{44.6} & \textbf{38.9} & \textbf{43.0} & \textbf{42.4} & \textbf{1.8 }& \textbf{1.4} & \textbf{1.7} & \textbf{1.0}
& \textbf{27.2} & \textbf{18.8} & \textbf{24.0} & \textbf{36.8} & \textbf{48.5} & \textbf{43.5} & \textbf{42.2} & \textbf{44.4} & 3.7 & \textbf{2.4} & 3.6 & 5.5 \\
\hline
& \multicolumn{36}{c|}{$\alpha=0.10$} \\
& \multicolumn{12}{c}{$N=200$} & \multicolumn{12}{c}{$N=1000$} & \multicolumn{12}{c|}{$N=5000$}\\
\cline{2-37}
&\multicolumn{4}{c}{$Q^1_M(\nu)$} &  \multicolumn{4}{c|}{$Q^2_M(\nu)$} &  \multicolumn{4}{c|}{$Q^3_M(\nu)$} & \multicolumn{4}{c}{$Q^1_M(\nu)$} &  \multicolumn{4}{c|}{$Q^2_M(\nu)$} & \multicolumn{4}{c|}{$Q^3_M(\nu)$} & \multicolumn{4}{c}{$Q^1_M(\nu)$} &  \multicolumn{4}{c|}{$Q^2_M(\nu)$} & \multicolumn{4}{c|}{$Q^3_M(\nu)$}\\
\hline
\diagbox{$M$}{$\nu$} & 1 & 2 & 3 & 4 & 1 & 2 & 3 & 4 & 1 & 2 & 3 & 4 & 1 & 2 & 3 & 4 & 1 & 2 & 3 & 4 & 1 & 2 & 3 & 4 & 1 & 2 & 3 & 4 & 1 & 2 & 3 & 4 & 1 & 2 & 3 & 4\\
\hline
1 & n.a. & n.a. & n.a. & n.a. & \textbf{66.7} & \textbf{57.3} & \textbf{60.2} & \textbf{63.1} & \textbf{7.5} & \textbf{7.9} & \textbf{7.2} & \textbf{7.4} & n.a. & n.a. & n.a. & n.a. & \textbf{66.0} & \textbf{64.1} & \textbf{67.8} & \textbf{68.0} & 9.7 & 9.2 & \textbf{7.8} & 10.9
& n.a. & n.a. & n.a. & n.a & \textbf{72.0} & \textbf{68.0} & \textbf{67.8} & \textbf{66.7} & 8.3 & 8.7 & 8.5 & 9.4\\
2 & \textbf{74.6} & \textbf{70.6} & \textbf{66.3} & \textbf{75.5} & \textbf{68.6} & \textbf{61.2} & \textbf{59.4} & \textbf{64.1} & \textbf{6.4} & \textbf{4.2} & \textbf{4.2} & \textbf{6.7} & \textbf{80.0} & \textbf{77.2} & \textbf{72.1} & \textbf{78.0} & \textbf{69.8} & \textbf{68.4} & \textbf{69.2} & \textbf{70.6} & 9.6 & \textbf{7.9} & \textbf{6.4} & 8.4
& \textbf{83.3} & \textbf{80.6} & \textbf{74.9} & \textbf{80.2} & \textbf{75.5} & \textbf{71.2} & \textbf{69.7} & \textbf{69.9} & 8.6 & 9.0 & 9.4 & \textbf{8.0} \\
3 & \textbf{67.9} & \textbf{62.0} & \textbf{59.9} & \textbf{63.1} & \textbf{64.4} & \textbf{54.6} & \textbf{55.1} & \textbf{58.2} & \textbf{4.7} & \textbf{2.9} & \textbf{4.9} & \textbf{3.7} & \textbf{74.8} & \textbf{68.0} & \textbf{68.6} & \textbf{69.2} & \textbf{65.7} & \textbf{62.5} & \textbf{64.8} & \textbf{67.5} & 9.7 & \textbf{7.1} & 8.6 & 8.8
& \textbf{75.5} & \textbf{68.0} & \textbf{73.6} & \textbf{74.0} & \textbf{72.0} & \textbf{66.1} & \textbf{65.7} & \textbf{66.7} & 9.3 & 8.9 & 9.1 & 10.2 \\
6 & \textbf{55.5} & \textbf{52.0} & \textbf{50.4} & \textbf{56.0} & \textbf{57.7} & \textbf{48.0} & \textbf{50.9} & \textbf{47.6} & \textbf{2.4} & \textbf{1.5 }& \textbf{2.2} & \textbf{2.4} &
\textbf{61.6} & \textbf{58.9} & \textbf{57.4} & \textbf{60.9} & \textbf{59.2} & \textbf{54.0} & \textbf{57.1} & \textbf{57.4} & \textbf{7.2} & \textbf{6.7} & \textbf{5.9} & \textbf{7.9}
& \textbf{64.1} & \textbf{62.9} & \textbf{59.2} & \textbf{63.8} & \textbf{65.0} & \textbf{57.6} & \textbf{58.0} & \textbf{61.2} & 8.5 & 9.7 & 9.4 & 10.2 \\
8 & \textbf{50.5} & \textbf{46.7} & \textbf{45.7} & \textbf{49.3} & \textbf{54.7 }& \textbf{45.3} & \textbf{46.9} & \textbf{45.5} & \textbf{1.0} & \textbf{1.1} & \textbf{0.8} & \textbf{0.7} & \textbf{61.0} &\textbf{ 50.9} & \textbf{51.8} & \textbf{53.8} & \textbf{54.5} & \textbf{51.4} & \textbf{53.7} & \textbf{53.7} & \textbf{7.0} & \textbf{4.5} & \textbf{6.4} & \textbf{5.3}
& \textbf{60.6} & \textbf{57.5} & \textbf{56.4} & \textbf{60.1} & \textbf{62.0} & \textbf{55.1} & \textbf{54.3} & \textbf{56.8} & \textbf{7.7} & 8.6 & 8.2 & 9.1 \\
10 & \textbf{47.2} & \textbf{41.8} & \textbf{44.7} & \textbf{46.9} & \textbf{51.8} & \textbf{42.5} & \textbf{44.7} & \textbf{44.0} & \textbf{0.6} & \textbf{0.3} & \textbf{0.4} & \textbf{0.6} & \textbf{53.6} & \textbf{50.0} & \textbf{49.7} & \textbf{54.1} & \textbf{52.9} & \textbf{47.8} & \textbf{51.8} & \textbf{52.0} & \textbf{5.1} & \textbf{5.1} & \textbf{5.1} & \textbf{5.7}
& \textbf{37.3} & \textbf{29.4} & \textbf{33.7} & \textbf{45.9} & \textbf{58.6} & \textbf{52.6 }& \textbf{52.5} & \textbf{54.7} & 8.3 & \textbf{6.0} & 8.1 & 10.2 \\ \hline
\end{tabular}
\end{table}

\begin{table}[H]
\caption{\scriptsize{Empirical levels at the $5\%$ and $10\%$ significance levels for the global portmanteau test statistics $Q_M$ for $M = 1,2,3,6,8,10$. The simulated model is a weak PVAR given by DGP$_2$. The numbers of years are equal to $N = 200$, $1000$ and $5000$. $\mathrm{{Q}}_M^{1}(\nu)$ is the global portmanteau test statistics defined by~(37) in~\cite{UD09}. $\mathrm{{Q}}_M^{2}(\nu)$ is the global portmanteau test statistics defined by~(28) in~\citet{DL13}. $\mathrm{{Q}}_M^{3}(\nu)$ is the global portmanteau test statistics defined by~(\ref{QMast}).}}\label{global_test_portmanteau_weak_const}
\begin{center}
\tiny
\setlength{\tabcolsep}{0.8mm}
\begin{tabular}{|c|ccc|ccc|ccc|ccc|ccc|ccc|} \hline
    & \multicolumn{9}{c|}{$\alpha=0.05$} & \multicolumn{9}{c|}{$\alpha=0.10$} \\
    & \multicolumn{3}{c}{$N=200$} & \multicolumn{3}{c}{$N=1000$} & \multicolumn{3}{c|}{$N=5000$}
    & \multicolumn{3}{c}{$N=200$} & \multicolumn{3}{c}{$N=1000$} & \multicolumn{3}{c|}{$N=5000$}\\ \hline
$M$ & $Q^1_M$ & $Q^2_M$ & $Q^3_M$ & $Q^1_M$ & $Q^2_M$ & $Q^3_M$ & $Q^1_M$ & $Q^2_M$ & $Q^3_M$ & $Q^1_M$ & $Q^2_M$ & $Q^3_M$ & $Q^1_M$ & $Q^2_M$ & $Q^3_M$ & $Q^1_M$ & $Q^2_M$ & $Q^3_M$\\ \hline
1 & n.a. & \textbf{90.8} & \textbf{1.3} & n.a. & \textbf{94.5} & \textbf{3.4} & n.a. & \textbf{95.7} & 4.5 & n.a. & \textbf{93.4} & \textbf{5.5} & n.a. & \textbf{96.9} & \textbf{7.7} & n.a. & \textbf{97.3} & \textbf{7.9} \\
2 & \textbf{96.2} & \textbf{91.0} & \textbf{1.2} & \textbf{98.9} & \textbf{95.9} & \textbf{2.9} & \textbf{99.3} & \textbf{97.7} & \textbf{3.4} & \textbf{98.5} & \textbf{94.6} & \textbf{3.4} & \textbf{99.5} & \textbf{97.1}  & \textbf{6.2} & \textbf{99.6} & \textbf{98.6} & 8.9 \\
3 & \textbf{91.3} & \textbf{86.6} & \textbf{0.3} & \textbf{97.9} & \textbf{92.8} & \textbf{2.0} & \textbf{98.3} & \textbf{96.0} & \textbf{2.7} & \textbf{93.7} & \textbf{91.8} & \textbf{2.2} & \textbf{98.7} & \textbf{95.3} & \textbf{6.2} & \textbf{99.0} & \textbf{97.4} & \textbf{8.0} \\
6 & \textbf{81.0} & \textbf{75.9} & \textbf{0.0} & \textbf{88.8} & \textbf{85.5} & \textbf{1.8} & \textbf{91.8} & \textbf{90.8} & \textbf{2.9} & \textbf{85.9} & \textbf{82.7} & \textbf{1.1} & \textbf{92.1}  & \textbf{90.4} & \textbf{5.1} & \textbf{94.7} & \textbf{93.8} & \textbf{7.5} \\
8 & \textbf{73.6} & \textbf{72.0} & \textbf{0.0} & \textbf{82.2} & \textbf{81.2} & \textbf{1.8} & \textbf{87.6} & \textbf{86.6} & 4.1 & \textbf{80.1} & \textbf{79.7} & \textbf{0.2 }& \textbf{87.6} & \textbf{86.6} & \textbf{4.8} & \textbf{92.5} & \textbf{90.3} & 8.1 \\
10 & \textbf{65.9} & \textbf{68.1} & \textbf{0.0} & \textbf{80.3} & \textbf{78.3} & \textbf{1.2} & \textbf{55.1} & \textbf{83.8} & 3.8 & \textbf{73.2} & \textbf{75.3} & \textbf{0.0} & \textbf{86.1} & \textbf{84.3} & \textbf{3.8} & \textbf{66.0} & \textbf{89.1} & 8.9 \\ \hline
\end{tabular}
\end{center}
\end{table}

As expected, for the weak PVAR models (unconstrained and constrained), the relative
rejection frequencies of the tests defined by~\citet{UD09} and by~\citet{DL13} are definitely outside the significant limits. Based on these two aforementioned tests, an overly complicated PVAR model will be used leading to misinterpretation and loss of efficiency in terms of linear predictions. It should be noted that, when the sample size is small ($N=200$) and $M\geq 2$, the relative rejection frequencies of our proposed test are far from the nominal $5\%$. This is not surprising that results improved when $N$ becomes large enough.

\subsection{Empirical power}
\noindent In this part, we simulated $1000$ independent trajectories of size $N=1000$ and $5000$ of a weak PVAR(2) defined by
\begin{equation}\label{pvar2_model}
\bfY_{ns+\nu} = \bfPhi_1(\nu)\bfY_{ns+\nu-1} + \bfPhi_2(\nu)\bfY_{ns+\nu-2} + \bfepsilon_{ns+\nu}, \quad \bfepsilon_{ns+\nu} =
\begin{pmatrix}
\eta_{1,ns+\nu}\eta_{1,ns+\nu-1}\eta_{1,ns+\nu-2} \\
\eta_{2,ns+\nu}\eta_{2,ns+\nu-1}\eta_{2,ns+\nu-2}
\end{pmatrix},
\end{equation}
where $\bfeta_{t} =
\begin{pmatrix}
\eta_{1,t}\\
\eta_{2,t}
\end{pmatrix}
$ iid $\mathcal{N}(\bfzero,\bfI_2)$. We considered the case of four seasons, that is $\nu=4$. The autoregressive coefficients of the model are given by:
\begin{eqnarray*}
\bfPhi_1(1) &=& \left( \begin{array}{cc}
                   0.50 & 0.30 \\
                   0.10 & 0.20
                   \end{array}
            \right), \;
\bfPhi_1(2) = \left( \begin{array}{cc}
                   0.42 & 0.24 \\
                   -0.20 & 0.50
                   \end{array}
            \right), \;\\
\bfPhi_1(3) &=& \left( \begin{array}{cc}
                   -0.80 & 0.20 \\
                   0.60 & 0.70
                   \end{array}
            \right), \;
\bfPhi_1(4) = \left( \begin{array}{cc}
                   -0.30 & 0.50 \\
                   0.90 & -0.20
                   \end{array}
            \right), \;\\
\bfPhi_2(1) &=& \left( \begin{array}{cc}
                   0.30 & 0.00 \\
                   0.00 & 0.20
                   \end{array}
            \right), \;
\bfPhi_2(2) = \left( \begin{array}{cc}
                   0.20 & 0.00 \\
                   0.00 & -0.30
                   \end{array}
            \right), \;\\
\bfPhi_2(3) &=& \left( \begin{array}{cc}
                   -0.30 & 0.00 \\
                   0.00 & 0.20
                   \end{array}
            \right), \;
\bfPhi_2(4) = \left( \begin{array}{cc}
                   -0.30 & 0.00 \\
                   0.00 & -0.20
                   \end{array}
            \right). \;
\end{eqnarray*}

For each of these $1000$ replications we fitted a PVAR(1) model and perform the following tests: $\mathrm{{Q}}_M^{1}(\nu)$ portmanteau test statistics defined by~(35) in~\citet{UD09}, $\mathrm{{Q}}_M^{2}(\nu)$ portmanteau test defined by~(26) in~\citet{DL13} and the $\mathrm{{Q}}_M^{3}(\nu)$ portmanteau test defined by~(\ref{QMnuast}) based on $M=1,2,3,6,8,10$ residual autocorrelations. The adequacy of the PVAR(1) model is rejected when the p-value is less than $5\%$. Table~\ref{power} displays the relative frequency of rejection over the $1000$ replications. As emphasized by~\citet{FR07}, one could think that our proposed portmanteau test is less powerful than the other two tests. For weak PVAR models, we have seen that the actual level of the tests proposed by~\citet{UD09} or by~\citet{DL13} is generally much greater than the $5\%$ nominal level (see Table~\ref{test_portmanteau_weak_unc}).

\begin{table}[H]
\caption{\scriptsize{Empirical power (in $\%$) of the $\mathrm{{Q}}_M^{1}(\nu)$ portmanteau test statistics defined by~(35) in~\citet{UD09}, the $\mathrm{{Q}}_M^{2}(\nu)$ portmanteau test defined by~(26) in~\citet{DL13} and the $\mathrm{{Q}}_M^{3}(\nu)$ portmanteau test defined by~(\ref{QMnuast}) at the $5\%$ significance level. The simulated model is a weak PVAR given by eq.~(\ref{pvar2_model}). The numbers of years are equal to $N = 1000$ and $5000$}.
}\label{power}
\centering
\tiny
\setlength{\tabcolsep}{0.8mm}
\begin{tabular}{|c|cccc|cccc|cccc|cccc|cccc|cccc|cccc|}
\hline
& \multicolumn{24}{c|}{$\alpha=0.05$} \\
& \multicolumn{12}{c}{$N=1000$} & \multicolumn{12}{c|}{$N=5000$}\\
\cline{2-25}
& \multicolumn{4}{c}{$Q^1_M(\nu)$} &  \multicolumn{4}{c|}{$Q^2_M(\nu)$} & \multicolumn{4}{c|}{$Q^3_M(\nu)$} & \multicolumn{4}{c}{$Q^1_M(\nu)$} &  \multicolumn{4}{c|}{$Q^2_M(\nu)$} & \multicolumn{4}{c|}{$Q^3_M(\nu)$}\\
\hline
\diagbox{$M$}{$\nu$}& 1 & 2 & 3 & 4 & 1 & 2 & 3 & 4 & 1 & 2 & 3 & 4 & 1 & 2 & 3 & 4 & 1 & 2 & 3 & 4 & 1 & 2 & 3 & 4\\
\hline
1 & n.a.  & n.a.  & n.a.  & n.a.  & 100.0 & 100.0 & 100.0 & 100.0 & 67.6 & 100.0 & 100.0 & 99.6 &
n.a. & n.a. & n.a. & n.a. & 100.0 & 100.0 & 100.0 & 100.0 & 99.9 & 100.0 & 100.0 & 100.0 \\
2 & 100.0 & 100.0 & 100.0 & 100.0 & 100.0 & 100.0 & 100.0 & 100.0 & 74.8 & 100.0 & 100.0 & 100.0
& 100.0 & 100.0 & 100.0 & 100.0 & 100.0 & 100.0 & 100.0 & 100.0 & 100.0 & 100.0 & 100.0 & 100.0 \\
3 & 100.0 & 100.0 & 100.0 & 100.0 & 100.0 & 100.0 & 100.0 & 100.0 & 79.5 & 99.7 & 100.0 & 99.6
& 100.0 & 100.0 & 100.0 & 100.0 & 100.0 & 100.0 & 100.0 & 100.0 & 100.0 & 100.0 & 100.0 & 100.0 \\
6 & 99.9 & 100.0 & 100.0 & 100.0 & 99.9 & 100.0 & 100.0 & 100.0 & 79.0 & 100.0 & 100.0 & 99.4 &
100.0 & 100.0 & 100.0 & 100.0 & 100.0 & 100.0 & 100.0 & 100.0 & 100.0 & 100.0 & 100.0 & 100.0 \\
8 & 99.7 & 100.0 & 100.0 & 100.0 & 99.7 & 100.0 & 100.0 & 100.0 & 76.7 & 100.0 & 99.9 & 99.4
& 100.0 & 100.0 & 100.0 & 100.0 & 100.0 & 100.0 & 100.0 & 100.0 & 100.0 & 100.0 & 100.0 & 100.0 \\
10 & 99.6 & 100.0 & 100.0 & 100.0 & 99.6 & 100.0 & 100.0 & 100.0 & 76.3 & 100.0 & 100.0 & 98.9
& 100.0 & 100.0 & 100.0 & 100.0 & 100.0 & 100.0 & 100.0 & 100.0 & 100.0 & 100.0 & 100.0 & 100.0 \\
\hline
\end{tabular}
\end{table}

\section{Application to real data}\label{real}
\subsection{The CAC 40 and DAX indices as an illustrative example}
\noindent We illustrate here the new portmanteau test statistics with a real data set derived from finance. The variables are the daily returns of two European stock market indices: CAC 40 (Paris) and DAX (Frankfurt), from March $3$, $1990$ to March $10$, $2022$. The data were obtained from \textit{Yahoo Finance}. Because of the legal holidays, many weeks comprise less than five observations. We preferred removing the entire weeks when there was less than five data available, giving a bivariate time series of sample size equal to $7060$. The period $\nu=5$ is naturally selected. These data were analyzed by~\citet{BMU23}.

\noindent A PVAR model of order $1$ is fitted to the bivariate series of observations:
$$\bfY_{ns+\nu} = \bfPhi(\nu)\bfY_{ns+\nu-1} +\bfepsilon_{ns+\nu}\quad \nu = 1,\ldots,5,$$
where
$\bfY_t = \left(r_t^1,r_t^2\right)^\top$ and $r_t^1$, $r_t^2$ represent the log-return of CAC 40 and DAX respectively.
The log-return is defined as $r_t = 100\times\ln{(I_t/I_{t-1})}$ where $I_t$ represents the value of the index at time $t$. Seasonal means are first removed from the series, meaning that a model is formulated by examining $\bfY_{ns+\nu}-\bfmu(\nu)$.

\begin{table}[H]
\scriptsize
\caption{Least squares estimators used to fit the log-returns of CAC 40 and DAX data to a bivariate PVAR model with $\nu=5$; the $\hat{\sigma}_\mathrm{S}$ and $\hat{\sigma}_\mathrm{{SP}}$ represent the standard errors in the strong case and for our proposed estimators in the weak case; pval$_S$ and pval$_\mathrm{{SP}}$ correspond to the $p$-values of the $t$-statistic of  $\hat{\bfbeta}$.}\label{CAC40_DAX}
\centering
\begin{tabular}{rrrrrr}
  \hline
    & $\hat{\beta}$ & $\hat{\sigma}_\mathrm{S}$ & $\hat{\sigma}_\mathrm{{SP}}$ & $\text{pval}_\mathrm{S}$ & $\text{pval}_\mathrm{{SP}}$ \\ 
  \hline
  1 & -0.0349 & 0.0707 & 0.1456 & 0.6220 & 0.8107 \\ 
  2 & 0.0153 & 0.0731 & 0.1480 & 0.8346 & 0.9180 \\ 
  3 & -0.0070 & 0.0706 & 0.0992 & 0.9215 & 0.9441 \\ 
  4 & -0.0378 & 0.0729 & 0.1019 & 0.6044 & 0.7106 \\ 
  5 & -0.0506 & 0.0399 & 0.0524 & 0.2045 & 0.3339 \\ 
  6 & -0.0270 & 0.0420 & 0.0663 & 0.5214 & 0.6843 \\ 
  7 & -0.0020 & 0.0386 & 0.0551 & 0.9591 & 0.9713 \\ 
  8 & -0.0246 & 0.0407 & 0.0742 & 0.5450 & 0.7399 \\ 
  9 & -0.3001 & 0.0532 & 0.1296 & 0.0000 & 0.0207 \\ 
  10 & -0.1256 & 0.0545 & 0.0781 & 0.0214 & 0.1078 \\ 
  11 & 0.2605 & 0.0505 & 0.1381 & 0.0000 & 0.0595 \\ 
  12 & 0.0360 & 0.0517 & 0.0607 & 0.4869 & 0.5535 \\ 
  13 & -0.1498 & 0.0548 & 0.1020 & 0.0063 & 0.1422 \\ 
  14 & -0.0744 & 0.0551 & 0.0639 & 0.1767 & 0.2445 \\ 
  15 & 0.1862 & 0.0538 & 0.1021 & 0.0006 & 0.0683 \\ 
  16 & 0.0955 & 0.0541 & 0.0620 & 0.0778 & 0.1240 \\ 
  17 & -0.0227 & 0.0521 & 0.0681 & 0.6627 & 0.7385 \\ 
  18 & -0.0055 & 0.0522 & 0.0670 & 0.9156 & 0.9343 \\ 
  19 & 0.0694 & 0.0520 & 0.0739 & 0.1823 & 0.3480 \\ 
  20 & 0.0420 & 0.0521 & 0.0734 & 0.4202 & 0.5677 \\ 
   \hline
\end{tabular}
\end{table}

\noindent We present in Table~\ref{CAC40_DAX} the estimated parameters $\hat{\bfbeta}=(\hat{\bfbeta}(1),\dots,\hat{\bfbeta}(5))^\top$ and their estimated standard error proposed in the strong case ($\hat{\sigma}_\mathrm{S}$) and the weakly consistent estimators proposed ($\hat{\sigma}_\mathrm{{SP}}$). In the weak case, none of them are significant at the $2\%$ level. This is in accordance with the results of~\cite{frs11} who showed that the log-returns of these two European stock market indices constitute a weak periodic white noises.


First, we apply portmanteau tests $Q_M^1(\nu)$ and $Q^2_M(\nu)$ for checking the hypothesis that the CAC 40 and DAX log-returns constitute a white noise. The $P$-values are reported in Table~\ref{pvalueQM}. Since the $P$-values of these two tests are very small, the white-noise hypothesis is rejected at the nominal level $\alpha = 2\%$. This is not surprising because these tests require the iid assumption and it is well known that the strong white-noise model is not adequate for these series~\citep{BMS18}. In contrast, the white-noise hypothesis is not rejected by our test $Q^3_M(\nu)$ since the statistic is not larger than the critical values. To summarize, the outputs of Table~\ref{pvalueQM} are in accordance with the common belief that these series are not strong white noises, but could be weak white noises.

\begin{table}[H]
\begin{center}
\scriptsize
\caption{$P$-values of the portmanteau test statistics used to check
a bivariate PVAR model with
$\nu=5$.}\label{pvalueQM}
\label{pvaluePVAR}
\setlength{\tabcolsep}{1.25mm}
\begin{tabular}{|c|ccccc|ccccc|ccccc|}
\hline
& \multicolumn{5}{c}{$Q^1_M(\nu)$} & \multicolumn{5}{c}{$Q^2_M(\nu)$} & \multicolumn{5}{c|}{$Q^3_M(\nu)$} \\
\hline
\diagbox{$M$}{$\nu$} & 1 & 2 & 3 & 4 & 5 & 1 & 2 & 3 & 4 & 5 & 1 & 2 & 3 & 4 & 5 \\
\hline
  1 & 0.51 & 0.06 & \textbf{0.00} & \textbf{0.00} & 0.34 & 0.51 & 0.06 & \textbf{0.00} & \textbf{0.00} & 0.34 & 0.82 & 0.41 & 0.10 & 0.45 & 0.69 \\ 
  2 & 0.04 & \textbf{0.00} & \textbf{0.00} & \textbf{0.00} & 0.02 & 0.04 & \textbf{0.00} & \textbf{0.00} & \textbf{0.00} & 0.02 & 0.13 & 0.02 & 0.09 & 0.32 & 0.07 \\ 
  3 & 0.07 & \textbf{0.00} & \textbf{0.00} & \textbf{0.00} & 0.03 & 0.07 & \textbf{0.00} & \textbf{0.00} & \textbf{0.00} & 0.03 & 0.67 & 0.09 & 0.14 & 0.43 & 0.47 \\ 
  6 & 0.14 & \textbf{0.00} & \textbf{0.00} & \textbf{0.00} & \textbf{0.00} & 0.14 & \textbf{0.00} & \textbf{0.00} & \textbf{0.00} & 0.00 & 0.79 & 0.11 & 0.12 & 0.46 & 0.45 \\ 
  8 & 0.13 & \textbf{0.00} & \textbf{0.00} & \textbf{0.00} & \textbf{0.00} & 0.13 & \textbf{0.00} & \textbf{0.00} & \textbf{0.00} & \textbf{0.00} & 0.81 & 0.09 & 0.17 & 0.50 & 0.44 \\ 
 10 & 0.05 & \textbf{0.00} & \textbf{0.00} & \textbf{0.00} & \textbf{0.00} & 0.05 & \textbf{0.00} & \textbf{0.00} & \textbf{0.00} & \textbf{0.00} & 0.74 & 0.12 & 0.19 & 0.46 & 0.51 \\ 
\hline
\end{tabular}
\end{center}
\end{table}

\begin{figure}[h] 
  \begin{minipage}[b]{0.48\linewidth}
    \includegraphics[width=.95\linewidth]{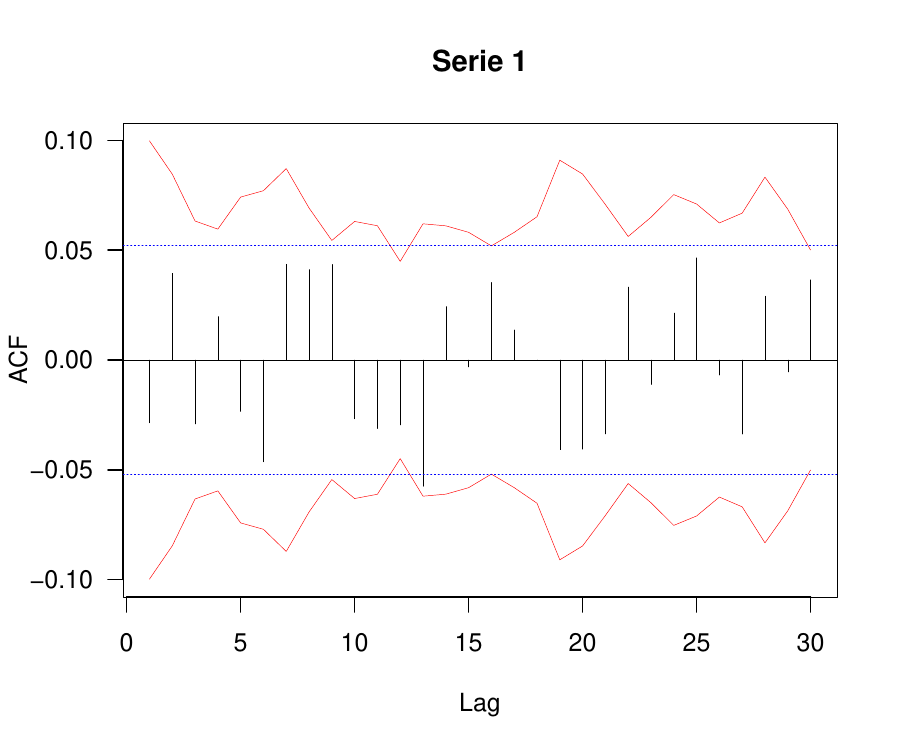} 
  \end{minipage} 
  \begin{minipage}[b]{0.48\linewidth}
    \includegraphics[width=.95\linewidth]{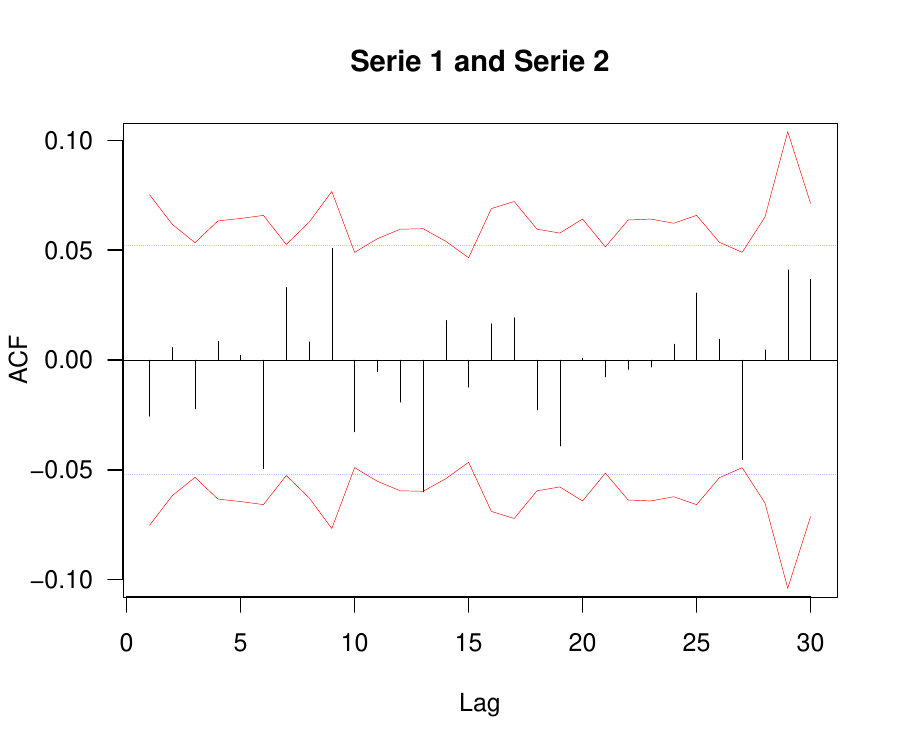} 
  \end{minipage} 
  \begin{minipage}[b]{0.48\linewidth}
    \includegraphics[width=.95\linewidth]{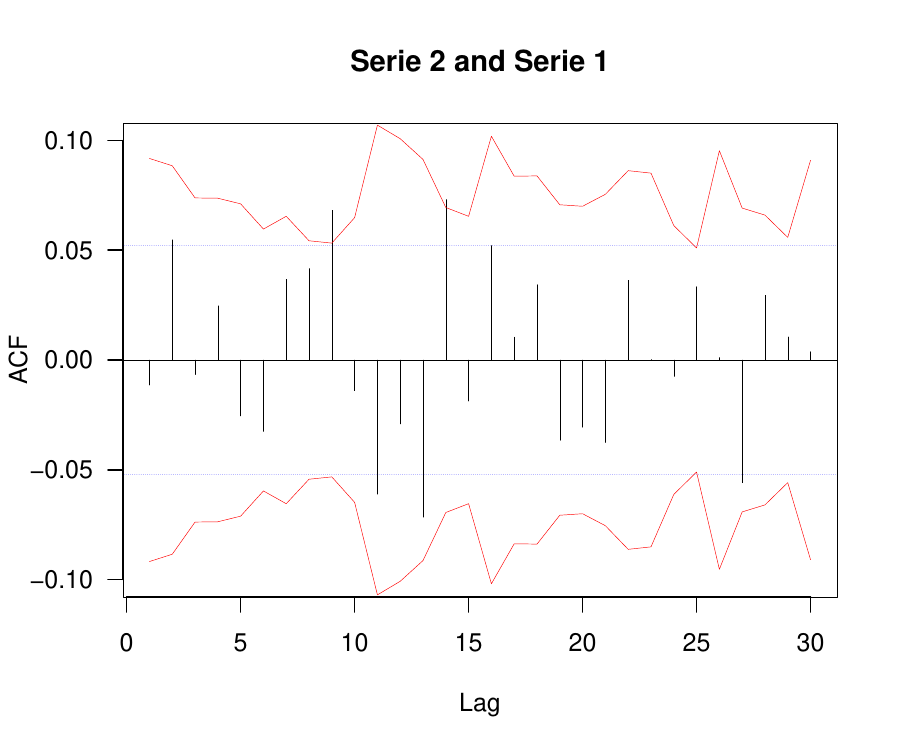} 
  \end{minipage}
  \hfill
  \begin{minipage}[b]{0.48\linewidth}
    \includegraphics[width=.95\linewidth]{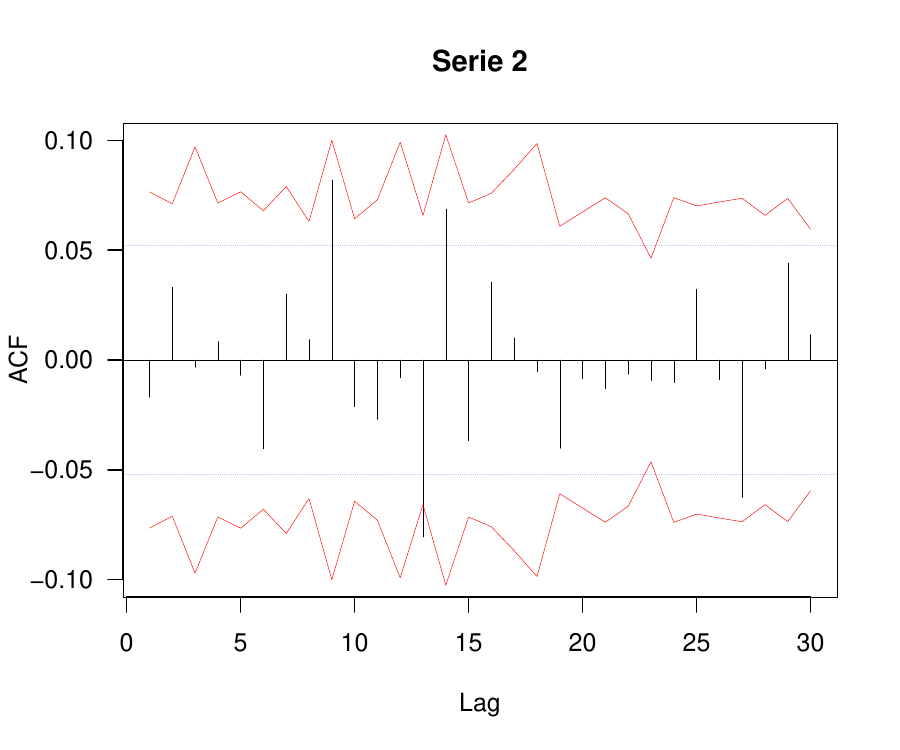} 
  \end{minipage} 
  \caption{Autocorrelations of the log-returns of CAC 40 and DAX for the first period (serie 1 corresponds to CAC 40).}\label{autocorrfig} 
\end{figure}

In view of Figure~\ref{autocorrfig}, the diagnostic checking of residuals does not indicate any inadequacy, approximately $95\%$ of the sample autocorrelations should lie between the bands shown as solid lines (red color).

\subsection{The monthly flows series as an illustrative example}\label{river}
\begin{figure}[ht]
  \centering
  \includegraphics[width=0.9\textwidth]{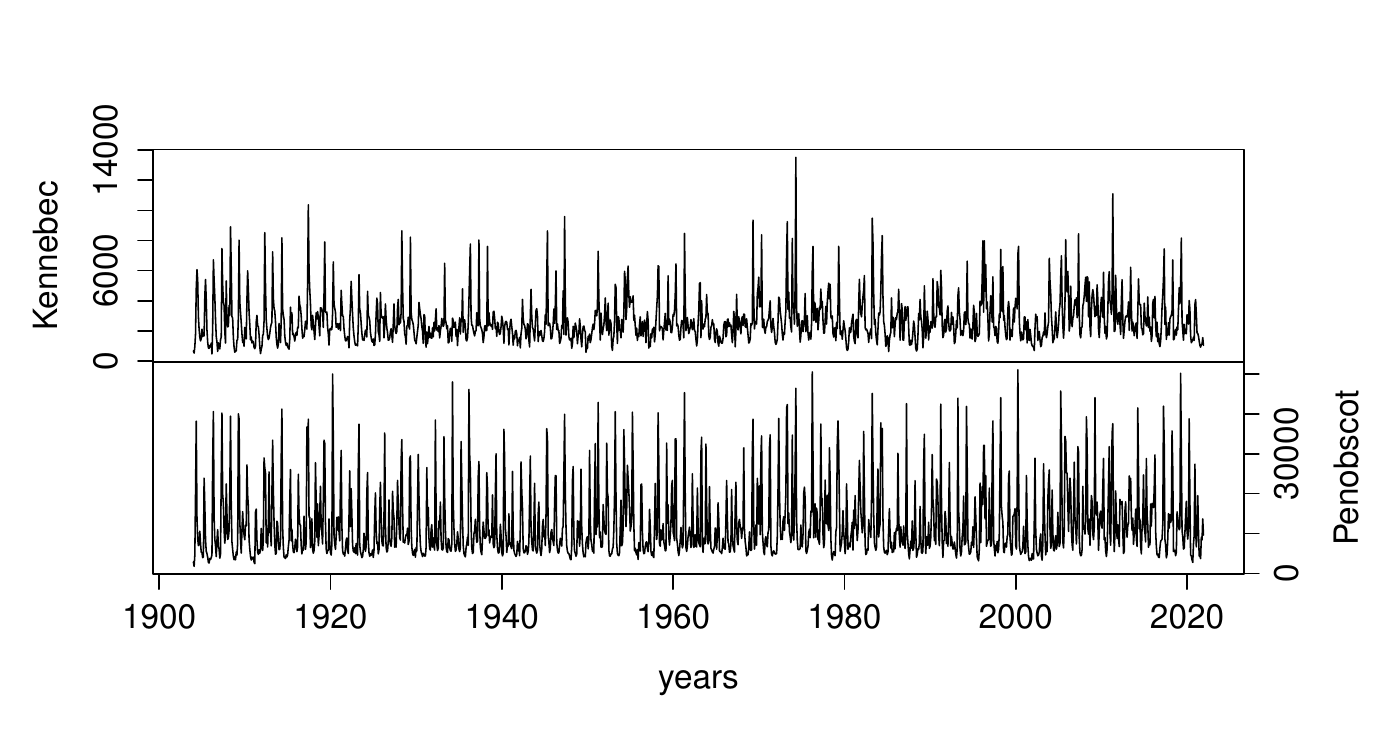}
  \caption{Monthly flows of Kennebec and Penobscot river}\label{river_flows}
\end{figure}

\noindent We illustrate here the effectiveness of the proposed methodology in river flow analysis. Data related to two rivers of the same length, with different means of annual flows, located in Maine, will be examined. They consist of:
\begin{itemize}
\item flows of Kennebec river, measured at The Forks, Maine;
\item flows of Penobscot river, measured at West Enfield, Maine.
\end{itemize}
Data are obtained from daily discharge measurements in cubic feet per second from January $1904$ to December $2022$ (USGS Surface-Water Monthly Statistics). Daily data flows are then transformed in monthly data, giving a bivariate time series of sample size equal to $1428$ ($119$ years). The period $\nu=12$ is naturally selected. Table~\ref{means_months} presents the means of Kennebec river and Piscataquis river, their p-values (in parentheses) and their estimated
standard errors (in bracquets). The monthly flows series from $1904$ to $2021$ are plotted in Figure~\ref{river_flows}. 

\noindent To ensure that residuals of the fitted model are approximately normally distributed and homoscedastic, a logarithmic transformation has been used~\citep{MG13}. Seasonal means are first removed from the series, meaning that a model is formulated by examining $\bfY_{ns+\nu}-\bfmu_{\nu}$. A PVAR model of order $1$ is fitted to the bivariate series of observations:
$$\bfY_{ns+\nu} = \bfPhi(\nu)\bfY_{ns+\nu-1} +\bfepsilon_{ns+\nu}\quad \nu = 1,\ldots,12,$$
where $\bfY_t = \left(y_t^1,y_t^2\right)^\top$ and $y_t^1$, $y_t^2$ represent the log-monthly flows of Kennebec river and Penobscot river respectively. The last year of the data set ($12$ observations) has been omitted for the estimation of the periodic model. 

\begin{table}[H]
\tiny
\caption{Periodic means, their p-values (in parentheses) and their estimated standard errors (in brackets)}\label{means_months}
\centering
\begin{tabular}{|l|c|c|c|c|c|c|c|c|c|c|c|c|}
  \hline
  month & 1 & 2 & 3 & 4 & 5 & 6 & 7 & 8 & 9 & 10 & 11 & 12 \\
  \hline
  Kennebec & 2502.950 & 2731.793 & 2480.923 & 3045.939 & 4672.663 & 3293.709 & 2628.069 & 2366.628 & 2280.790 & 2036.861 & 1965.089 & 2243.970 \\ 
           & (0.013) & (0.010) & (0.038) & (0.077) & (0.081) & (0.027) & (0.012) & (0.001) & (0.002) & (0.011) & (0.100) & (0.075) \\ 
           & [993.061] & [1039.111] & [1180.615] & [1705.674] & [2656.274] & [1468.866] & [1024.917] & [696.912] & [711.597] & [784.091] & [1183.652] & [1250.200] \\ 
           \hline
  Penobscot & 8597.168 & 7804.975 & 11225.471 & 30028.739 & 23391.815 & 11594.899 & 7772.168 & 6509.739 & 6724.05 & 9098.076 & 12496.613 & 11785.25 \\
          & (0.012) & (0.025) & (0.066) & (0.001) & (0.019) & (0.041) & (0.024) & (0.012) & (0.050) & (0.087) & (0.056) & (0.075) \\
          & [3380.514] & [3449.669] & [6057.943] & [8823.026] & [9830.225] & [5623.413] & [3402.226] & [2549.783] & [3388.915] & [5270.996] & [6465.082] & [6561.644] \\ 
  \hline
\end{tabular}
\end{table}

We present in Table~\ref{Kennebec_Penobscot} the estimated parameters $\hat{\bfbeta}=(\hat{\bfbeta}(1),\dots,\hat{\bfbeta}(12))^\top$ and their estimated standard error proposed in the strong case denoted $\hat{\sigma}_\mathrm{S}$ and the weakly consistent estimators proposed denoted by $\hat{\sigma}_\mathrm{{SP}}$. In the weak case, $46$ of them are significant at the $1\%$ level. This is in accordance with the results of~\cite{WGVM05} who showed that the weak PVAR model is adequate in modelling monthly flows. Surprisingly, in the strong case, no parameters should be estimated, meaning that the log-monthly flows of these two series constitute a strong periodic white noise.

\begin{table}[H]
\scriptsize
\caption{Least squares estimators used to fit the log of monthly flows of Kennebec and Piscataquis river data to a bivariate PVAR model with $\nu=12$; the $\hat{\sigma}_\mathrm{S}$ and $\hat{\sigma}_\mathrm{{SP}}$ represent the standard errors in the strong case and for our proposed estimators in the weak case; pval$_S$ and pval$_\mathrm{{SP}}$ correspond to the $p$-values of the $t$-statistic of  $\hat{\bfbeta}$. The p-values less than $1\%$ are in bold.}\label{Kennebec_Penobscot}
\centering
\begin{tabular}{rrrrrr}
  \hline
    & $\hat{\beta}$ & $\hat{\sigma}_\mathrm{S}$ & $\hat{\sigma}_\mathrm{{SP}}$ & $\text{pval}_\mathrm{S}$ & $\text{pval}_\mathrm{{SP}}$ \\ 
  \hline
  1 & 0.5363 & 0.6827 & 0.0020 & 0.4338 & \textbf{0.0000} \\ 
  2 & 0.1335 & 0.7934 & 0.0028 & 0.8667 & \textbf{0.0000} \\ 
  3 & 0.2764 & 0.6150 & 0.0020 & 0.6539 & \textbf{0.0000} \\ 
  4 & 0.3530 & 0.7147 & 0.0028 & 0.6223 & \textbf{0.0000} \\ 
  5 & 0.5958 & 0.6004 & 0.0016 & 0.3231 & \textbf{0.0000} \\ 
  6 & 0.1698 & 0.7467 & 0.0022 & 0.8206 & \textbf{0.0000} \\ 
  7 & 0.4258 & 0.6995 & 0.0016 & 0.5440 & \textbf{0.0000} \\ 
  8 & 0.6790 & 0.8701 & 0.0022 & 0.4368 & \textbf{0.0000} \\ 
  9 & 0.7391 & 0.8167 & 0.0015 & 0.3674 & \textbf{0.0000} \\ 
  10 & -0.0212 & 1.0866 & 0.0048 & 0.9844 & \textbf{0.0000} \\ 
  11 & 0.3461 & 0.8708 & 0.0015 & 0.6917 & \textbf{0.0000} \\ 
  12 & 0.6536 & 1.1584 & 0.0048 & 0.5737 & \textbf{0.0000} \\ 
  13 & 0.2841 & 0.8580 & 0.0083 & 0.7412 & \textbf{0.0000} \\ 
  14 & 0.1496 & 0.5828 & 0.0048 & 0.7978 & \textbf{0.0000} \\ 
  15 & 0.5654 & 0.9338 & 0.0083 & 0.5461 & \textbf{0.0000} \\ 
  16 & 0.1417 & 0.6342 & 0.0048 & 0.8236 & \textbf{0.0000} \\ 
  17 & 0.5685 & 1.5459 & 0.0224 & 0.7138 & \textbf{0.0000} \\ 
  18 & 0.0851 & 1.2271 & 0.0166 & 0.9448 & \textbf{0.0000} \\ 
  19 & 0.2967 & 2.7414 & 0.0224 & 0.9140 & \textbf{0.0000} \\ 
  20 & 0.2136 & 2.1759 & 0.0166 & 0.9220 & \textbf{0.0000} \\ 
  21 & 0.2868 & 0.8727 & 0.0078 & 0.7430 & \textbf{0.0000} \\ 
  22 & 0.1394 & 0.8647 & 0.0077 & 0.8722 & \textbf{0.0000 }\\ 
  23 & 0.0884 & 1.3000 & 0.0078 & 0.9459 & \textbf{0.0000} \\ 
  24 & 0.2476 & 1.2881 & 0.0077 & 0.8479 & \textbf{0.0000} \\ 
  25 & 0.4377 & 1.1100 & 0.0037 & 0.6941 & \textbf{0.0000} \\ 
  26 & -0.0066 & 1.0916 & 0.0032 & 0.9952 & 0.0425 \\ 
  27 & 0.0676 & 1.1568 & 0.0037 & 0.9535 & \textbf{0.0000} \\ 
  28 & 0.5137 & 1.1376 & 0.0032 & 0.6524 & \textbf{0.0000} \\ 
  29 & 0.1828 & 0.8642 & 0.0024 & 0.8329 & \textbf{0.0000} \\ 
  30 & -0.2056 & 0.9667 & 0.0024 & 0.8320 & \textbf{0.0000} \\ 
  31 & 0.3248 & 0.8756 & 0.0024 & 0.7113 & \textbf{0.0000} \\ 
  32 & 0.6588 & 0.9795 & 0.0024 & 0.5026 & \textbf{0.0000} \\ 
  33 & 0.1557 & 0.9342 & 0.0022 & 0.8679 & \textbf{0.0000} \\ 
  34 & -0.3257 & 1.0695 & 0.0019 & 0.7613 & \textbf{0.0000} \\ 
  35 & 0.3184 & 0.8235 & 0.0022 & 0.6997 & \textbf{0.0000} \\ 
  36 & 0.8045 & 0.9427 & 0.0019 & 0.3953 & \textbf{0.0000} \\ 
  37 & 0.7662 & 0.9462 & 0.0022 & 0.4198 & \textbf{0.0000} \\ 
  38 & -0.1897 & 1.5668 & 0.0056 & 0.9038 & \textbf{0.0000} \\ 
  39 & 0.1641 & 0.7149 & 0.0022 & 0.8189 & \textbf{0.0000} \\ 
  40 & 0.9073 & 1.1838 & 0.0056 & 0.4450 & \textbf{0.0000} \\ 
  41 & 0.4522 & 1.1802 & 0.0053 & 0.7023 & \textbf{0.0000} \\ 
  42 & 0.0015 & 1.2551 & 0.0059 & 0.9991 & 0.8027 \\ 
  43 & 0.5145 & 0.8184 & 0.0053 & 0.5308 & \textbf{0.0000} \\ 
  44 & 0.7333 & 0.8704 & 0.0059 & 0.4013 & \textbf{0.0000} \\ 
  45 & 0.5088 & 0.7600 & 0.0025 & 0.5045 & \textbf{0.0000} \\ 
  46 & 0.1563 & 1.0850 & 0.0073 & 0.8857 & \textbf{0.0000} \\ 
  47 & 0.2903 & 0.7218 & 0.0025 & 0.6883 & \textbf{0.0000} \\ 
  48 & 0.4929 & 1.0305 & 0.0073 & 0.6334 & \textbf{0.0000} \\ 
   \hline
\end{tabular}
\end{table}

\noindent Then, in a second step, the reduced PVAR model was estimated with constraints on the autoregressive parameters with non-significant $\text{pval}_\mathrm{{SP}}$ in Table~\ref{Kennebec_Penobscot}. The least squares estimators of the final model are presented 
in Table~\ref{LSKennebec_Penobscot}.
\begin{table}[H]
\scriptsize
\begin{center}
\caption{\label{LSKennebec_Penobscot}Least squares estimators used to fit the log-monthly flows data to a bivariate PVAR(1) model with $\nu=12$. Constraints on the parameters are given in the Table~(\ref{Kennebec_Penobscot}).}
\[
\hat{\bfPhi}(1) =
\left(
   \begin{array}{cc}
0.548     & 0.279    \\
0.153   &  0.358  \\
   \end{array}
\right),
\hat{\bfPhi}(2) =
\left(
   \begin{array}{cc}
 0.596  & 0.426  \\
 0.169  & 0.679 \\
   \end{array}
\right),
\hat{\bfPhi}(3) =
\left(
   \begin{array}{cc}
 0.739  & 0.346  \\
 -0.021  & 0.654 \\
   \end{array}
\right),
\hat{\bfPhi}(4) =
\left(
   \begin{array}{cc}
 0.285  & 0.565  \\
 0.149      & 0.142 \\
   \end{array}
\right),
\hat{\bfPhi}(5) =
\left(
   \begin{array}{cc}
 0.568  & 0.297  \\
 0.085      & 0.214 \\
   \end{array}
\right),
\hat{\bfPhi}(6) =
\left(
   \begin{array}{cc}
 0.287  & 0.088  \\
 0.139      &  0.247 \\
   \end{array}
\right),
\]
\[
\hat{\bfPhi}(7) =
\left(
   \begin{array}{cc}
0.441        & 0.065    \\
0            & 0.508  \\
   \end{array}
\right),
\hat{\bfPhi}(8) =
\left(
   \begin{array}{cc}
 0.183   & 0.345  \\
 -0.206  & 0.659 \\
   \end{array}
\right),
\hat{\bfPhi}(9) =
\left(
   \begin{array}{cc}
 0.156  & 0.318  \\
 -0.326 & 0.804 \\
   \end{array}
\right),
\hat{\bfPhi}(10) =
\left(
   \begin{array}{cc}
 0.766  & 0.164  \\
 -0.189 & 0.907 \\
   \end{array}
\right),
\hat{\bfPhi}(11) =
\left(
   \begin{array}{cc}
 0.452  & 0.515  \\
 0  & 0.734 \\
   \end{array}
\right),
\hat{\bfPhi}(12) =
\left(
   \begin{array}{cc}
 0.509  & 0.290  \\
 0.156  & 0.493 \\
   \end{array}
\right)
\]
\end{center}
\end{table}

The residual analysis was done and the portmanteau test statistics $Q_M^1(\nu)$, $Q^2_M(\nu)$ and $Q^3_M(\nu)$ were calculated. As far as short time series were used, we will interpret the results with caution. First, we apply portmanteau tests $Q_M^1(\nu)$ and $Q^2_M(\nu)$ for checking the hypothesis that the Kennebec and Penobscot log-monthly flows constitute a white noise. The $P$-values are reported in Table~\ref{pvalueQM_example2}. Since the $P$-values of these two tests are small, the white-noise hypothesis is rejected at the nominal level $\alpha = 1\%$. \cite{WGVM05} show that clear evidences are found for the existence of a nonlinear phenomenon of the variance behaviour, in the residual series from linear models fitted to daily and monthly streamflow processes of the upper Yellow River, China. The major cause of this effect is the seasonal variation in variance of the residual series. In contrast, the white-noise hypothesis is not rejected by our test $Q^3_M(\nu)$ since the statistic is not larger than the critical values.

\begin{table}[H]
\begin{center}
\scriptsize
\caption{$P$-values of the portmanteau test statistics used to check
a bivariate PVAR model of order $1$ with
$\nu=12$.}\label{pvalueQM_example2}
\setlength{\tabcolsep}{1.25mm}
\begin{tabular}{|c|cccccccccccc|}
\hline
& \multicolumn{12}{c|}{$Q^1_M(\nu)$}\\ 
\hline
\diagbox{$M$}{$\nu$} & 1 & 2 & 3 & 4 & 5 & 6 & 7 & 8 & 9 & 10 & 11 & 12 \\
\hline
  1 & n.a. & n.a. & n.a. & n.a. & n.a. & n.a. & n.a. & n.a. & n.a. & n.a. & n.a. & n.a. \\
  2 & 0.03 & 0.05 & 0.04 & 0.08 & \textbf{0.00} & 0.20 & 0.96 & 0.53 & 0.20 & 0.07 & 0.56 & 0.58 \\ 
  3 & 0.03 & 0.09 & 0.04 & 0.22 & \textbf{0.00} & 0.53 & 0.83 & 0.10 & 0.30 & 0.13 & 0.73 & 0.40 \\ 
  6 & 0.24 & 0.08 & 0.14 & 0.09 & \textbf{0.00} & 0.56 & 0.24 & 0.28 & 0.58 & 0.01 & 0.45 & 0.17 \\ 
  8 & 0.42 & 0.16 & 0.26 & 0.16 & \textbf{0.00} & 0.65 & 0.06 & 0.13 & 0.12 & 0.01 & 0.18 & 0.07 \\ 
 10 & 0.58 & 0.26 & 0.40 & 0.07 & 0.01 & 0.60 & 0.03 & 0.12 & 0.22 & 0.02 & 0.16 & 0.28 \\ 
\hline
& \multicolumn{12}{c|}{$Q^2_M(\nu)$}\\
\hline
  1 & \textbf{0.00} & 0.18 & 0.06 & 0.64 & 0.09 & 0.03 & 0.87 & 0.23 & 0.46 & 0.06 & 0.70 & 0.46 \\ 
  2 & 0.04 & 0.10 & 0.07 & 0.13 & \textbf{0.00} & 0.22 & 0.96 & 0.58 & 0.24 & 0.08 & 0.63 & 0.67 \\ 
  3 & 0.03 & 0.12 & 0.06 & 0.25 & \textbf{0.00} & 0.51 & 0.83 & 0.10 & 0.32 & 0.13 & 0.73 & 0.43 \\ 
  6 & 0.23 & 0.08 & 0.14 & 0.08 & \textbf{0.00} & 0.53 & 0.23 & 0.26 & 0.56 & 0.01 & 0.44 & 0.16 \\ 
  8 & 0.41 & 0.15 & 0.24 & 0.15 & \textbf{0.00} & 0.61 & 0.06 & 0.12 & 0.11 & 0.01 & 0.17 & 0.06 \\ 
 10 & 0.56 & 0.25 & 0.38 & 0.06 & 0.01 & 0.56 & 0.03 & 0.10 & 0.19 & 0.02 & 0.15 & 0.25 \\ 
\hline
& \multicolumn{12}{c|}{$Q^3_M(\nu)$}\\
\hline
  1 & 0.02 & 0.13 & 0.22 & 0.55 & 0.13 & 0.03 & 0.71 & 0.35 & 0.56 & 0.11 & 0.65 & 0.48 \\ 
  2 & 0.12 & 0.08 & 0.13 & 0.11 & 0.01 & 0.22 & 0.91 & 0.65 & 0.37 & 0.12 & 0.61 & 0.71 \\ 
  3 & 0.13 & 0.10 & 0.09 & 0.27 & 0.02 & 0.35 & 0.64 & 0.12 & 0.46 & 0.14 & 0.68 & 0.56 \\ 
  6 & 0.31 & 0.11 & 0.16 & 0.16 & 0.04 & 0.36 & 0.27 & 0.24 & 0.56 & 0.02 & 0.40 & 0.27 \\ 
  8 & 0.42 & 0.27 & 0.24 & 0.19 & 0.05 & 0.44 & 0.14 & 0.12 & 0.29 & 0.02 & 0.22 & 0.16 \\ 
 10 & 0.54 & 0.36 & 0.35 & 0.14 & 0.08 & 0.41 & 0.11 & 0.09 & 0.37 & 0.01 & 0.23 & 0.34 \\ 
\hline
\end{tabular}
\end{center}
\end{table}

\begin{figure}[h] 
  \begin{minipage}[b]{0.48\linewidth}
    \includegraphics[width=.95\linewidth]{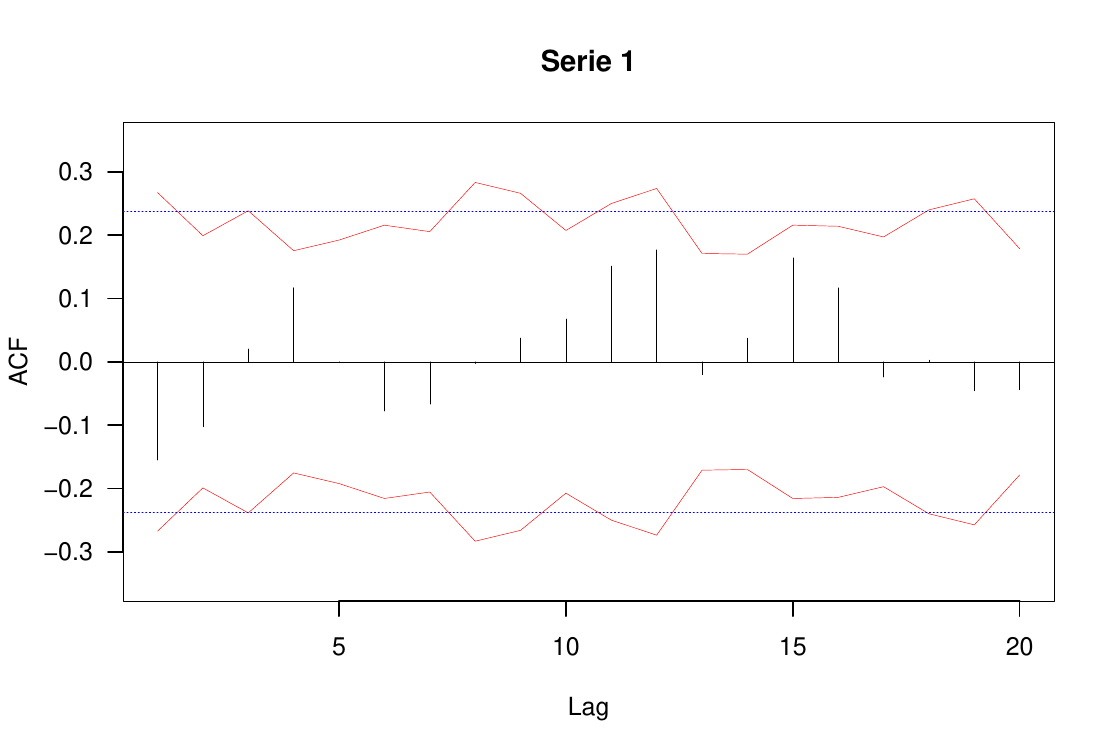} 
  \end{minipage} 
  \begin{minipage}[b]{0.48\linewidth}
    \includegraphics[width=.95\linewidth]{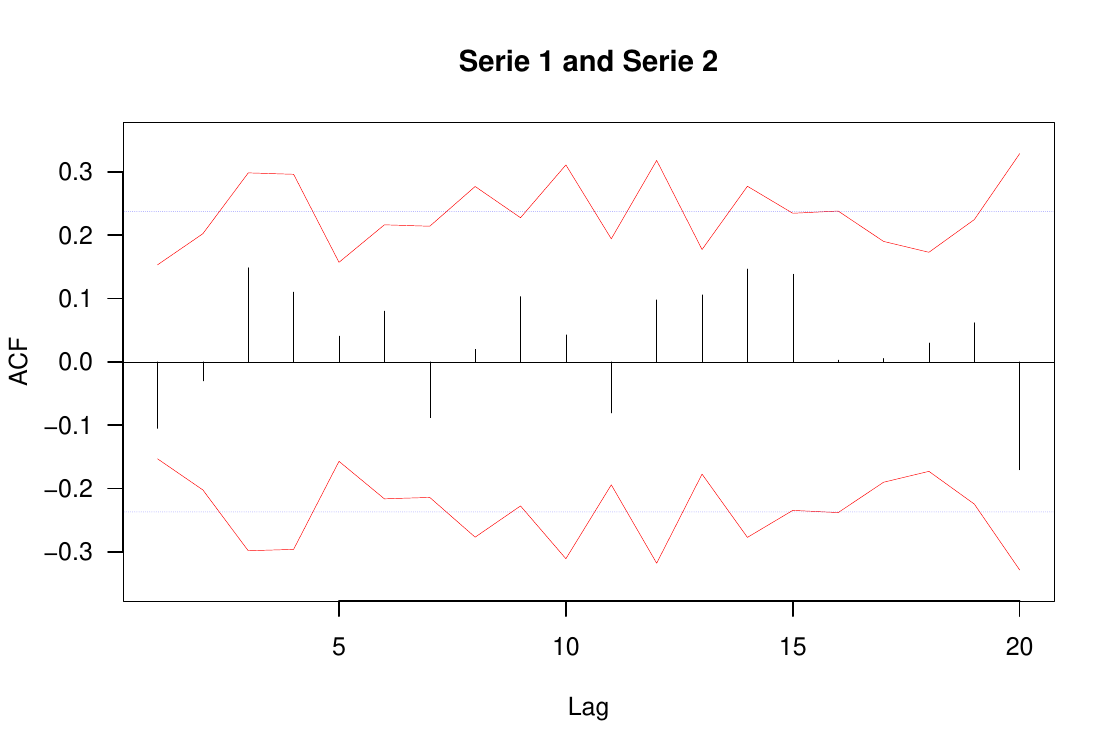} 
  \end{minipage} 
  \begin{minipage}[b]{0.48\linewidth}
    \includegraphics[width=.95\linewidth]{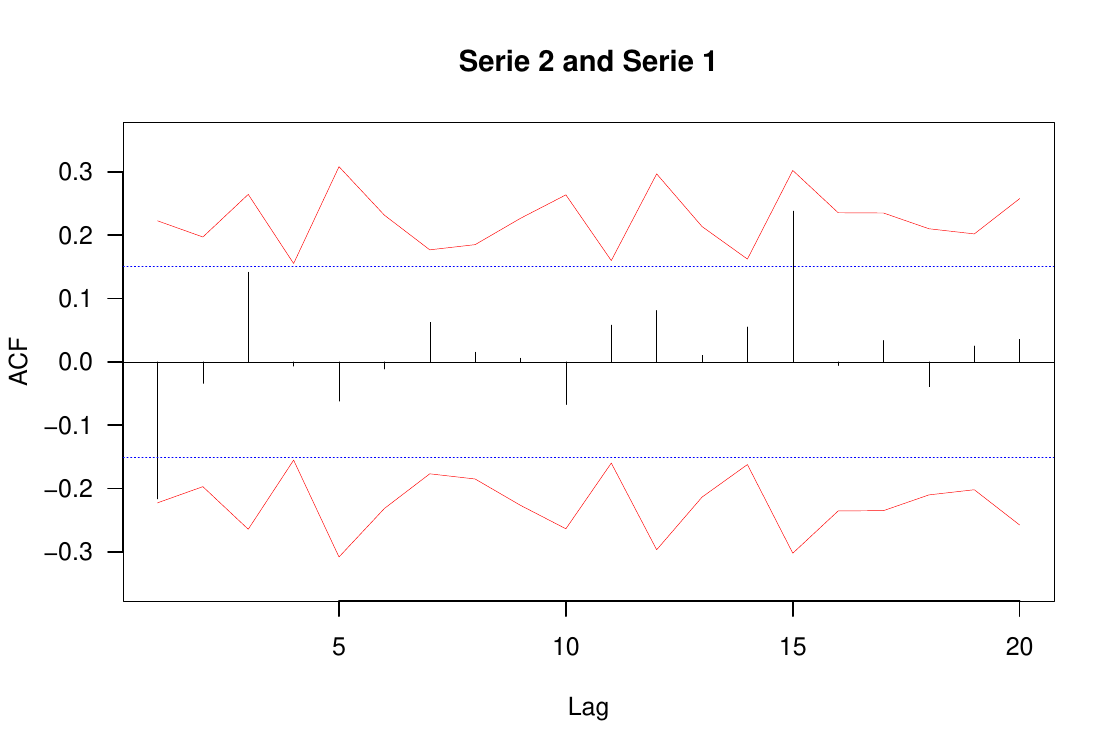} 
  \end{minipage}
  \hfill
  \begin{minipage}[b]{0.48\linewidth}
    \includegraphics[width=.95\linewidth]{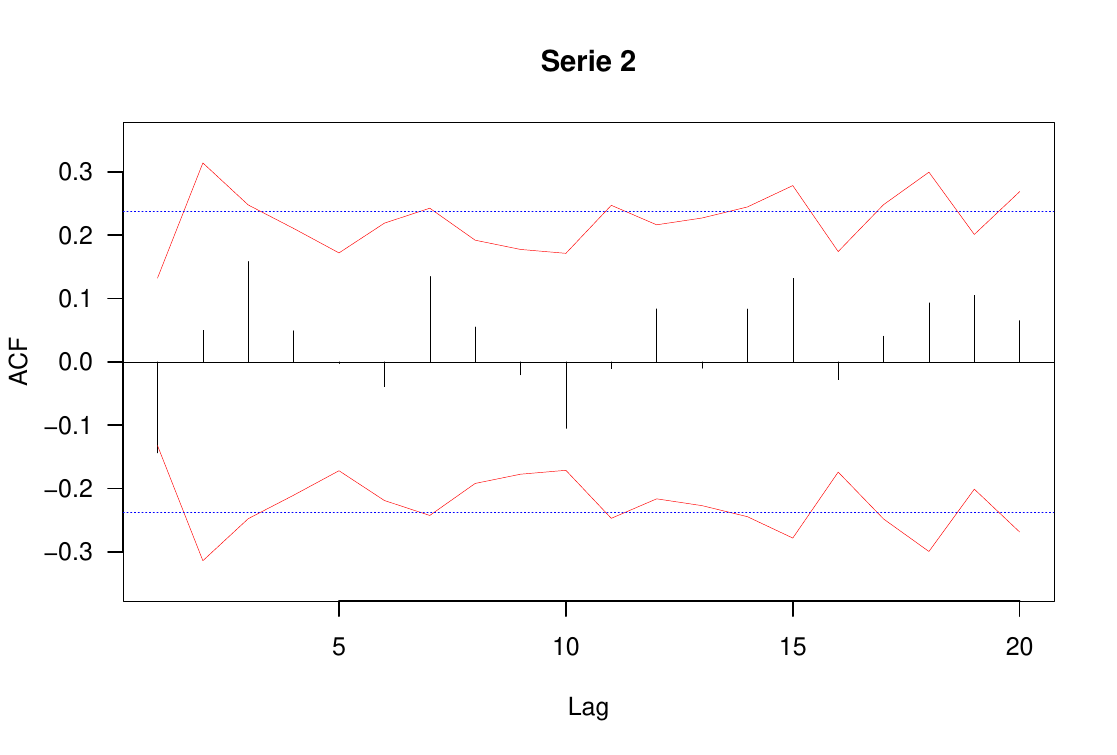} 
  \end{minipage} 
  \caption{Autocorrelations of the log-monthly flows of Kennebec and Penobscot river for the first period (serie 1 corresponds to Kennebec).}\label{autocorrfig_river} 
\end{figure}

In view of Figure~\ref{autocorrfig_river}, the diagnostic checking of residuals does not indicate any inadequacy, approximately $99\%$ of the sample autocorrelations should lie between the bands shown as solid lines (red color).

\begin{figure}[h]
\centering
\includegraphics[width=.7\textwidth]{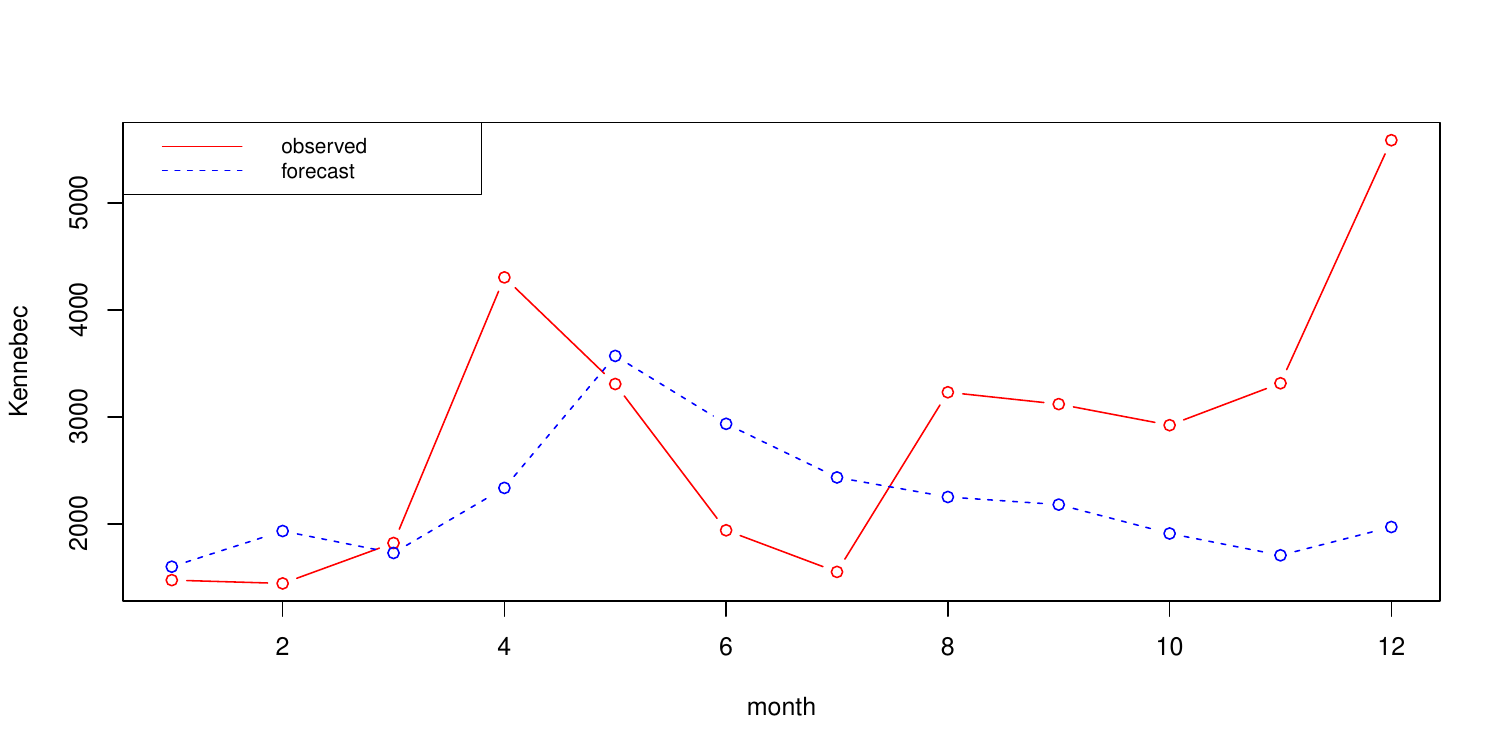}
\includegraphics[width=.7\textwidth]{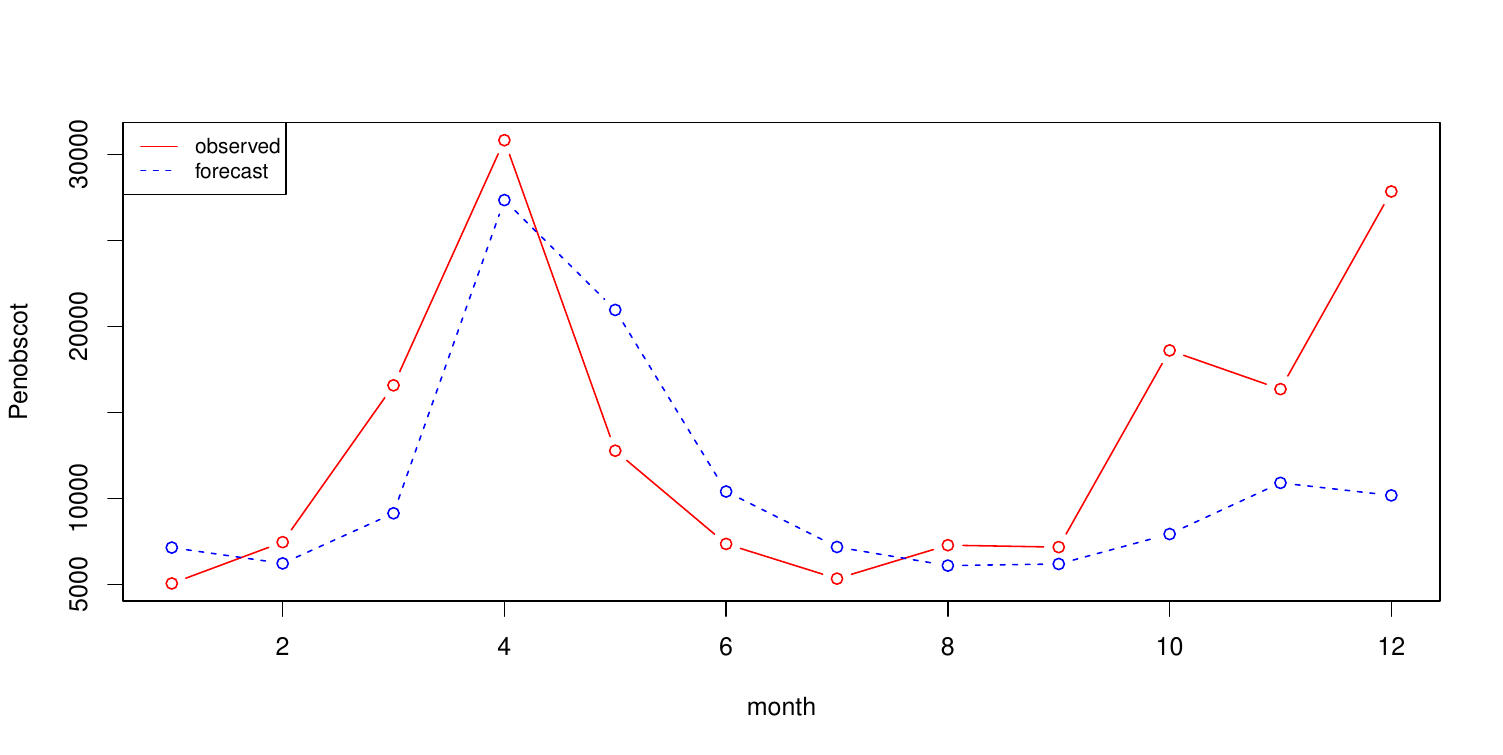}
\caption{Twelve-month forecasts (dashed line) based on $118$ years of monthly data. The last year of the data set (solid line) was not used in the forecast.}\label{prediction_river} 
\end{figure}

The weak PVAR model for Kennebec and Penobscot river data is used to generate
one-step-ahead forecasts for the monthly flow series (Figure~\ref{prediction_river}). The differences in the forecasts can be explained by a possible change around $2002$.

\begin{figure}[h]
\centering
\includegraphics[width=.7\textwidth]{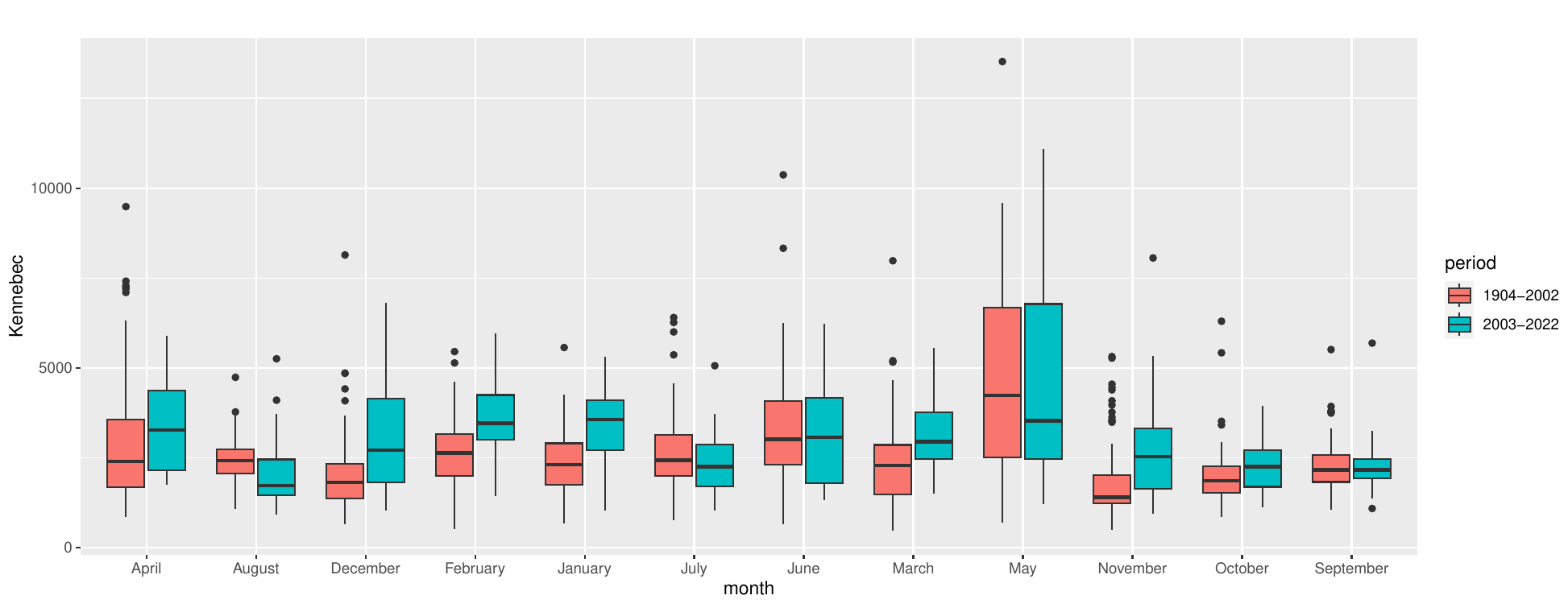}
\includegraphics[width=.7\textwidth]{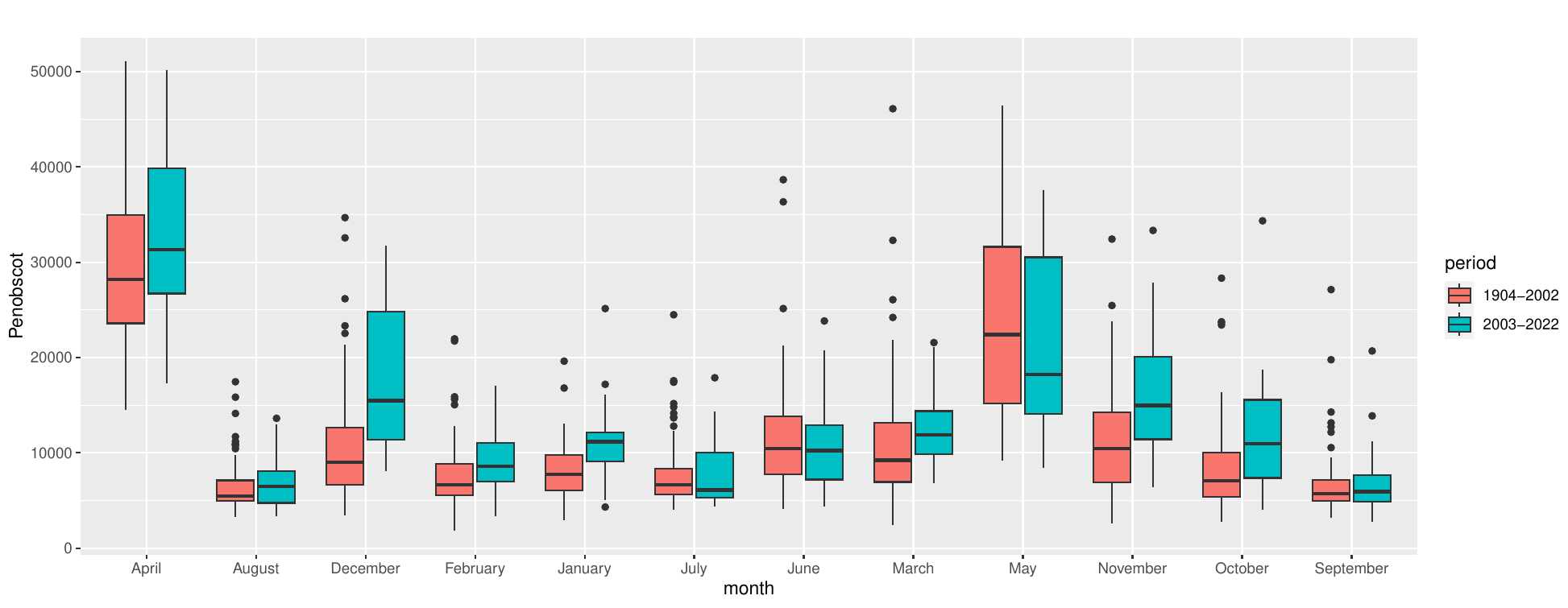}
\caption{Box-plots of the monthly flows before and after $2002$.}\label{boxplot_river} 
\end{figure}

\begin{figure}[h]
\centering
\includegraphics[width=.7\textwidth]{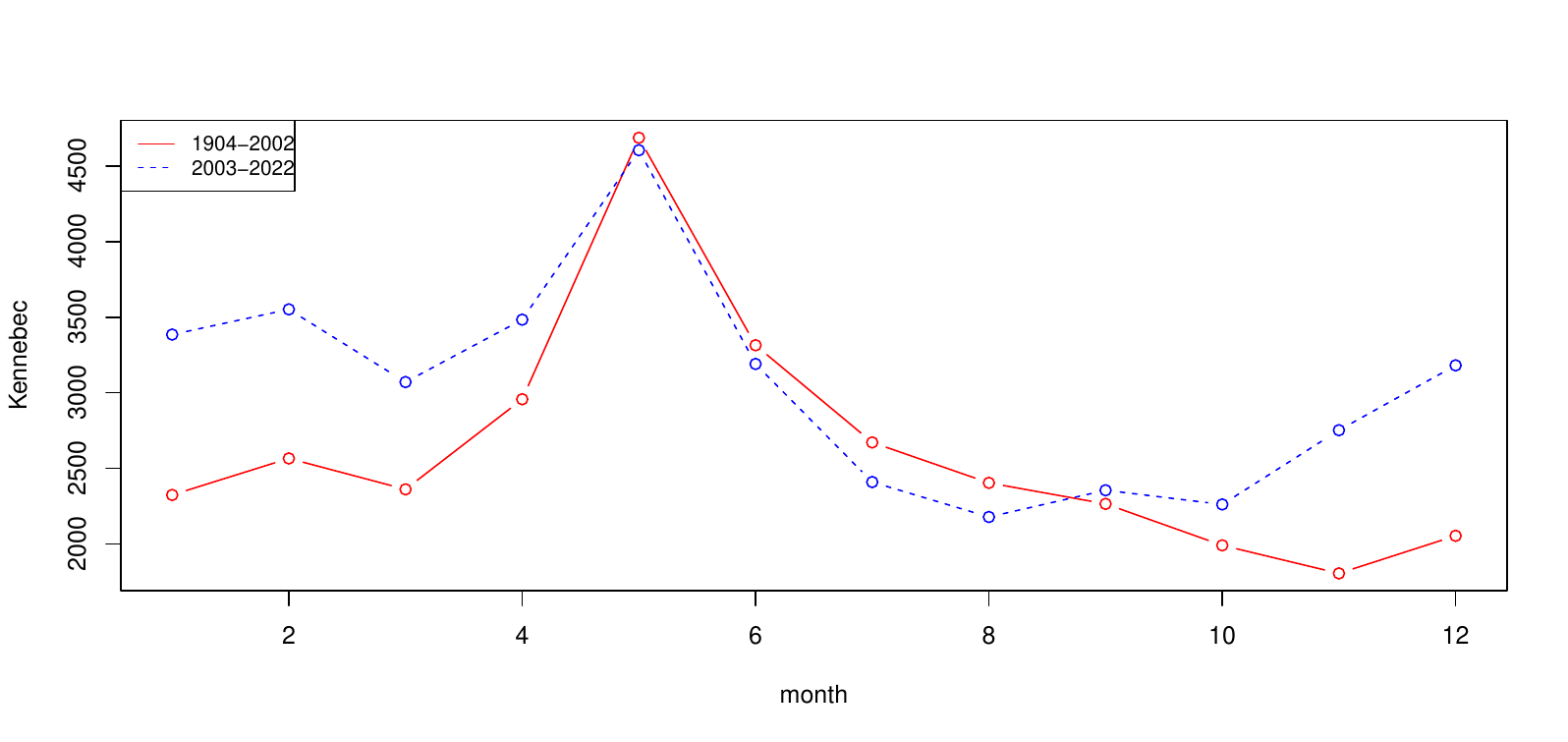}
\includegraphics[width=.7\textwidth]{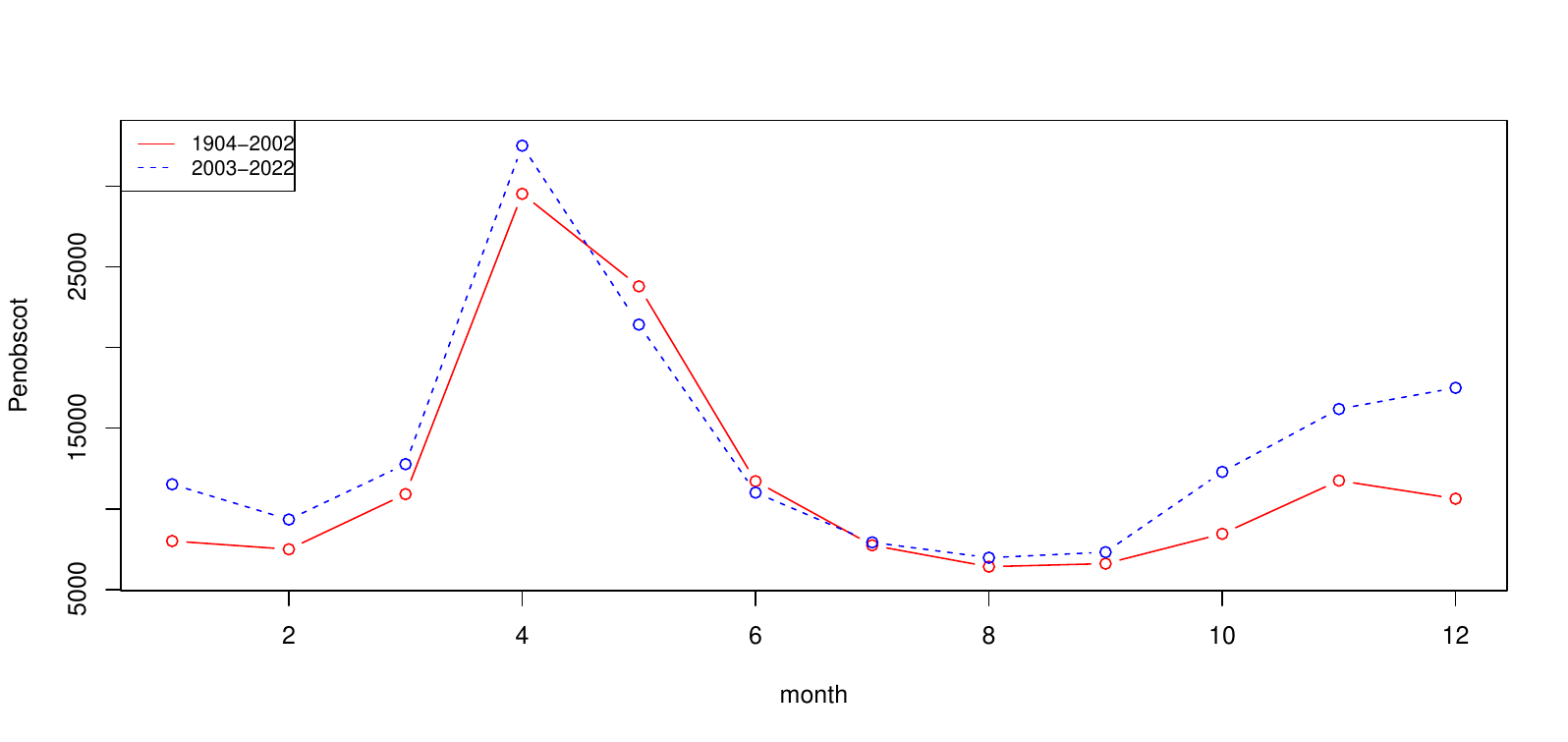}
\caption{Means of the monthly flows before and after $2002$.}\label{means_river} 
\end{figure}

To ascertain the changes in the time series data, we calculate the $12$ seasonal means for all the years up to $2002$, and from $2002$ onwards. We plot the seasonal means for each segment in Figure~\ref{means_river}. We observe that the seasonal means increased from October to March after $2002$, compared with the previous period.
\section{Conclusion}\label{conclusion}
\noindent In this paper we derive the asymptotic distribution of residual empirical autocovariances and autocorrelations in the class of PVAR models under weak assumptions on the noise. We establish the asymptotic distribution of the modified Ljung-Box portmanteau test statistics for PVAR models with non-independent innovations. This asymptotic distribution is quite different from the usual chi-squared approximation used under iid assumptions on the noise. Therefore, the proposed test is more difficult to implement because its critical values have to be computed from the data, whereas those of the standard versions are simply given in a $\chi^2$-table. For this reason the use of our proposed test for short series is not recommended. In simulation experiments, our proposed test has a satisfactory behavior for series of length $N\geq 1000$ under dependent errors.

The proposed test statistics were illustrated in a small simulation study. A comparison with the statistics proposed by~\citet{UD09} and by~\citet{DL13} was made. From our simulation experiments, we demonstrated that the proposed portmanteau test statistics have reasonable finite sample performance, at least for the models considered in our study. Therefore, the proposed test performs adequately in terms of level in the strong PVAR setting in both the unconstrained and constrained cases when the number of observations is large, but less so when fewer observations are available, likely as a downstream consequence of needing to estimate the additional quantities. Under non-independent errors, it appears that the test statistics proposed by~\citet{UD09} or by~\citet{DL13} are generally unreliable, overrejecting severally, while our proposed test statistics offers satisfactory levels in most cases. Even for independent errors, our test may be preferable to that proposed by~\citet{UD09}, when the number $M$ of autocorrelations is small. 

We applied our methodology using real data. 
From our analysis, the strong PVAR model was clearly rejected. Using our approach, a reasonable model for modeling the above data has been proposed and checked, using the portmanteau test statistics described in this paper.

\appendix
\section{Appendix : Proofs of the main results}\label{app}
\subsection{Proof of Proposition \ref{loijointeUn}}
To establish the asymptotic distribution of
$N^{1/2}\bfc_{\hat{\bfepsilon}}(\nu)$,
we first note that:
\begin{eqnarray}\label{epseps}
  \bfc_{\bfepsilon}(l; \nu) = \vec\left\{N^{-1}\sum_{n=l}^{N-1} (\bfepsilon_{ns+\nu}  \bfepsilon_{ns+\nu-l}^{\top}) \right\}
                            = N^{-1} \sum_{n=l}^{N-1} (\bfepsilon_{ns+\nu-l} \otimes \bfepsilon_{ns+\nu}),
\end{eqnarray}
where we used the fact that $\vec(\bfa \bfb^{\top}) = \bfb \otimes \bfa$, where $\bfa$ and $\bfb$
are two arbitrary vectors. Note that the mixing assumption \textbf{(A3)} will entail the asymptotic normality of $\bfc_{{\bfepsilon}}(\nu)$ and $\bfc_{{\bfepsilon}}$. Since  $\{\bfepsilon_t \}$ represents a weak periodic white noise, it follows that $\mathbb{E}\{ \bfc_{\bfepsilon}(l; \nu) \} = \bfzero$. Thus for $l,l'\geq1$ and $\nu,\nu'=1,\dots,s$ we have:
\begin{eqnarray*}
  N\cov\{ \bfc_{\bfepsilon}(l; \nu), \bfc_{\bfepsilon}(l'; \nu') \} &=& N\mathbb{E}\{ \bfc_{\bfepsilon}(l; \nu) \bfc_{\bfepsilon}^{\top}(l'; \nu') \}, \\
                            &\to&\bfV_{\nu\nu'}(l,l')=: \sum_{h=-\infty}^{\infty}
                            \mathbb{E}(\bfepsilon_{ns+\nu-l} \bfepsilon_{(n-h)s+\nu'-l'}^{\top} \otimes\bfepsilon_{ns+\nu} \bfepsilon_{(n-h)s+\nu'}^{\top}),\quad\text{as}\quad N\to\infty,
\end{eqnarray*}
by stationary and  the dominated convergence theorem.
The existence of the matrix $\bfV_{\nu\nu'}(l,l')$ is a consequence of Assumption \textbf{(A3)} and the \cite{D68} inequality.

Applying the central limit theorem (CLT in what follows) for mixing processes (see \cite{H84}) to the vector\\%
$\left((\bfepsilon_{ns+\nu-1} \otimes \bfepsilon_{ns+\nu})^\top,\dots,(\bfepsilon_{ns+\nu-M} \otimes \bfepsilon_{ns+\nu})^\top\right)^\top$, we obtain:
\begin{equation}
\label{cepsilonnu}
 N^{1/2} \bfc_{\bfepsilon}(\nu) \cd N_{d^2M}(\bfzero, \bfV_{\nu\nu}),
\end{equation}
where
$\bfV_{\nu\nu}$
corresponds to the
$(d^2 M) \times (d^2 M)$
  matrix given by:
\begin{equation*}
  \bfV_{\nu\nu} = \left[\bfV_{\nu\nu}(l,l')\right]_{1\leq l,l'\leq M}.
\end{equation*}
Applying the CLT  for mixing processes (see \cite{H84}) to (\ref{betahatnu}), we directly obtain:
\begin{equation*}
 N^{1/2}\{ \hat{\bfbeta}(\nu) - \bfbeta(\nu) \} \cd N_{d^2p(\nu)}(\bfzero, \bfTheta_{\nu\nu}),
\end{equation*}
which gives \eqref{th1c}.
From \eqref{betahatnuOmega} and for $\nu,\nu'=1,\dots,s$ we deduce that:
\[
 \lim_{N \rightarrow \infty}
 N \cov\{ \hat{\bfbeta}(\nu)  - \bfbeta(\nu),  \bfc_{\bfepsilon}(\nu') \} =
 \{ \bfOmega^{-1}(\nu) \otimes \bfI_d \}
 \lim_{N \rightarrow \infty}\mathbb{E}
 \left[ \vec\{ \bfE(\nu)\bfX^{\top}(\nu) \} \bfc_{\bfepsilon}^{\top}(\nu') \right].
\]
Vectorizing,
\begin{eqnarray}\label{EX}
  \vec\{ \bfE(\nu) \bfX^{\top}(\nu) \} = \sum_{n=0}^{N-1}\vec\{ \bfepsilon_{ns+\nu} \bfX^{\top}_n(\nu) \}  = \sum_{n=0}^{N-1} \bfX_n(\nu) \otimes \bfepsilon_{ns+\nu}.
\end{eqnarray}
Consequently, in view of~(\ref{epseps}) and~(\ref{EX}) we have:
\begin{eqnarray*}
 \lim_{N \rightarrow \infty}\mathbb{E} \left[ \vec\{ \bfE(\nu) \bfX^{\top}(\nu) \} \bfc^{\top}_{\bfepsilon}(l;\nu') \right] &=&
 \lim_{N \rightarrow \infty}  N^{-1}\sum_{n=0}^{N-1}\sum_{n'=l}^{N-1}\mathbb{E}\left[ \bfX_n(\nu) \bfepsilon_{n's+\nu'-l}^{\top}
                             \otimes \bfepsilon_{ns+\nu}\bfepsilon_{n's+\nu'}^\top\right], \\
&=& \lim_{N \rightarrow \infty}  N^{-1}\sum_{h=1-N}^{N-1}(N-|h|)\mathbb{E}\left[ \bfX_n(\nu) \bfepsilon_{(n-h)s+\nu'-l}^{\top}
                             \otimes \bfepsilon_{ns+\nu}\bfepsilon_{(n-h)s+\nu'}^\top\right], \\
&=& \sum_{h=-\infty}^{\infty}\mathbb{E}\left[ \bfX_n(\nu) \bfepsilon_{(n-h)s+\nu'-l}^{\top}
                             \otimes \bfepsilon_{ns+\nu}\bfepsilon_{(n-h)s+\nu'}^\top\right],
\end{eqnarray*}
by the stationary and dominated convergence theorem. Now, collecting the $\{ d^2 p(\nu) \} \times d^2$ matrix $\bfG^{\mathrm{U}}_{\nu\nu'}(\cdot,l)$, for $l=1,\ldots,M$, in a
$\{ d^2 p(\nu) \} \times (Md^2)$ matrix leads to the expression:
\begin{equation*}
 \bfG^{\mathrm{U}}_{\nu\nu'}=: \left[ \bfG^{\mathrm{U}}_{\nu\nu'}(\cdot,1), \ldots, \bfG^{\mathrm{U}}_{\nu\nu'}(\cdot,M) \right].
\end{equation*}
where
\begin{align*}
\bfG^{\mathrm{U}}_{\nu\nu'}(\cdot,l) = :\sum_{h=-\infty}^{\infty}\mathbb{E}\left[\bfOmega^{-1}(\nu) \bfX_n(\nu) \bfepsilon_{(n-h)s+\nu'-l}^{\top}
                             \otimes \bfepsilon_{ns+\nu}\bfepsilon_{(n-h)s+\nu'}^\top\right].
\end{align*}
The joint asymptotic normality of $N^{1/2} \{ \hat{\bfbeta}^\top - \bfbeta^\top,\bfc_{\bfepsilon}^\top \}^\top$ follows using the same kind of manipulations as those for a single season $\nu$. We also have
\begin{eqnarray}\nonumber
  N^{1/2} \left( \hat{\bfbeta}^\top - \bfbeta^\top,\bfc_{\bfepsilon}^\top \right)^\top
  &\cd& \mathcal{N}_{sd^2(p(\nu)+M)}\left(\bfzero, \left(\begin{array}{cc}
                           \bfTheta & \bfG^{\mathrm{U}}\\
                           \left(\bfG^{\mathrm{U}}\right)^{\top} & \bfV
                           \end{array}
                     \right) \right),
\end{eqnarray}
where the asymptotic covariance matrix $\bfTheta$ is given by \eqref{th1d}, $\bfV = \left[\bfV_{\nu\nu'}(M)\right]_{1\leq \nu,\nu'\leq s}$ and $\bfG^{\mathrm{U}} = \left[\bfG^{\mathrm{U}}_{\nu\nu'}\right]_{1\leq \nu,\nu'\leq s}$.

It is clear that the existence of the above matrices  is a consequence of Assumption \textbf{(A3)} and the \cite{D68} inequality, which completes the proof of Proposition~\ref{loijointeUn}.$\hfill\square$

\subsection{Proof of Theorem~\ref{loi_res_gamUn}}\label{proof_loi_res}
The proof is divided in two steps.
\subsubsection{Step 1: Taylor's expansion of $\sqrt{N}\bfc_{\hat{\bfepsilon}}(\nu)$ and $\sqrt{N}\bfr_{\hat{\bfepsilon}}(\nu)$}
The aim of this step is to prove \eqref{hat_gamma} and \eqref{hat_rho}.

By expanding $\bfc_{\dot{\bfepsilon}}(\nu)$ in a Taylor expansion around ${\bfbeta}(\nu)$ and evaluating at the point $\dot{\bfbeta}(\nu) = \hat{\bfbeta}(\nu)$, for $\nu=1,\dots,s$ we obtain the following development:
\[
\bfc_{\hat{\bfepsilon}}(\nu) = \bfc_{\bfepsilon}(\nu) +
                    \frac{\partial \bfc_{{\bfepsilon}}(\nu)}{\partial \bfbeta^{\top}(\nu)}
                    \left\{
                    \hat{\bfbeta}(\nu) - \bfbeta(\nu)
                    \right\} + \mathrm{o}_{\mathbb{P}}(N^{-1/2}),
\]
where $\partial \bfc_{{\bfepsilon}}(\nu)/\partial \bfbeta^{\top}(\nu)$ corresponds to a $(d^2 M) \times \{ d^2p(\nu) \}$ matrix satisfying
\[
\left( \partial \bfc_{{\bfepsilon}}(\nu)/\partial \bfbeta^{\top}(\nu) \right)^{\top} =
\left(  (\partial \bfc_{{\bfepsilon}}(1;\nu)/\partial \bfbeta^{\top}(\nu))^{\top}, \ldots,
        (\partial \bfc_{{\bfepsilon}}(M;\nu)/\partial \bfbeta^{\top}(\nu))^{\top} \right).
\]
We have
\[
  \frac{\partial \bfc_{{\bfepsilon}}(\nu)}{\partial \bfbeta^{\top}(\nu)} \cp \bfUpsilon_M(\nu)=-\mathbb{E}\left\{\left(\begin{array}{c}
\bfepsilon_{ns+\nu-1}\\\vdots\\\bfepsilon_{ns+\nu-M}\end{array}\right)\otimes
\bfX_n^\top(\nu)\otimes \bfI_d\right\},
\]
by the ergodic theorem. The Taylor expansion \eqref{hat_gamma} of $N^{1/2}\bfc_{\hat{\bfepsilon}}(\nu)$ is proved.

We now come back to the vector $N^{1/2}\bfr_{\hat{\bfepsilon}}(\nu)$.
From a Taylor expansion about $\bfbeta$ of $\bfC_{\hat{\bfepsilon}}(0; \nu)$ we
have, $\bfC_{\hat{\bfepsilon}}(0; \nu)=\bfC_{{\bfepsilon}}(0; \nu)+\mathrm{O}_{\mathbb{P}}(N^{-1/2}).$ Moreover,
$\sqrt{N}(\bfC_{{\bfepsilon}}(0; \nu)-\mathbb{E} (\bfC_{{\bfepsilon}}(0; \nu)))=\mathrm{O}_{\mathbb{P}}(1)$ by the CLT for
mixing processes (see \cite{H84}) to the process $(\bfepsilon_{ns+\nu}\otimes\bfepsilon_{ns+\nu})_{n\in\mathbb{Z}}$. Thus $\sqrt{N}\{\bfD_{\hat{\bfepsilon}}(\nu)  \otimes
                                    \bfD_{\hat{\bfepsilon}}(\nu)-
                                    \bfD^{\ast}_{{\bfepsilon}}(\nu)  \otimes
                                    \bfD^{\ast}_{{\bfepsilon}}(\nu)\}=\mathrm{O}_{\mathbb{P}}(1)$, where $\bfD^{\ast}_{{\bfepsilon}}(\nu) =
\diag\left( \bfGamma_{{\bfepsilon},11}^{1/2}(0; \nu),\ldots,\bfGamma_{{\bfepsilon},dd}^{1/2}(0; \nu) \right)$. Using (\ref{chatgammaU}) and the ergodic Theorem,
we obtain
\begin{align*}
N&\left(\bfr_{\dot{\bfepsilon}}(h; \nu) - \vec\left\{ \bfD_{{\bfepsilon}}^{\ast -1}(\nu)\bfC_{\dot{\bfepsilon}}(h; \nu)
                                    \bfD_{{\bfepsilon}}^{\ast -1}(\nu-h) \right\}\right)
                                    =N\left\{\left( \bfD_{\dot{\bfepsilon}}^{-1}(\nu-h)  \otimes
                                    \bfD_{\dot{\bfepsilon}}^{-1}(\nu) \right)-
                                     \left( \bfD_{{\bfepsilon}}^{\ast -1}(\nu-h)  \otimes
                                    \bfD_{{\bfepsilon}}^{\ast -1}(\nu) \right)
                                     \right\}\bfc_{\dot{\bfepsilon}}(h;\nu)\\&=\left( \bfD_{\dot{\bfepsilon}}(\nu-h)  \otimes
                                    \bfD_{\dot{\bfepsilon}}(\nu) \right)^{-1}\sqrt{N}\left\{ \bfD_{{\bfepsilon}}^{\ast}(\nu-h)  \otimes
                                    \bfD_{{\bfepsilon}}^{\ast}(\nu)
                                    -\bfD_{\dot{\bfepsilon}}(\nu-h)  \otimes
                                    \bfD_{\dot{\bfepsilon}}(\nu)\right\}\left( \bfD_{{\bfepsilon}}^{\ast}(\nu-h)  \otimes
                                    \bfD_{{\bfepsilon}}^{\ast}(\nu) \right)^{-1}\sqrt{N}\bfc_{\dot{\bfepsilon}}(h;\nu)\\&=\mathrm{O}_{\mathbb{P}}(1),
\end{align*}
by applying the CLT for mixing processes (see \cite{H84}) to the process $(\bfepsilon_{ns+\nu-l}\otimes\bfepsilon_{ns+\nu})_{n\in\mathbb{Z}}$, $l\geq1$.

In the previous equalities, we also use $\vec(\bfA\bfB\bfC)=(\bfC^\top\otimes\bfA)\vec(\bfB)$ and $(\bfA\otimes \bfB)^{-1}=\bfA^{-1}\otimes \bfB^{-1}$ when $\bfA$ and $\bfB$ are invertible.

It follows that:
\begin{align*}\nonumber
\sqrt{N}\bfr_{\hat{\bfepsilon}}(\nu)
&=\sqrt{N}\left(\left\{\left( \bfD_{\hat{\bfepsilon}}^{-1}(\nu-1)  \otimes
                                    \bfD_{\hat{\bfepsilon}}^{-1}(\nu) \right)
                                     \bfc_{\hat{\bfepsilon}}(1;\nu)\right\}^\top,\dots,\left\{\left( \bfD_{\hat{\bfepsilon}}^{-1}(\nu-M)  \otimes
                                    \bfD_{\hat{\bfepsilon}}^{-1}(\nu) \right)
                                     \bfc_{\hat{\bfepsilon}}(M;\nu)\right\}^\top
\right)^\top
\\ \nonumber&=\diag\left[ \bfD_{\hat{\bfepsilon}}^{-1}(\nu-1)  \otimes
                                    \bfD_{\hat{\bfepsilon}}^{-1}(\nu),\dots,\bfD_{\hat{\bfepsilon}}^{-1}(\nu-M)  \otimes
                                    \bfD_{\hat{\bfepsilon}}^{-1}(\nu) \right]
\sqrt{N}\bfc_{\hat{\bfepsilon}}(\nu) \\ \label{hat_rhobis}&=
\diag\left[ \bfD_{{\bfepsilon}}^{\ast -1}(\nu-1)  \otimes
                                    \bfD_{{\bfepsilon}}^{\ast -1}(\nu),\dots,\bfD_{{\bfepsilon}}^{\ast -1}(\nu-M)  \otimes
                                    \bfD_{{\bfepsilon}}^{\ast -1}(\nu) \right]
\sqrt{N}\bfc_{\dot{\bfepsilon}}(\nu)+\mathrm{o}_{\mathbb{P}}(1),
\end{align*}
and the Taylor expansion \eqref{hat_rho} of $\sqrt{N}\bfr_{\hat{\bfepsilon}}(\nu)$ is proved. The proof of the Taylor expansion  of $\sqrt{N}\bfr_{\hat{\bfepsilon}}$ follows using the same kind of arguments as the Taylor expansion \eqref{hat_rho} of $\sqrt{N}\bfr_{\hat{\bfepsilon}}(\nu)$.
This ends our first step.

The next step deals with the asymptotic distributions of $\sqrt{N}\bfc_{\hat{\bfepsilon}}(\nu)$ and $\sqrt{N}\bfr_{\hat{\bfepsilon}}(\nu)$.
\subsubsection{Step 2: asymptotic distributions of $\sqrt{N}\bfc_{\hat{\bfepsilon}}(\nu)$ and $\sqrt{N}\bfr_{\hat{\bfepsilon}}(\nu)$}
In view of \eqref{hat_gamma},
$N^{1/2}\bfc_{\hat{\bfepsilon}}(\nu)$
and
$N^{1/2} \left(  \bfc_{\bfepsilon}(\nu) +
                    \bfUpsilon_M(\nu)
                    \left\{
                    \hat{\bfbeta}(\nu) - \bfbeta(\nu)
                    \right\}\right)$
have the same asymptotic distribution.
The  joint asymptotic distribution of $N^{1/2} \left( \hat{\bfbeta}(\nu) - \bfbeta(\nu)\right)$
and $N^{1/2} \bfc_{\bfepsilon}(\nu)$ shows that $N^{1/2} \bfc_{\hat{\bfepsilon}}(\nu)$ has a limiting normal distribution with mean zero and covariance matrix
\begin{align*}
\lim_{N\rightarrow\infty}\mathrm{Var}\left(\sqrt{N} \bfc_{\hat{\bfepsilon}}(\nu)\right)&=\lim_{N\rightarrow\infty}\mathrm{Var}\left(\sqrt{N} \bfc_{\bfepsilon}(\nu)\right)+\bfUpsilon_M(\nu)\lim_{\rightarrow\infty}\mathrm{Var}\left(\sqrt{N}\left( \hat{\bfbeta}(\nu) - \bfbeta(\nu)\right)\right)\bfUpsilon^{\top}_M(\nu)\\
&\qquad+\bfUpsilon_M(\nu)\lim_{N\rightarrow\infty}\mathrm{Cov}\left(\sqrt{N}\left( \hat{\bfbeta}(\nu) - \bfbeta(\nu)\right),\sqrt{N} \bfc_{\bfepsilon}(\nu)\right)\\&\qquad+\lim_{N\rightarrow\infty}\mathrm{Cov}\left(\sqrt{N} \bfc_{\bfepsilon}(\nu),\sqrt{N}\left( \hat{\bfbeta}(\nu) - \bfbeta(\nu)\right)\right)\bfUpsilon^{\top}_M(\nu)\\
&=\bfV_{\nu\nu}+\bfUpsilon_M(\nu)\bfTheta_{\nu\nu}\bfUpsilon^{\top}_M(\nu)+\bfUpsilon_M(\nu)\bfG^{\mathrm{U}}_{\nu\nu}+\left(\bfG^{\mathrm{U}}_{\nu\nu}\right)^{\top}\bfUpsilon^{\top}_M(\nu)=:\bfDelta^{\mathrm{U}}_{\nu\nu},
\end{align*}
which prove \eqref{chatgammaU}.
We also deduce the asymptotic covariances for $\nu\neq\nu'$
\begin{align*}
\lim_{N\rightarrow\infty}\mathrm{Cov}\left(\sqrt{N} \bfc_{\hat{\bfepsilon}}(\nu),\sqrt{N} \bfc_{\hat{\bfepsilon}}(\nu')\right)&=
\lim_{N\rightarrow\infty}\mathrm{Cov}\left(\sqrt{N} \bfc_{\bfepsilon}(\nu),\sqrt{N} \bfc_{\bfepsilon}(\nu')\right)\\&\quad+\bfUpsilon_M(\nu)\lim_{N\rightarrow\infty}\mathrm{Cov}\left(\sqrt{N}\left( \hat{\bfbeta}(\nu) - \bfbeta(\nu)\right),\sqrt{N}\left( \hat{\bfbeta}(\nu') - \bfbeta(\nu')\right)\right)\bfUpsilon^{\top}_M(\nu')\\
&\quad+\bfUpsilon_M(\nu)\lim_{N\rightarrow\infty}\mathrm{Cov}\left(\sqrt{N}\left( \hat{\bfbeta}(\nu) - \bfbeta(\nu)\right),\sqrt{N} \bfc_{\bfepsilon}(\nu')\right)
\\&\quad+\lim_{N\rightarrow\infty}\mathrm{Cov}\left(\sqrt{N} \bfc_{\bfepsilon}(\nu),\sqrt{N}\left( \hat{\bfbeta}(\nu') - \bfbeta(\nu')\right)\right)\bfUpsilon^{\top}_M(\nu')\\
&=\bfV_{\nu\nu'}+\bfUpsilon_M(\nu)\bfTheta_{\nu\nu'}\bfUpsilon^{\top}_M(\nu')+\bfUpsilon_M(\nu)\bfG^{\mathrm{U}}_{\nu\nu'}+\left(\bfG^{\mathrm{U}}_{\nu\nu'}\right)^{\top}\bfUpsilon^{\top}_M(\nu')
\\&=:\bfDelta^{\mathrm{U}}_{\nu\nu'}.
\end{align*}
In view of \eqref{hat_rho}, we have
\begin{align*}
\lim_{N\rightarrow\infty}&\mathrm{Var}\left(\sqrt{N} \bfr_{\hat{\bfepsilon}}(\nu)\right)=\lim_{N\rightarrow\infty}\mathrm{Var}\left(\left(\diag\left[\left\{ \bfD^{\ast}_{{\bfepsilon}}(\nu-1),\dots,\bfD^{\ast}_{{\bfepsilon}}(\nu-M)\right\}  \otimes
                                    \bfD^{\ast}_{{\bfepsilon}}(\nu) \right]\right)^{-1}\sqrt{N} \bfc_{\hat{\bfepsilon}}(\nu)\right),\\&=
\left(\diag\left[\left\{ \bfD^{\ast}_{{\bfepsilon}}(\nu-1),\dots,\bfD^{\ast}_{{\bfepsilon}}(\nu-M)\right\}  \otimes
                                    \bfD^{\ast}_{{\bfepsilon}}(\nu) \right]\right)^{-1}\bfDelta^{\mathrm{U}}_{\nu\nu}
\left(\diag\left[\left\{ \bfD^{\ast}_{{\bfepsilon}}(\nu-1),\dots,\bfD^{\ast}_{{\bfepsilon}}(\nu-M)\right\}  \otimes
                                    \bfD^{\ast}_{{\bfepsilon}}(\nu) \right]\right)^{-1},
                                    \\&
=:\bfnabla^{\mathrm{U}}_{\nu\nu}, \qquad\text{for}\quad\nu=1,\dots,s,
\end{align*}
which prove \eqref{chatrhoU}. For $\nu\neq\nu'$ we also have:
\begin{align*}
\bfnabla^{\mathrm{U}}_{\nu\nu'} = \left(\diag\left[\left\{ \bfD^{\ast}_{{\bfepsilon}}(\nu-1),\dots,\bfD^{\ast}_{{\bfepsilon}}(\nu-M)\right\}  \otimes
                                    \bfD^{\ast}_{{\bfepsilon}}(\nu) \right]\right)^{-1}\bfDelta^{\mathrm{U}}_{\nu\nu'}
\left(\diag\left[\left\{ \bfD^{\ast}_{{\bfepsilon}}(\nu'-1),\dots,\bfD^{\ast}_{{\bfepsilon}}(\nu'-M)\right\}  \otimes
                                    \bfD^{\ast}_{{\bfepsilon}}(\nu') \right]\right)^{-1}.
\end{align*}
The asymptotic normality   of$\sqrt{N}\bfc_{\hat{\bfepsilon}}$ and $\sqrt{N}\bfr_{\hat{\bfepsilon}}$ follow using the same kind of arguments as \eqref{chatgammaU} and \eqref{chatrhoU} of $\sqrt{N}\bfc_{\hat{\bfepsilon}}(\nu)$ and $\sqrt{N}\bfr_{\hat{\bfepsilon}}(\nu)$. We then have:
\begin{align*}
  N^{1/2}\bfc_{\hat{\bfepsilon}} \cd \mathcal{N}_{sd^2M}\left(\bfzero, \bfDelta^{\mathrm{U}}\right), 
\end{align*}
where the asymptotic covariance matrix $\bfDelta^{\mathrm{U}}$ is a block matrix, with the asymptotic variances given by $\bfDelta^{\mathrm{U}}_{\nu\nu}$, for  $\nu=1,\dots,s$  and the asymptotic covariances given by
$$\bfDelta^{\mathrm{U}}_{\nu\nu'}=\bfV_{\nu\nu'}+\bfUpsilon_M(\nu)\bfTheta_{\nu\nu'}\bfUpsilon^{\top}_M(\nu')+\bfUpsilon_M(\nu)\bfG^{\mathrm{U}}_{\nu\nu'}+\left(\bfG^{\mathrm{U}}_{\nu\nu'}\right)^{\top}\bfUpsilon^{\top}_M(\nu'),\quad \nu\neq\nu'.$$
We also have
\begin{align*}
  N^{1/2}\bfr_{\hat{\bfepsilon}} \cd \mathcal{N}_{sd^2M}\left(\bfzero, \bfnabla^{\mathrm{U}}\right), 
\end{align*}
where the asymptotic covariance matrix $\bfnabla^{\mathrm{U}}= [\bfnabla^{\mathrm{U}}_{\nu\nu'}]_{1\leq \nu,\nu'\leq s}$ is a block matrix.

This ends our second step and the proof is completed.$\hfill\square$

\subsection{Proof of Corollary~\ref{bobo}}
We suppose that $H_1(\nu)$ holds true. In view of \eqref{hat_rho}, one may rewrite the above arguments in order to prove that there exists a nonsingular matrix $\bfL^{\mathrm{U\, or \,R}}_{\nu\nu'}$ such that 
\begin{equation}\label{tt}
\sqrt{N} \left( \bfr_{\hat{\bfepsilon}}(\nu) - \bfrho_{\bfepsilon}(\nu) \right) 
=
\bfK^{-1}(\nu)
\sqrt{N}\left(\bfc_{\hat{\bfepsilon}}(\nu)  - \bfgamma_{\bfepsilon}(\nu) \right)  + \mathrm{o}_{\mathbb P} (1)
 \xrightarrow[N\to\infty]{\mathrm{d}} \mathcal{N}_{d^2(p(\nu)+M)}\left(\bfzero,\bfL^{\mathrm{U\, or \,R}}_{\nu\nu'} \right) \ ,
\end{equation}
with $\bfK(\nu)= \diag\left[\left\{ \bfD^\ast_{{\bfepsilon}}(\nu-1),\dots,\bfD^\ast_{{\bfepsilon}}(\nu-M)\right\}  \otimes
                                    \bfD^\ast_{{\bfepsilon}}(\nu) \right]$ and the matrix $\bfL^{\mathrm{U\, or \,R}}_{\nu\nu'}$ is given by $$\bfL^{\mathrm{U\, or \,R}}_{\nu\nu'}=\bfK^{-1}(\nu)\left[ \bfV^\ast_{\nu\nu}+\bfUpsilon^\ast_M(\nu)\bfTheta_{\nu\nu}\bfUpsilon^{\ast\,\top}_M(\nu)+ \bfUpsilon^\ast_M(\nu)\bfG^{\ast\,\mathrm{U\, or \,R}}_{\nu\nu}+\left(\bfG^{\mathrm{\ast\,U\, or \,R}}_{\nu\nu}\right)^{\top}\bfUpsilon^{\ast\,\top}_M(\nu)\right]\bfK^{-1}(\nu),$$ where the matrices $\bfV^\ast_{\nu\nu}$ and  $\bfG^{\ast\,\mathrm{U\, or \,R}}_{\nu\nu}$ are obtained from the asymptotic distribution of 
\begin{eqnarray}\nonumber
  N^{1/2} \left( \hat{\bfbeta}^\top(\nu) - \bfbeta^\top(\nu),\bfc_{\bfepsilon}^\top(\nu)  - \bfgamma^\top_{\bfepsilon}(\nu)\right)^\top
  & \xrightarrow[N\to\infty]{\mathrm{d}}& \mathcal{N}_{d^2(p(\nu)+M)}\left(\bfzero, \begin{pmatrix} \bfTheta_{\nu\nu} & \bfG^{\ast\,\mathrm{U\, or \,R}}_{\nu\nu} \\ \left(\bfG^{\ast\,\mathrm{U\, or \,R}}_{\nu\nu}\right)^\top & \bfV^\ast_{\nu\nu} \end{pmatrix}\right),
\end{eqnarray}
The matrix $\bfUpsilon^{\ast}_M(\nu)$ is given by 
\begin{align*}
 \bfUpsilon^{\ast}_M(\nu) &:= \mathbb{E}\left[\left(\begin{array}{c}
\bfepsilon_{ns+\nu-1}\\\vdots\\\bfepsilon_{ns+\nu-M}\end{array}\right) \otimes \dfrac{\partial \bfepsilon_{ns+\nu} }{\partial \bfbeta^{\top}(\nu)}+ 
\left(\begin{array}{c}
\dfrac{\partial \bfepsilon_{ns+\nu-1}}{\partial \bfbeta^{\top}(\nu)}\\\vdots\\\dfrac{\partial \bfepsilon_{ns+\nu-M}}{\partial \bfbeta^{\top}(\nu)}\end{array}\right)
\otimes \bfepsilon_{ns+\nu}\right]
\quad\text{where}\quad\dfrac{\partial \bfepsilon_{ns+\nu} }{\partial \bfbeta^{\top}(\nu)}=-\bfX_n^\top(\nu)\otimes \bfI_d.
\end{align*}
We point out the fact that under $H_1(\nu)$, $$ \mathbb{E}\left[ \dfrac{\partial \bfepsilon_{ns+\nu-h}}{\partial \bfbeta^{\top}(\nu)}\otimes \bfepsilon_{ns+\nu}\right] \neq \bfzero,\quad\text{for}\quad h=1,\ldots,M$$ whereas it vanishes under $H_0(\nu)$. 
Thus we have
\begin{equation*}
  \frac{\partial \bfc_{{\bfepsilon}}(\nu)}{\partial \bfbeta^{\top}(\nu)} \cp \bfUpsilon^\ast_M(\nu)=-\mathbb{E}\left\{\left(\begin{array}{c}
\bfepsilon_{ns+\nu-1}\\\vdots\\\bfepsilon_{ns+\nu-M}\end{array}\right)\otimes
\bfX_n^\top(\nu)\otimes \bfI_d+ 
\left(\begin{array}{c}
\bfX_n^\top(\nu-1)\otimes \bfI_d\\\vdots\\
\bfX_n^\top(\nu-M)\otimes \bfI_d\end{array}\right)
\otimes \bfepsilon_{ns+\nu}\right\}.
\end{equation*}
Now we write 
\begin{align*}\hat{\bfN}^{-1/2}(\nu)\sqrt{N} \bfr_{\hat{\bfepsilon}}(\nu)    & =  \hat{\bfN}^{-1/2}(\nu)\sqrt{N} \left( \bfr_{\hat{\bfepsilon}}(\nu) - \bfrho_{\bfepsilon}(\nu) \right) +  \hat{\bfN}^{-1/2}(\nu)\sqrt{N} \bfrho_{\bfepsilon}(\nu)  \\ 
& =  {\bfN}^{-1/2}(\nu)\bfK^{-1}(\nu) \left[
\sqrt{N}\left(\bfc_{\hat{\bfepsilon}}(\nu)  - \bfgamma_{\bfepsilon}(\nu) \right)   +  \sqrt{N}  \bfgamma_{\bfepsilon}(\nu) \right]  + \mathrm{o}_{\mathbb P} (1)\ . 
\end{align*}
where $\hat{\bfN}(\nu)=\left(\diag\left[\bfR_{\hat{\bfepsilon}}(0; \nu-1),\dots, \bfR_{\hat{\bfepsilon}}(0; \nu-M)\right]\right) \otimes \bfR_{\hat{\bfepsilon}}(0; \nu)$ and $\bfN(\nu)=\left(\diag\left[\bfR^\ast_{{\bfepsilon}}(0; \nu-1),\dots, \bfR^\ast_{{\bfepsilon}}(0; \nu-M)\right]\right) \otimes \bfR^\ast_{{\bfepsilon}}(0; \nu)$ with 
$\bfR^\ast_{\bfepsilon}(0; \nu)=\left(\bfD^\ast_{\bfepsilon}(\nu)\right)^{-1}\bfGamma_{\bfepsilon}(0; \nu)
                                    \left(\bfD^\ast_{\bfepsilon}(\nu)\right)^{-1}$.              
Then it holds that 
\begin{align}\label{opp}
\bfQ_M(\nu) &= N\bfr_{\hat{\bfepsilon}}^{\top}(\nu) \hat{\bfN}^{-1}(\nu)\bfr_{\hat{\bfepsilon}}(\nu)=\left(\hat{\bfN}^{-1/2}(\nu)\sqrt{N} \bfr_{\hat{\bfepsilon}}(\nu)\right)^\top\times \left(\hat{\bfN}^{-1/2}(\nu)\sqrt{N} \bfr_{\hat{\bfepsilon}}(\nu)\right)
\nonumber \\ 
 & =    N\left(\bfc_{\hat{\bfepsilon}}(\nu)  - \bfgamma_{\bfepsilon}(\nu) \right)^\top \bfK^{-1}(\nu){\bfN}^{-1}(\nu)\bfK^{-1}(\nu) \left(\bfc_{\hat{\bfepsilon}}(\nu)  - \bfgamma_{\bfepsilon}(\nu) \right) \nonumber \\&\qquad + 2 N \left(\bfc_{\hat{\bfepsilon}}(\nu)  - \bfgamma_{\bfepsilon}(\nu) \right)^\top \bfK^{-1}(\nu){\bfN}^{-1}(\nu)\bfK^{-1}(\nu)  \bfgamma_{\bfepsilon}(\nu)  + N \bfgamma^\top_{\bfepsilon}(\nu) \bfK^{-1}(\nu){\bfN}^{-1}(\nu)\bfK^{-1}(\nu) \bfgamma_{\bfepsilon}(\nu)  + \mathrm{o}_{\mathbb P} (1) 
\end{align}
By the ergodic theroem, $ \left(\bfc_{\hat{\bfepsilon}}(\nu)  - \bfgamma_{\bfepsilon}(\nu) \right)^\top \bfK^{-1}(\nu){\bfN}^{-1}(\nu)\bfK^{-1}(\nu)  \bfgamma_{\bfepsilon}(\nu) = \mathrm{o}_{\mathbb P} (1)$. Using  \citet[Lemma 17.1]{vdv}, the convergence \eqref{tt} implies that 
$$N\left(\bfc_{\hat{\bfepsilon}}(\nu)  - \bfgamma_{\bfepsilon}(\nu) \right)^\top \bfK^{-1}(\nu){\bfN}^{-1}(\nu)\bfK^{-1}(\nu) \left(\bfc_{\hat{\bfepsilon}}(\nu)  - \bfgamma_{\bfepsilon}(\nu) \right) \xrightarrow[N\to\infty]{\mathrm{d}} Z_M^{\nu}(\xi_{Md^2}^{\mathrm{U\, or \, R}}(\nu))=\sum_{i=1}^{Md^2}\xi_{i,Md^2}^{\mathrm{U\, or \, R}}(\nu) Z_i^2$$ 
 where $\xi_{Md^2}^{\mathrm{U\, or\, R}}(\nu)=(\xi_{1,Md^2}^{\mathrm{U\, or \, R}}(\nu),\dots,\xi_{Md^2,Md^2}^{\mathrm{U\, or \, R}}(\nu))^\top$ is the vector of the eigenvalues of the matrix ${\bfN}^{-1/2}(\nu)\bfL^{\mathrm{U\, or \,R}}_{\nu\nu'}{\bfN}^{-1/2}(\nu)$ and  $(Z_i)_{1\le i\le Md^2}$ are i.i.d. with $\mathcal N(0,1)$ laws. 
Reporting these convergences in \eqref{opp}, we deduce that 
\begin{align*}
\frac{\bfQ_M(\nu)}{N}
 & =   \left(\bfc_{\hat{\bfepsilon}}(\nu)  - \bfgamma_{\bfepsilon}(\nu) \right)^\top \bfK^{-1}(\nu){\bfN}^{-1}(\nu)\bfK^{-1}(\nu) \left(\bfc_{\hat{\bfepsilon}}(\nu)  - \bfgamma_{\bfepsilon}(\nu) \right)  \\&\qquad + 2 \left(\bfc_{\hat{\bfepsilon}}(\nu)  - \bfgamma_{\bfepsilon}(\nu) \right)^\top \bfK^{-1}(\nu){\bfN}^{-1}(\nu)\bfK^{-1}(\nu)  \bfgamma_{\bfepsilon}(\nu)  +  \bfgamma^\top_{\bfepsilon}(\nu) \bfK^{-1}(\nu){\bfN}^{-1}(\nu)\bfK^{-1}(\nu) \bfgamma_{\bfepsilon}(\nu)  + \mathrm{o}_{\mathbb P} \left(\frac{1}{N}\right)  \\ 
 & \xrightarrow[N\to\infty]{\mathbb{P}} \bfgamma^\top_{\bfepsilon}(\nu) \bfK^{-1}(\nu){\bfN}^{-1}(\nu)\bfK^{-1}(\nu) \bfgamma_{\bfepsilon}(\nu)=\bfrho^\top_{\bfepsilon}(\nu) {\bfN}^{-1}(\nu) \bfrho_{\bfepsilon}(\nu).
\end{align*}
The same calculations hold for $\bfQ^\ast_M(\nu)$ and the corollary is proved. $\hfill\square$

\noindent
{\bf Acknowledgements:}
We sincerely thank the
anonymous reviewers and editor for helpful remarks.




\end{document}